\documentclass[10pt]{amsart}
\usepackage[english,russian]{babel}
\usepackage[cp1251]{inputenc}
\usepackage{amsmath}
\usepackage{amssymb}
\usepackage{amsfonts}

\usepackage[linktocpage=true, colorlinks=true, linkcolor=blue, citecolor=blue, urlcolor=blue]{hyperref}

\begin{document}
\renewcommand{\refname}{References}
\renewcommand\contentsname{Contents}

\thispagestyle{empty}

\title[Exact Calculation of the Mean-Square Error]
{Exact Calculation of the Mean-Square Error in the
Method of Expansion
of Iterated Ito Stochastic integrals Based 
on Generalized Multiple Fourier Series}
\author[D.F. Kuznetsov]{Dmitriy F. Kuznetsov}
\address{Dmitriy Feliksovich Kuznetsov
\newline\hphantom{iii} Peter the Great Saint-Petersburg Polytechnic University,
\newline\hphantom{iii} Polytechnicheskaya ul., 29,
\newline\hphantom{iii} 195251, Saint-Petersburg, Russia}%
\email{sde\_kuznetsov@inbox.ru}

\maketitle {\small
\begin{quote}
\vspace{5mm}
\noindent{\sc Abstract.} 
The article is devoted to the developement of the
method of expansion and mean-square approximation 
of iterated Ito stochastic integrals based on generalized
multiple Fourier series converging in the sense of norm
in the space $L_2([t, T]^k)$ ($k$ is the multiplicity
of the iterated Ito stochastic integral).
We obtain the
exact and approximate 
expressions for the mean-square
approximation error of iterated Ito stochastic integrals of
multiplicity $k$ ($k\in\mathbb{N}$) from
the stochastic Taylor--Ito expansion. 
As a result, we do not need to use redundant terms
of expansions of iterated Ito stochastic integrals
that complicate the numerical methods
for Ito stochastic differential equations.
Moreover, we proved the convergence with propability 1 
for the method of expansion 
of iterated Ito stochastic integrals based on generalized
multiple Fourier series for the cases of multiple Fourier--Legendre series
and multiple trigonometric Fourier series.
Mean-square approximation of iterated Stratonovich
stochastic integrals is also considered in the article.
The results of the article can be applied to the 
high-order strong
numerical methods for Ito stochastic differential equations
as well as for non-commutative semilinear
stochastic partial differential equations with multiplicative 
trace class noise.

\medskip
\end{quote}
}

\vspace{6mm}

\linespread{1.6}

\tableofcontents

\linespread{1.0}

\vspace{3mm}

\section{Introduction}

\vspace{5mm}

In this article we develop the method of expansion
and mean-square approximation of iterated Ito stochastic integrals
based on generalized multiple Fourier series converging
in the sense of norm in the space 
$L_2([t, T]^k)$ ($k$ is the multiplicity
of the iterated Ito stochastic integral),
which was proposed and developed by the author of this work
\cite{1}-\cite{new-2023a} (also see related publications
\cite{new-art-1}-\cite{old-art-4}).
Further, this method is referred to as the method of generalized 
multiple Fourier series.

The question of how to estimate or calculate 
exactly the mean-square approximation error of iterated Ito stochastic 
integrals
for the method of generalized 
multiple Fourier series composes the subject of the article.
From the one side the mentioned question is essentially difficult
for the case of a multidimensional Wiener process, because of we need to take
into account all possible combinations of components of the
multidimensional Wiener
process. From the other side an effective solution of the 
mentioned problem allows us to construct
more economical numerical methods for Ito stochastic differential equations
than in \cite{Zapad-1}-\cite{Zapad-3}.

The results of the article (also see \cite{1}-\cite{new-2023a} and 
related publications \cite{new-art-1}-\cite{old-art-4})
will be useful for the
implementation
of high-order strong numerical methods for 
Ito stochastic differential equations
as well as for non-commutative semilinear
stochastic partial differential equations with multiplicative trace class 
noise. The latter methods are constucted, for example, in 
\cite{Zapad-4}, \cite{Zapad-5}.

\vspace{5mm}

\section{Method of Generalized Multiple Fourier Series. The Case 
of Complete Orthonormal Systems  
of Continuous Functions in the Space $L_2([t,T])$ 
and Continuous Weight Functions $\psi_1(\tau),\ldots,\psi_k(\tau)$
}

\vspace{5mm}

Let $(\Omega,$ ${\rm F},$ ${\sf P})$ be a complete probability space, let 
$\{{\rm F}_t, t\in[0,T]\}$ be a nondecreasing right-continous family of 
$\sigma$-subfields of ${\rm F},$
and let ${\bf f}_t$ be a standard $m$-dimensional Wiener 
stochastic process, which is
${\rm F}_t$-measurable for any $t\in[0, T].$ We assume that the components
${\bf f}_{t}^{(i)}$ $(i=1,\ldots,m)$ of this process are independent.

Let us consider 
the following iterated Ito 
stochastic integrals

\vspace{-2mm}
\begin{equation}
\label{sodom20}
J[\psi^{(k)}]_{T,t}=\int\limits_t^T\psi_k(t_k) \ldots \int\limits_t^{t_{2}}
\psi_1(t_1) d{\bf w}_{t_1}^{(i_1)}\ldots
d{\bf w}_{t_k}^{(i_k)},
\end{equation}

\vspace{2mm}
\noindent
where $\psi_l(\tau)$ $(l=1,\ldots,k)$ are
nonrandom functions at the interval $[t, T]$,
${\bf w}_{\tau}^{(i)}={\bf f}_{\tau}^{(i)}$
for $i=1,\ldots,m$ and
${\bf w}_{\tau}^{(0)}=\tau,$ 
$i_1,\ldots,i_k=0, 1,\ldots,m.$

In addition, suppose that 
every $\psi_l(\tau)$ $(l=1,\ldots,k)$ is a continuous 
nonrandom function at the interval $[t, T]$
(the case $\psi_1(\tau),\ldots,\psi_k(\tau)\in L_2([t, T])$
will be considered in Sect.~5).

Let us define the following function on the hypercube $[t, T]^k$

\vspace{1mm}
\begin{equation}
\label{ppp}
K(t_1,\ldots,t_k)=
\begin{cases}
\psi_1(t_1)\ldots \psi_k(t_k),\ &t_1<\ldots<t_k\\
~\\
~\\
0,\ &\hbox{\rm otherwise}
\end{cases}
\ \ \ =\ \ \ \ 
\prod\limits_{l=1}^k
\psi_l(t_l)\  \prod\limits_{l=1}^{k-1}\  {\bf 1}_{\{t_l<t_{l+1}\}},\ 
\end{equation}

\vspace{5mm}
\noindent
where $t_1,\ldots,t_k\in [t, T]$ for $k\ge 2$ and 
$K(t_1)\equiv\psi_1(t_1)$ for $t_1\in[t, T].$ Here 
${\bf 1}_A$ denotes the indicator of the set $A$.

Suppose that $\{\phi_j(x)\}_{j=0}^{\infty}$
is a complete orthonormal system of functions in 
the space $L_2([t, T])$. 
The function $K(t_1,\ldots,t_k)$ is piecewise continuous in the 
hypercube $[t, T]^k.$
At this situation it is well known that the generalized 
multiple Fourier series 
of $K(t_1,\ldots,t_k)\in L_2([t, T]^k)$ is converging 
to $K(t_1,\ldots,t_k)$ in the hypercube $[t, T]^k$ in 
the mean-square sense, i.e.

\vspace{1mm}
\begin{equation}
\label{sos1z}
\hbox{\vtop{\offinterlineskip\halign{
\hfil#\hfil\cr
{\rm lim}\cr
$\stackrel{}{{}_{p_1,\ldots,p_k\to \infty}}$\cr
}} }\Biggl\Vert
K(t_1,\ldots,t_k)-
\sum_{j_1=0}^{p_1}\ldots \sum_{j_k=0}^{p_k}
C_{j_k\ldots j_1}\prod_{l=1}^{k} \phi_{j_l}(t_l)\Biggr\Vert_{L_2([t,T]^k)}=0,
\end{equation}

\vspace{4mm}
\noindent
where

\vspace{-2mm}
\begin{equation}
\label{ppppa}
C_{j_k\ldots j_1}=\int\limits_{[t,T]^k}
K(t_1,\ldots,t_k)\prod_{l=1}^{k}\phi_{j_l}(t_l)dt_1\ldots dt_k
\end{equation}

\vspace{3mm}
\noindent
is the Fourier coefficient and

$$
\left\Vert f\right\Vert_{L_2([t,T]^k)}=\left(\int\limits_{[t,T]^k}
f^2(t_1,\ldots,t_k)dt_1\ldots dt_k\right)^{1/2}
$$

\vspace{5mm}
\noindent
is a norm in the space $L_2([t,T]^k)$.

Consider the partition $\{\tau_j\}_{j=0}^N$ of the interval $[t,T]$ such that

\begin{equation}
\label{1111}
t=\tau_0<\ldots <\tau_N=T,\ \ \
\Delta_N=
\hbox{\vtop{\offinterlineskip\halign{
\hfil#\hfil\cr
{\rm max}\cr
$\stackrel{}{{}_{0\le j\le N-1}}$\cr
}} }\Delta\tau_j\to 0\ \ \hbox{if}\ \ N\to \infty,\ \ \ 
\Delta\tau_j=\tau_{j+1}-\tau_j.
\end{equation}

\vspace{3mm}

{\bf Theorem 1} \cite{1} (2006) (also see
\cite{2}-\cite{new-2023a}). 
{\it Suppose that
every $\psi_l(\tau)$ $(l=1,\ldots, k)$ is a continuous nonrandom 
function on 
$[t, T]$ and
$\{\phi_j(x)\}_{j=0}^{\infty}$ is a complete orthonormal system  
of continuous functions in the space $L_2([t,T]).$ 
Then

$$
J[\psi^{(k)}]_{T,t}\  =\ 
\hbox{\vtop{\offinterlineskip\halign{
\hfil#\hfil\cr
{\rm l.i.m.}\cr
$\stackrel{}{{}_{p_1,\ldots,p_k\to \infty}}$\cr
}} }\sum_{j_1=0}^{p_1}\ldots\sum_{j_k=0}^{p_k}
C_{j_k\ldots j_1}\Biggl(
\prod_{l=1}^k\zeta_{j_l}^{(i_l)}\ -
\Biggr.
$$

\vspace{2mm}
\begin{equation}
\label{tyyyxxx}
-\ \Biggl.
\hbox{\vtop{\offinterlineskip\halign{
\hfil#\hfil\cr
{\rm l.i.m.}\cr
$\stackrel{}{{}_{N\to \infty}}$\cr
}} }\sum_{(l_1,\ldots,l_k)\in {\rm G}_k}
\phi_{j_{1}}(\tau_{l_1})
\Delta{\bf w}_{\tau_{l_1}}^{(i_1)}\ldots
\phi_{j_{k}}(\tau_{l_k})
\Delta{\bf w}_{\tau_{l_k}}^{(i_k)}\Biggr),
\end{equation}

\vspace{5mm}
\noindent
where

$$
{\rm G}_k={\rm H}_k\backslash{\rm L}_k,\ \ \
{\rm H}_k=\{(l_1,\ldots,l_k):\ l_1,\ldots,l_k=0,\ 1,\ldots,N-1\},
$$

$$
{\rm L}_k=\{(l_1,\ldots,l_k):\ l_1,\ldots,l_k=0,\ 1,\ldots,N-1;\
l_g\ne l_r\ (g\ne r);\ g, r=1,\ldots,k\},
$$

\vspace{5mm}
\noindent
${\rm l.i.m.}$ is a limit in the mean-square sense,
$i_1,\ldots,i_k=0,1,\ldots,m,$ 

\vspace{-2mm}
\begin{equation}
\label{rr23}
\zeta_{j}^{(i)}=
\int\limits_t^T \phi_{j}(s) d{\bf w}_s^{(i)}
\end{equation} 

\vspace{2mm}
\noindent
are independent standard Gaussian random variables
for various
$i$ or $j$ {\rm(}if $i\ne 0${\rm),}
$C_{j_k\ldots j_1}$ is the Fourier coefficient {\rm(\ref{ppppa}),}
$\Delta{\bf w}_{\tau_{j}}^{(i)}=
{\bf w}_{\tau_{j+1}}^{(i)}-{\bf w}_{\tau_{j}}^{(i)}$
$(i=0, 1,\ldots,m),$
$\left\{\tau_{j}\right\}_{j=0}^{N}$ is a partition of
$[t,T],$ which satisfies the condition {\rm (\ref{1111})}.}

\vspace{2mm}

{\bf Remark 1.} {\it Further {\rm (see Theorem {\rm 2})} we will use the 
following form of expansion {\rm (\ref{tyyyxxx})}

\begin{equation}
\label{tyyy}
J[\psi^{(k)}]_{T,t}=
\hbox{\vtop{\offinterlineskip\halign{
\hfil#\hfil\cr
{\rm l.i.m.}\cr
$\stackrel{}{{}_{p_1,\ldots,p_k\to \infty}}$\cr
}} }\sum_{j_1=0}^{p_1}\ldots\sum_{j_k=0}^{p_k}
C_{j_k\ldots j_1}\left(
\prod_{l=1}^k\zeta_{j_l}^{(i_l)}-S_{j_1,\ldots,j_k}^{(i_1\ldots i_k)}\right),
\end{equation}

\vspace{3mm}
\noindent
where

\vspace{-1mm}
\begin{equation}
\label{tyyy1}
S_{j_1,\ldots,j_k}^{(i_1\ldots i_k)}\ =\ 
\hbox{\vtop{\offinterlineskip\halign{
\hfil#\hfil\cr
{\rm l.i.m.}\cr
$\stackrel{}{{}_{N\to \infty}}$\cr
}} }\sum_{(l_1,\ldots,l_k)\in {\rm G}_k}
\phi_{j_{1}}(\tau_{l_1})
\Delta{\bf w}_{\tau_{l_1}}^{(i_1)}\ldots
\phi_{j_{k}}(\tau_{l_k})
\Delta{\bf w}_{\tau_{l_k}}^{(i_k)},
\end{equation}

\vspace{4mm}
\noindent
where notations are the same as in Theorem {\rm 1}.}

\vspace{2mm}

Note that the version of Theorem 1 for the Haar and
Rademacher--Walsh functions has been considered in 
\cite{1}-\cite{12}, \cite{arxiv-1}.
Some modifications of Theorem 1 for another types of iterated
stochastic integrals (including iterated stochastic integrals
with respect to the infinite-dimensional $Q$-Wiener process)
as well as for complete orthonormal with
weight $r(t_1)\ldots r(t_k)\ge 0$ systems of functions 
in the space $L_2([t, T]^k)$
can be found in \cite{arxiv-1}, \cite{arxiv-13}
\cite{arxiv-20}, \cite{arxiv-21} (also see \cite{1}-\cite{art-9},
\cite{arxiv-2}-\cite{arxiv-12},
\cite{arxiv-14}-\cite{arxiv-19}, \cite{arxiv-22}-\cite{new-2023a}).
Generalization of Theorem 1 for the case
of an arbitrary complete orthonormal systems  
of functions in the space $L_2([t,T])$ 
and $\psi_1(\tau),\ldots,\psi_k(\tau)\in L_2([t, T])$
will be considered in Sect.~5.

In order to evaluate the significance of Theorem 1 for practice we will
demonstrate its trans\-for\-med particular cases for 
$k=1,\ldots,6$ \cite{1}-\cite{new-2023a}

\begin{equation}
\label{a1}
J[\psi^{(1)}]_{T,t}
=\hbox{\vtop{\offinterlineskip\halign{
\hfil#\hfil\cr
{\rm l.i.m.}\cr
$\stackrel{}{{}_{p_1\to \infty}}$\cr
}} }\sum_{j_1=0}^{p_1}
C_{j_1}\zeta_{j_1}^{(i_1)},
\end{equation}

\vspace{2mm}
\begin{equation}
\label{leto5001}
J[\psi^{(2)}]_{T,t}
=\hbox{\vtop{\offinterlineskip\halign{
\hfil#\hfil\cr
{\rm l.i.m.}\cr
$\stackrel{}{{}_{p_1,p_2\to \infty}}$\cr
}} }\sum_{j_1=0}^{p_1}\sum_{j_2=0}^{p_2}
C_{j_2j_1}\Biggl(\zeta_{j_1}^{(i_1)}\zeta_{j_2}^{(i_2)}
-{\bf 1}_{\{i_1=i_2\ne 0\}}
{\bf 1}_{\{j_1=j_2\}}\Biggr),
\end{equation}

\vspace{5mm}
$$
J[\psi^{(3)}]_{T,t}=
\hbox{\vtop{\offinterlineskip\halign{
\hfil#\hfil\cr
{\rm l.i.m.}\cr
$\stackrel{}{{}_{p_1,\ldots,p_3\to \infty}}$\cr
}} }\sum_{j_1=0}^{p_1}\sum_{j_2=0}^{p_2}\sum_{j_3=0}^{p_3}
C_{j_3j_2j_1}\Biggl(
\zeta_{j_1}^{(i_1)}\zeta_{j_2}^{(i_2)}\zeta_{j_3}^{(i_3)}
-\Biggr.
$$
\begin{equation}
\label{leto5002}
\Biggl.-{\bf 1}_{\{i_1=i_2\ne 0\}}
{\bf 1}_{\{j_1=j_2\}}
\zeta_{j_3}^{(i_3)}
-{\bf 1}_{\{i_2=i_3\ne 0\}}
{\bf 1}_{\{j_2=j_3\}}
\zeta_{j_1}^{(i_1)}-
{\bf 1}_{\{i_1=i_3\ne 0\}}
{\bf 1}_{\{j_1=j_3\}}
\zeta_{j_2}^{(i_2)}\Biggr),
\end{equation}

\vspace{5mm}
$$
J[\psi^{(4)}]_{T,t}
=
\hbox{\vtop{\offinterlineskip\halign{
\hfil#\hfil\cr
{\rm l.i.m.}\cr
$\stackrel{}{{}_{p_1,\ldots,p_4\to \infty}}$\cr
}} }\sum_{j_1=0}^{p_1}\ldots\sum_{j_4=0}^{p_4}
C_{j_4\ldots j_1}\Biggl(
\prod_{l=1}^4\zeta_{j_l}^{(i_l)}
\Biggr.
-
$$
$$
-
{\bf 1}_{\{i_1=i_2\ne 0\}}
{\bf 1}_{\{j_1=j_2\}}
\zeta_{j_3}^{(i_3)}
\zeta_{j_4}^{(i_4)}
-
{\bf 1}_{\{i_1=i_3\ne 0\}}
{\bf 1}_{\{j_1=j_3\}}
\zeta_{j_2}^{(i_2)}
\zeta_{j_4}^{(i_4)}-
$$
$$
-
{\bf 1}_{\{i_1=i_4\ne 0\}}
{\bf 1}_{\{j_1=j_4\}}
\zeta_{j_2}^{(i_2)}
\zeta_{j_3}^{(i_3)}
-
{\bf 1}_{\{i_2=i_3\ne 0\}}
{\bf 1}_{\{j_2=j_3\}}
\zeta_{j_1}^{(i_1)}
\zeta_{j_4}^{(i_4)}-
$$
$$
-
{\bf 1}_{\{i_2=i_4\ne 0\}}
{\bf 1}_{\{j_2=j_4\}}
\zeta_{j_1}^{(i_1)}
\zeta_{j_3}^{(i_3)}
-
{\bf 1}_{\{i_3=i_4\ne 0\}}
{\bf 1}_{\{j_3=j_4\}}
\zeta_{j_1}^{(i_1)}
\zeta_{j_2}^{(i_2)}+
$$
$$
+
{\bf 1}_{\{i_1=i_2\ne 0\}}
{\bf 1}_{\{j_1=j_2\}}
{\bf 1}_{\{i_3=i_4\ne 0\}}
{\bf 1}_{\{j_3=j_4\}}
+
{\bf 1}_{\{i_1=i_3\ne 0\}}
{\bf 1}_{\{j_1=j_3\}}
{\bf 1}_{\{i_2=i_4\ne 0\}}
{\bf 1}_{\{j_2=j_4\}}+
$$
\begin{equation}
\label{leto5003}
+\Biggl.
{\bf 1}_{\{i_1=i_4\ne 0\}}
{\bf 1}_{\{j_1=j_4\}}
{\bf 1}_{\{i_2=i_3\ne 0\}}
{\bf 1}_{\{j_2=j_3\}}\Biggr),
\end{equation}

\vspace{5mm}
$$
J[\psi^{(5)}]_{T,t}
=\hbox{\vtop{\offinterlineskip\halign{
\hfil#\hfil\cr
{\rm l.i.m.}\cr
$\stackrel{}{{}_{p_1,\ldots,p_5\to \infty}}$\cr
}} }\sum_{j_1=0}^{p_1}\ldots\sum_{j_5=0}^{p_5}
C_{j_5\ldots j_1}\Biggl(
\prod_{l=1}^5\zeta_{j_l}^{(i_l)}
-\Biggr.
$$
$$
-
{\bf 1}_{\{i_1=i_2\ne 0\}}
{\bf 1}_{\{j_1=j_2\}}
\zeta_{j_3}^{(i_3)}
\zeta_{j_4}^{(i_4)}
\zeta_{j_5}^{(i_5)}-
{\bf 1}_{\{i_1=i_3\ne 0\}}
{\bf 1}_{\{j_1=j_3\}}
\zeta_{j_2}^{(i_2)}
\zeta_{j_4}^{(i_4)}
\zeta_{j_5}^{(i_5)}-
$$
$$
-
{\bf 1}_{\{i_1=i_4\ne 0\}}
{\bf 1}_{\{j_1=j_4\}}
\zeta_{j_2}^{(i_2)}
\zeta_{j_3}^{(i_3)}
\zeta_{j_5}^{(i_5)}-
{\bf 1}_{\{i_1=i_5\ne 0\}}
{\bf 1}_{\{j_1=j_5\}}
\zeta_{j_2}^{(i_2)}
\zeta_{j_3}^{(i_3)}
\zeta_{j_4}^{(i_4)}-
$$
$$
-
{\bf 1}_{\{i_2=i_3\ne 0\}}
{\bf 1}_{\{j_2=j_3\}}
\zeta_{j_1}^{(i_1)}
\zeta_{j_4}^{(i_4)}
\zeta_{j_5}^{(i_5)}-
{\bf 1}_{\{i_2=i_4\ne 0\}}
{\bf 1}_{\{j_2=j_4\}}
\zeta_{j_1}^{(i_1)}
\zeta_{j_3}^{(i_3)}
\zeta_{j_5}^{(i_5)}-
$$
$$
-
{\bf 1}_{\{i_2=i_5\ne 0\}}
{\bf 1}_{\{j_2=j_5\}}
\zeta_{j_1}^{(i_1)}
\zeta_{j_3}^{(i_3)}
\zeta_{j_4}^{(i_4)}
-{\bf 1}_{\{i_3=i_4\ne 0\}}
{\bf 1}_{\{j_3=j_4\}}
\zeta_{j_1}^{(i_1)}
\zeta_{j_2}^{(i_2)}
\zeta_{j_5}^{(i_5)}-
$$
$$
-
{\bf 1}_{\{i_3=i_5\ne 0\}}
{\bf 1}_{\{j_3=j_5\}}
\zeta_{j_1}^{(i_1)}
\zeta_{j_2}^{(i_2)}
\zeta_{j_4}^{(i_4)}
-{\bf 1}_{\{i_4=i_5\ne 0\}}
{\bf 1}_{\{j_4=j_5\}}
\zeta_{j_1}^{(i_1)}
\zeta_{j_2}^{(i_2)}
\zeta_{j_3}^{(i_3)}+
$$
$$
+
{\bf 1}_{\{i_1=i_2\ne 0\}}
{\bf 1}_{\{j_1=j_2\}}
{\bf 1}_{\{i_3=i_4\ne 0\}}
{\bf 1}_{\{j_3=j_4\}}\zeta_{j_5}^{(i_5)}+
{\bf 1}_{\{i_1=i_2\ne 0\}}
{\bf 1}_{\{j_1=j_2\}}
{\bf 1}_{\{i_3=i_5\ne 0\}}
{\bf 1}_{\{j_3=j_5\}}\zeta_{j_4}^{(i_4)}+
$$
$$
+
{\bf 1}_{\{i_1=i_2\ne 0\}}
{\bf 1}_{\{j_1=j_2\}}
{\bf 1}_{\{i_4=i_5\ne 0\}}
{\bf 1}_{\{j_4=j_5\}}\zeta_{j_3}^{(i_3)}+
{\bf 1}_{\{i_1=i_3\ne 0\}}
{\bf 1}_{\{j_1=j_3\}}
{\bf 1}_{\{i_2=i_4\ne 0\}}
{\bf 1}_{\{j_2=j_4\}}\zeta_{j_5}^{(i_5)}+
$$
$$
+
{\bf 1}_{\{i_1=i_3\ne 0\}}
{\bf 1}_{\{j_1=j_3\}}
{\bf 1}_{\{i_2=i_5\ne 0\}}
{\bf 1}_{\{j_2=j_5\}}\zeta_{j_4}^{(i_4)}+
{\bf 1}_{\{i_1=i_3\ne 0\}}
{\bf 1}_{\{j_1=j_3\}}
{\bf 1}_{\{i_4=i_5\ne 0\}}
{\bf 1}_{\{j_4=j_5\}}\zeta_{j_2}^{(i_2)}+
$$
$$
+
{\bf 1}_{\{i_1=i_4\ne 0\}}
{\bf 1}_{\{j_1=j_4\}}
{\bf 1}_{\{i_2=i_3\ne 0\}}
{\bf 1}_{\{j_2=j_3\}}\zeta_{j_5}^{(i_5)}+
{\bf 1}_{\{i_1=i_4\ne 0\}}
{\bf 1}_{\{j_1=j_4\}}
{\bf 1}_{\{i_2=i_5\ne 0\}}
{\bf 1}_{\{j_2=j_5\}}\zeta_{j_3}^{(i_3)}+
$$
$$
+
{\bf 1}_{\{i_1=i_4\ne 0\}}
{\bf 1}_{\{j_1=j_4\}}
{\bf 1}_{\{i_3=i_5\ne 0\}}
{\bf 1}_{\{j_3=j_5\}}\zeta_{j_2}^{(i_2)}+
{\bf 1}_{\{i_1=i_5\ne 0\}}
{\bf 1}_{\{j_1=j_5\}}
{\bf 1}_{\{i_2=i_3\ne 0\}}
{\bf 1}_{\{j_2=j_3\}}\zeta_{j_4}^{(i_4)}+
$$
$$
+
{\bf 1}_{\{i_1=i_5\ne 0\}}
{\bf 1}_{\{j_1=j_5\}}
{\bf 1}_{\{i_2=i_4\ne 0\}}
{\bf 1}_{\{j_2=j_4\}}\zeta_{j_3}^{(i_3)}+
{\bf 1}_{\{i_1=i_5\ne 0\}}
{\bf 1}_{\{j_1=j_5\}}
{\bf 1}_{\{i_3=i_4\ne 0\}}
{\bf 1}_{\{j_3=j_4\}}\zeta_{j_2}^{(i_2)}+
$$
$$
+
{\bf 1}_{\{i_2=i_3\ne 0\}}
{\bf 1}_{\{j_2=j_3\}}
{\bf 1}_{\{i_4=i_5\ne 0\}}
{\bf 1}_{\{j_4=j_5\}}\zeta_{j_1}^{(i_1)}+
{\bf 1}_{\{i_2=i_4\ne 0\}}
{\bf 1}_{\{j_2=j_4\}}
{\bf 1}_{\{i_3=i_5\ne 0\}}
{\bf 1}_{\{j_3=j_5\}}\zeta_{j_1}^{(i_1)}+
$$
\begin{equation}
\label{a5}
+\Biggl.
{\bf 1}_{\{i_2=i_5\ne 0\}}
{\bf 1}_{\{j_2=j_5\}}
{\bf 1}_{\{i_3=i_4\ne 0\}}
{\bf 1}_{\{j_3=j_4\}}\zeta_{j_1}^{(i_1)}\Biggr),
\end{equation}

\vspace{5mm}

$$
J[\psi^{(6)}]_{T,t}
=\hbox{\vtop{\offinterlineskip\halign{
\hfil#\hfil\cr
{\rm l.i.m.}\cr
$\stackrel{}{{}_{p_1,\ldots,p_6\to \infty}}$\cr
}} }\sum_{j_1=0}^{p_1}\ldots\sum_{j_6=0}^{p_6}
C_{j_6\ldots j_1}\Biggl(
\prod_{l=1}^6
\zeta_{j_l}^{(i_l)}
-\Biggr.
$$
$$
-
{\bf 1}_{\{i_1=i_6\ne 0\}}
{\bf 1}_{\{j_1=j_6\}}
\zeta_{j_2}^{(i_2)}
\zeta_{j_3}^{(i_3)}
\zeta_{j_4}^{(i_4)}
\zeta_{j_5}^{(i_5)}-
{\bf 1}_{\{i_2=i_6\ne 0\}}
{\bf 1}_{\{j_2=j_6\}}
\zeta_{j_1}^{(i_1)}
\zeta_{j_3}^{(i_3)}
\zeta_{j_4}^{(i_4)}
\zeta_{j_5}^{(i_5)}-
$$
$$
-
{\bf 1}_{\{i_3=i_6\ne 0\}}
{\bf 1}_{\{j_3=j_6\}}
\zeta_{j_1}^{(i_1)}
\zeta_{j_2}^{(i_2)}
\zeta_{j_4}^{(i_4)}
\zeta_{j_5}^{(i_5)}-
{\bf 1}_{\{i_4=i_6\ne 0\}}
{\bf 1}_{\{j_4=j_6\}}
\zeta_{j_1}^{(i_1)}
\zeta_{j_2}^{(i_2)}
\zeta_{j_3}^{(i_3)}
\zeta_{j_5}^{(i_5)}-
$$
$$
-
{\bf 1}_{\{i_5=i_6\ne 0\}}
{\bf 1}_{\{j_5=j_6\}}
\zeta_{j_1}^{(i_1)}
\zeta_{j_2}^{(i_2)}
\zeta_{j_3}^{(i_3)}
\zeta_{j_4}^{(i_4)}-
{\bf 1}_{\{i_1=i_2\ne 0\}}
{\bf 1}_{\{j_1=j_2\}}
\zeta_{j_3}^{(i_3)}
\zeta_{j_4}^{(i_4)}
\zeta_{j_5}^{(i_5)}
\zeta_{j_6}^{(i_6)}-
$$
$$
-
{\bf 1}_{\{i_1=i_3\ne 0\}}
{\bf 1}_{\{j_1=j_3\}}
\zeta_{j_2}^{(i_2)}
\zeta_{j_4}^{(i_4)}
\zeta_{j_5}^{(i_5)}
\zeta_{j_6}^{(i_6)}-
{\bf 1}_{\{i_1=i_4\ne 0\}}
{\bf 1}_{\{j_1=j_4\}}
\zeta_{j_2}^{(i_2)}
\zeta_{j_3}^{(i_3)}
\zeta_{j_5}^{(i_5)}
\zeta_{j_6}^{(i_6)}-
$$
$$
-
{\bf 1}_{\{i_1=i_5\ne 0\}}
{\bf 1}_{\{j_1=j_5\}}
\zeta_{j_2}^{(i_2)}
\zeta_{j_3}^{(i_3)}
\zeta_{j_4}^{(i_4)}
\zeta_{j_6}^{(i_6)}-
{\bf 1}_{\{i_2=i_3\ne 0\}}
{\bf 1}_{\{j_2=j_3\}}
\zeta_{j_1}^{(i_1)}
\zeta_{j_4}^{(i_4)}
\zeta_{j_5}^{(i_5)}
\zeta_{j_6}^{(i_6)}-
$$
$$
-
{\bf 1}_{\{i_2=i_4\ne 0\}}
{\bf 1}_{\{j_2=j_4\}}
\zeta_{j_1}^{(i_1)}
\zeta_{j_3}^{(i_3)}
\zeta_{j_5}^{(i_5)}
\zeta_{j_6}^{(i_6)}-
{\bf 1}_{\{i_2=i_5\ne 0\}}
{\bf 1}_{\{j_2=j_5\}}
\zeta_{j_1}^{(i_1)}
\zeta_{j_3}^{(i_3)}
\zeta_{j_4}^{(i_4)}
\zeta_{j_6}^{(i_6)}-
$$
$$
-
{\bf 1}_{\{i_3=i_4\ne 0\}}
{\bf 1}_{\{j_3=j_4\}}
\zeta_{j_1}^{(i_1)}
\zeta_{j_2}^{(i_2)}
\zeta_{j_5}^{(i_5)}
\zeta_{j_6}^{(i_6)}-
{\bf 1}_{\{i_3=i_5\ne 0\}}
{\bf 1}_{\{j_3=j_5\}}
\zeta_{j_1}^{(i_1)}
\zeta_{j_2}^{(i_2)}
\zeta_{j_4}^{(i_4)}
\zeta_{j_6}^{(i_6)}-
$$
$$
-
{\bf 1}_{\{i_4=i_5\ne 0\}}
{\bf 1}_{\{j_4=j_5\}}
\zeta_{j_1}^{(i_1)}
\zeta_{j_2}^{(i_2)}
\zeta_{j_3}^{(i_3)}
\zeta_{j_6}^{(i_6)}+
$$
$$
+
{\bf 1}_{\{i_1=i_2\ne 0\}}
{\bf 1}_{\{j_1=j_2\}}
{\bf 1}_{\{i_3=i_4\ne 0\}}
{\bf 1}_{\{j_3=j_4\}}
\zeta_{j_5}^{(i_5)}
\zeta_{j_6}^{(i_6)}+
{\bf 1}_{\{i_1=i_2\ne 0\}}
{\bf 1}_{\{j_1=j_2\}}
{\bf 1}_{\{i_3=i_5\ne 0\}}
{\bf 1}_{\{j_3=j_5\}}
\zeta_{j_4}^{(i_4)}
\zeta_{j_6}^{(i_6)}+
$$
$$
+
{\bf 1}_{\{i_1=i_2\ne 0\}}
{\bf 1}_{\{j_1=j_2\}}
{\bf 1}_{\{i_4=i_5\ne 0\}}
{\bf 1}_{\{j_4=j_5\}}
\zeta_{j_3}^{(i_3)}
\zeta_{j_6}^{(i_6)}
+
{\bf 1}_{\{i_1=i_3\ne 0\}}
{\bf 1}_{\{j_1=j_3\}}
{\bf 1}_{\{i_2=i_4\ne 0\}}
{\bf 1}_{\{j_2=j_4\}}
\zeta_{j_5}^{(i_5)}
\zeta_{j_6}^{(i_6)}+
$$
$$
+
{\bf 1}_{\{i_1=i_3\ne 0\}}
{\bf 1}_{\{j_1=j_3\}}
{\bf 1}_{\{i_2=i_5\ne 0\}}
{\bf 1}_{\{j_2=j_5\}}
\zeta_{j_4}^{(i_4)}
\zeta_{j_6}^{(i_6)}
+{\bf 1}_{\{i_1=i_3\ne 0\}}
{\bf 1}_{\{j_1=j_3\}}
{\bf 1}_{\{i_4=i_5\ne 0\}}
{\bf 1}_{\{j_4=j_5\}}
\zeta_{j_2}^{(i_2)}
\zeta_{j_6}^{(i_6)}+
$$
$$
+
{\bf 1}_{\{i_1=i_4\ne 0\}}
{\bf 1}_{\{j_1=j_4\}}
{\bf 1}_{\{i_2=i_3\ne 0\}}
{\bf 1}_{\{j_2=j_3\}}
\zeta_{j_5}^{(i_5)}
\zeta_{j_6}^{(i_6)}
+
{\bf 1}_{\{i_1=i_4\ne 0\}}
{\bf 1}_{\{j_1=j_4\}}
{\bf 1}_{\{i_2=i_5\ne 0\}}
{\bf 1}_{\{j_2=j_5\}}
\zeta_{j_3}^{(i_3)}
\zeta_{j_6}^{(i_6)}+
$$
$$
+
{\bf 1}_{\{i_1=i_4\ne 0\}}
{\bf 1}_{\{j_1=j_4\}}
{\bf 1}_{\{i_3=i_5\ne 0\}}
{\bf 1}_{\{j_3=j_5\}}
\zeta_{j_2}^{(i_2)}
\zeta_{j_6}^{(i_6)}
+
{\bf 1}_{\{i_1=i_5\ne 0\}}
{\bf 1}_{\{j_1=j_5\}}
{\bf 1}_{\{i_2=i_3\ne 0\}}
{\bf 1}_{\{j_2=j_3\}}
\zeta_{j_4}^{(i_4)}
\zeta_{j_6}^{(i_6)}+
$$
$$
+
{\bf 1}_{\{i_1=i_5\ne 0\}}
{\bf 1}_{\{j_1=j_5\}}
{\bf 1}_{\{i_2=i_4\ne 0\}}
{\bf 1}_{\{j_2=j_4\}}
\zeta_{j_3}^{(i_3)}
\zeta_{j_6}^{(i_6)}
+
{\bf 1}_{\{i_1=i_5\ne 0\}}
{\bf 1}_{\{j_1=j_5\}}
{\bf 1}_{\{i_3=i_4\ne 0\}}
{\bf 1}_{\{j_3=j_4\}}
\zeta_{j_2}^{(i_2)}
\zeta_{j_6}^{(i_6)}+
$$
$$
+
{\bf 1}_{\{i_2=i_3\ne 0\}}
{\bf 1}_{\{j_2=j_3\}}
{\bf 1}_{\{i_4=i_5\ne 0\}}
{\bf 1}_{\{j_4=j_5\}}
\zeta_{j_1}^{(i_1)}
\zeta_{j_6}^{(i_6)}
+
{\bf 1}_{\{i_2=i_4\ne 0\}}
{\bf 1}_{\{j_2=j_4\}}
{\bf 1}_{\{i_3=i_5\ne 0\}}
{\bf 1}_{\{j_3=j_5\}}
\zeta_{j_1}^{(i_1)}
\zeta_{j_6}^{(i_6)}+
$$
$$
+
{\bf 1}_{\{i_2=i_5\ne 0\}}
{\bf 1}_{\{j_2=j_5\}}
{\bf 1}_{\{i_3=i_4\ne 0\}}
{\bf 1}_{\{j_3=j_4\}}
\zeta_{j_1}^{(i_1)}
\zeta_{j_6}^{(i_6)}
+
{\bf 1}_{\{i_6=i_1\ne 0\}}
{\bf 1}_{\{j_6=j_1\}}
{\bf 1}_{\{i_3=i_4\ne 0\}}
{\bf 1}_{\{j_3=j_4\}}
\zeta_{j_2}^{(i_2)}
\zeta_{j_5}^{(i_5)}+
$$
$$
+
{\bf 1}_{\{i_6=i_1\ne 0\}}
{\bf 1}_{\{j_6=j_1\}}
{\bf 1}_{\{i_3=i_5\ne 0\}}
{\bf 1}_{\{j_3=j_5\}}
\zeta_{j_2}^{(i_2)}
\zeta_{j_4}^{(i_4)}
+
{\bf 1}_{\{i_6=i_1\ne 0\}}
{\bf 1}_{\{j_6=j_1\}}
{\bf 1}_{\{i_2=i_5\ne 0\}}
{\bf 1}_{\{j_2=j_5\}}
\zeta_{j_3}^{(i_3)}
\zeta_{j_4}^{(i_4)}+
$$
$$
+
{\bf 1}_{\{i_6=i_1\ne 0\}}
{\bf 1}_{\{j_6=j_1\}}
{\bf 1}_{\{i_2=i_4\ne 0\}}
{\bf 1}_{\{j_2=j_4\}}
\zeta_{j_3}^{(i_3)}
\zeta_{j_5}^{(i_5)}
+
{\bf 1}_{\{i_6=i_1\ne 0\}}
{\bf 1}_{\{j_6=j_1\}}
{\bf 1}_{\{i_4=i_5\ne 0\}}
{\bf 1}_{\{j_4=j_5\}}
\zeta_{j_2}^{(i_2)}
\zeta_{j_3}^{(i_3)}+
$$
$$
+
{\bf 1}_{\{i_6=i_1\ne 0\}}
{\bf 1}_{\{j_6=j_1\}}
{\bf 1}_{\{i_2=i_3\ne 0\}}
{\bf 1}_{\{j_2=j_3\}}
\zeta_{j_4}^{(i_4)}
\zeta_{j_5}^{(i_5)}
+
{\bf 1}_{\{i_6=i_2\ne 0\}}
{\bf 1}_{\{j_6=j_2\}}
{\bf 1}_{\{i_3=i_5\ne 0\}}
{\bf 1}_{\{j_3=j_5\}}
\zeta_{j_1}^{(i_1)}
\zeta_{j_4}^{(i_4)}+
$$
$$
+
{\bf 1}_{\{i_6=i_2\ne 0\}}
{\bf 1}_{\{j_6=j_2\}}
{\bf 1}_{\{i_4=i_5\ne 0\}}
{\bf 1}_{\{j_4=j_5\}}
\zeta_{j_1}^{(i_1)}
\zeta_{j_3}^{(i_3)}
+
{\bf 1}_{\{i_6=i_2\ne 0\}}
{\bf 1}_{\{j_6=j_2\}}
{\bf 1}_{\{i_3=i_4\ne 0\}}
{\bf 1}_{\{j_3=j_4\}}
\zeta_{j_1}^{(i_1)}
\zeta_{j_5}^{(i_5)}+
$$
$$
+
{\bf 1}_{\{i_6=i_2\ne 0\}}
{\bf 1}_{\{j_6=j_2\}}
{\bf 1}_{\{i_1=i_5\ne 0\}}
{\bf 1}_{\{j_1=j_5\}}
\zeta_{j_3}^{(i_3)}
\zeta_{j_4}^{(i_4)}
+
{\bf 1}_{\{i_6=i_2\ne 0\}}
{\bf 1}_{\{j_6=j_2\}}
{\bf 1}_{\{i_1=i_4\ne 0\}}
{\bf 1}_{\{j_1=j_4\}}
\zeta_{j_3}^{(i_3)}
\zeta_{j_5}^{(i_5)}+
$$
$$
+
{\bf 1}_{\{i_6=i_2\ne 0\}}
{\bf 1}_{\{j_6=j_2\}}
{\bf 1}_{\{i_1=i_3\ne 0\}}
{\bf 1}_{\{j_1=j_3\}}
\zeta_{j_4}^{(i_4)}
\zeta_{j_5}^{(i_5)}
+
{\bf 1}_{\{i_6=i_3\ne 0\}}
{\bf 1}_{\{j_6=j_3\}}
{\bf 1}_{\{i_2=i_5\ne 0\}}
{\bf 1}_{\{j_2=j_5\}}
\zeta_{j_1}^{(i_1)}
\zeta_{j_4}^{(i_4)}+
$$
$$
+
{\bf 1}_{\{i_6=i_3\ne 0\}}
{\bf 1}_{\{j_6=j_3\}}
{\bf 1}_{\{i_4=i_5\ne 0\}}
{\bf 1}_{\{j_4=j_5\}}
\zeta_{j_1}^{(i_1)}
\zeta_{j_2}^{(i_2)}
+
{\bf 1}_{\{i_6=i_3\ne 0\}}
{\bf 1}_{\{j_6=j_3\}}
{\bf 1}_{\{i_2=i_4\ne 0\}}
{\bf 1}_{\{j_2=j_4\}}
\zeta_{j_1}^{(i_1)}
\zeta_{j_5}^{(i_5)}+
$$
$$
+
{\bf 1}_{\{i_6=i_3\ne 0\}}
{\bf 1}_{\{j_6=j_3\}}
{\bf 1}_{\{i_1=i_5\ne 0\}}
{\bf 1}_{\{j_1=j_5\}}
\zeta_{j_2}^{(i_2)}
\zeta_{j_4}^{(i_4)}
+
{\bf 1}_{\{i_6=i_3\ne 0\}}
{\bf 1}_{\{j_6=j_3\}}
{\bf 1}_{\{i_1=i_4\ne 0\}}
{\bf 1}_{\{j_1=j_4\}}
\zeta_{j_2}^{(i_2)}
\zeta_{j_5}^{(i_5)}+
$$
$$
+
{\bf 1}_{\{i_6=i_3\ne 0\}}
{\bf 1}_{\{j_6=j_3\}}
{\bf 1}_{\{i_1=i_2\ne 0\}}
{\bf 1}_{\{j_1=j_2\}}
\zeta_{j_4}^{(i_4)}
\zeta_{j_5}^{(i_5)}
+
{\bf 1}_{\{i_6=i_4\ne 0\}}
{\bf 1}_{\{j_6=j_4\}}
{\bf 1}_{\{i_3=i_5\ne 0\}}
{\bf 1}_{\{j_3=j_5\}}
\zeta_{j_1}^{(i_1)}
\zeta_{j_2}^{(i_2)}+
$$
$$
+
{\bf 1}_{\{i_6=i_4\ne 0\}}
{\bf 1}_{\{j_6=j_4\}}
{\bf 1}_{\{i_2=i_5\ne 0\}}
{\bf 1}_{\{j_2=j_5\}}
\zeta_{j_1}^{(i_1)}
\zeta_{j_3}^{(i_3)}
+
{\bf 1}_{\{i_6=i_4\ne 0\}}
{\bf 1}_{\{j_6=j_4\}}
{\bf 1}_{\{i_2=i_3\ne 0\}}
{\bf 1}_{\{j_2=j_3\}}
\zeta_{j_1}^{(i_1)}
\zeta_{j_5}^{(i_5)}+
$$
$$
+
{\bf 1}_{\{i_6=i_4\ne 0\}}
{\bf 1}_{\{j_6=j_4\}}
{\bf 1}_{\{i_1=i_5\ne 0\}}
{\bf 1}_{\{j_1=j_5\}}
\zeta_{j_2}^{(i_2)}
\zeta_{j_3}^{(i_3)}
+
{\bf 1}_{\{i_6=i_4\ne 0\}}
{\bf 1}_{\{j_6=j_4\}}
{\bf 1}_{\{i_1=i_3\ne 0\}}
{\bf 1}_{\{j_1=j_3\}}
\zeta_{j_2}^{(i_2)}
\zeta_{j_5}^{(i_5)}+
$$
$$
+
{\bf 1}_{\{i_6=i_4\ne 0\}}
{\bf 1}_{\{j_6=j_4\}}
{\bf 1}_{\{i_1=i_2\ne 0\}}
{\bf 1}_{\{j_1=j_2\}}
\zeta_{j_3}^{(i_3)}
\zeta_{j_5}^{(i_5)}
+
{\bf 1}_{\{i_6=i_5\ne 0\}}
{\bf 1}_{\{j_6=j_5\}}
{\bf 1}_{\{i_3=i_4\ne 0\}}
{\bf 1}_{\{j_3=j_4\}}
\zeta_{j_1}^{(i_1)}
\zeta_{j_2}^{(i_2)}+
$$
$$
+
{\bf 1}_{\{i_6=i_5\ne 0\}}
{\bf 1}_{\{j_6=j_5\}}
{\bf 1}_{\{i_2=i_4\ne 0\}}
{\bf 1}_{\{j_2=j_4\}}
\zeta_{j_1}^{(i_1)}
\zeta_{j_3}^{(i_3)}
+
{\bf 1}_{\{i_6=i_5\ne 0\}}
{\bf 1}_{\{j_6=j_5\}}
{\bf 1}_{\{i_2=i_3\ne 0\}}
{\bf 1}_{\{j_2=j_3\}}
\zeta_{j_1}^{(i_1)}
\zeta_{j_4}^{(i_4)}+
$$
$$
+
{\bf 1}_{\{i_6=i_5\ne 0\}}
{\bf 1}_{\{j_6=j_5\}}
{\bf 1}_{\{i_1=i_4\ne 0\}}
{\bf 1}_{\{j_1=j_4\}}
\zeta_{j_2}^{(i_2)}
\zeta_{j_3}^{(i_3)}
+
{\bf 1}_{\{i_6=i_5\ne 0\}}
{\bf 1}_{\{j_6=j_5\}}
{\bf 1}_{\{i_1=i_3\ne 0\}}
{\bf 1}_{\{j_1=j_3\}}
\zeta_{j_2}^{(i_2)}
\zeta_{j_4}^{(i_4)}+
$$
$$
+
{\bf 1}_{\{i_6=i_5\ne 0\}}
{\bf 1}_{\{j_6=j_5\}}
{\bf 1}_{\{i_1=i_2\ne 0\}}
{\bf 1}_{\{j_1=j_2\}}
\zeta_{j_3}^{(i_3)}
\zeta_{j_4}^{(i_4)}-
$$
$$
-
{\bf 1}_{\{i_6=i_1\ne 0\}}
{\bf 1}_{\{j_6=j_1\}}
{\bf 1}_{\{i_2=i_5\ne 0\}}
{\bf 1}_{\{j_2=j_5\}}
{\bf 1}_{\{i_3=i_4\ne 0\}}
{\bf 1}_{\{j_3=j_4\}}-
$$
$$
-
{\bf 1}_{\{i_6=i_1\ne 0\}}
{\bf 1}_{\{j_6=j_1\}}
{\bf 1}_{\{i_2=i_4\ne 0\}}
{\bf 1}_{\{j_2=j_4\}}
{\bf 1}_{\{i_3=i_5\ne 0\}}
{\bf 1}_{\{j_3=j_5\}}-
$$
$$
-
{\bf 1}_{\{i_6=i_1\ne 0\}}
{\bf 1}_{\{j_6=j_1\}}
{\bf 1}_{\{i_2=i_3\ne 0\}}
{\bf 1}_{\{j_2=j_3\}}
{\bf 1}_{\{i_4=i_5\ne 0\}}
{\bf 1}_{\{j_4=j_5\}}-
$$
$$
-
{\bf 1}_{\{i_6=i_2\ne 0\}}
{\bf 1}_{\{j_6=j_2\}}
{\bf 1}_{\{i_1=i_5\ne 0\}}
{\bf 1}_{\{j_1=j_5\}}
{\bf 1}_{\{i_3=i_4\ne 0\}}
{\bf 1}_{\{j_3=j_4\}}-
$$
$$
-
{\bf 1}_{\{i_6=i_2\ne 0\}}
{\bf 1}_{\{j_6=j_2\}}
{\bf 1}_{\{i_1=i_4\ne 0\}}
{\bf 1}_{\{j_1=j_4\}}
{\bf 1}_{\{i_3=i_5\ne 0\}}
{\bf 1}_{\{j_3=j_5\}}-
$$
$$
-
{\bf 1}_{\{i_6=i_2\ne 0\}}
{\bf 1}_{\{j_6=j_2\}}
{\bf 1}_{\{i_1=i_3\ne 0\}}
{\bf 1}_{\{j_1=j_3\}}
{\bf 1}_{\{i_4=i_5\ne 0\}}
{\bf 1}_{\{j_4=j_5\}}-
$$
$$
-
{\bf 1}_{\{i_6=i_3\ne 0\}}
{\bf 1}_{\{j_6=j_3\}}
{\bf 1}_{\{i_1=i_5\ne 0\}}
{\bf 1}_{\{j_1=j_5\}}
{\bf 1}_{\{i_2=i_4\ne 0\}}
{\bf 1}_{\{j_2=j_4\}}-
$$
$$
-
{\bf 1}_{\{i_6=i_3\ne 0\}}
{\bf 1}_{\{j_6=j_3\}}
{\bf 1}_{\{i_1=i_4\ne 0\}}
{\bf 1}_{\{j_1=j_4\}}
{\bf 1}_{\{i_2=i_5\ne 0\}}
{\bf 1}_{\{j_2=j_5\}}-
$$
$$
-
{\bf 1}_{\{i_3=i_6\ne 0\}}
{\bf 1}_{\{j_3=j_6\}}
{\bf 1}_{\{i_1=i_2\ne 0\}}
{\bf 1}_{\{j_1=j_2\}}
{\bf 1}_{\{i_4=i_5\ne 0\}}
{\bf 1}_{\{j_4=j_5\}}-
$$
$$
-
{\bf 1}_{\{i_6=i_4\ne 0\}}
{\bf 1}_{\{j_6=j_4\}}
{\bf 1}_{\{i_1=i_5\ne 0\}}
{\bf 1}_{\{j_1=j_5\}}
{\bf 1}_{\{i_2=i_3\ne 0\}}
{\bf 1}_{\{j_2=j_3\}}-
$$
$$
-
{\bf 1}_{\{i_6=i_4\ne 0\}}
{\bf 1}_{\{j_6=j_4\}}
{\bf 1}_{\{i_1=i_3\ne 0\}}
{\bf 1}_{\{j_1=j_3\}}
{\bf 1}_{\{i_2=i_5\ne 0\}}
{\bf 1}_{\{j_2=j_5\}}-
$$
$$
-
{\bf 1}_{\{i_6=i_4\ne 0\}}
{\bf 1}_{\{j_6=j_4\}}
{\bf 1}_{\{i_1=i_2\ne 0\}}
{\bf 1}_{\{j_1=j_2\}}
{\bf 1}_{\{i_3=i_5\ne 0\}}
{\bf 1}_{\{j_3=j_5\}}-
$$
$$
-
{\bf 1}_{\{i_6=i_5\ne 0\}}
{\bf 1}_{\{j_6=j_5\}}
{\bf 1}_{\{i_1=i_4\ne 0\}}
{\bf 1}_{\{j_1=j_4\}}
{\bf 1}_{\{i_2=i_3\ne 0\}}
{\bf 1}_{\{j_2=j_3\}}-
$$
$$
-
{\bf 1}_{\{i_6=i_5\ne 0\}}
{\bf 1}_{\{j_6=j_5\}}
{\bf 1}_{\{i_1=i_2\ne 0\}}
{\bf 1}_{\{j_1=j_2\}}
{\bf 1}_{\{i_3=i_4\ne 0\}}
{\bf 1}_{\{j_3=j_4\}}-
$$
\begin{equation}
\label{a6}
\Biggl.-
{\bf 1}_{\{i_6=i_5\ne 0\}}
{\bf 1}_{\{j_6=j_5\}}
{\bf 1}_{\{i_1=i_3\ne 0\}}
{\bf 1}_{\{j_1=j_3\}}
{\bf 1}_{\{i_2=i_4\ne 0\}}
{\bf 1}_{\{j_2=j_4\}}\Biggr),
\end{equation}

\vspace{5mm}
\noindent
where ${\bf 1}_A$ is the indicator of the set $A$.

The cases $k=7$ and $k>7$ are considered in 
\cite{1}-\cite{new-2023a} (also see Sect.~5).

\vspace{5mm}

\section{Exact Calculation of the Mean-Square Approximation Error in The Method of 
Generalized Multiple Fourier Series. The Case of Complete Orthonormal Systems  
of Continuous Functions in the Space $L_2([t,T])$ 
and Continuous Weight Functions $\psi_1(\tau),\ldots,\psi_k(\tau)$}

\vspace{5mm}

{\bf Theorem 2} \cite{arxiv-3} (also see \cite{11}-\cite{12}, \cite{art-1}). 
{\it Suppose that
every $\psi_l(\tau)$ $(l=1,\ldots, k)$ is a continuous nonrandom function on 
$[t, T]$ and
$\{\phi_j(x)\}_{j=0}^{\infty}$ is a complete orthonormal system  
of continuous functions in the space $L_2([t,T]).$ Then

\vspace{1mm}
$$
{\sf M}\left\{\left(J[\psi^{(k)}]_{T,t}-
J[\psi^{(k)}]_{T,t}^p\right)^2\right\}
= \int\limits_{[t,T]^k} K^2(t_1,\ldots,t_k)
dt_1\ldots dt_k - 
$$

\vspace{2mm}
\begin{equation}
\label{tttr11}
- \sum_{j_1=0}^{p}\ldots\sum_{j_k=0}^{p}
C_{j_k\ldots j_1}
{\sf M}\left\{J[\psi^{(k)}]_{T,t}
\sum\limits_{(j_1,\ldots,j_k)}
\int\limits_t^T \phi_{j_k}(t_k)
\ldots
\int\limits_t^{t_{2}}\phi_{j_{1}}(t_{1})
d{\bf f}_{t_1}^{(i_1)}\ldots
d{\bf f}_{t_k}^{(i_k)}\right\},
\end{equation}

\vspace{5mm}
\noindent
where
$$
J[\psi^{(k)}]_{T,t}=\int\limits_t^T\psi_k(t_k) \ldots \int\limits_t^{t_{2}}
\psi_1(t_1) d{\bf f}_{t_1}^{(i_1)}\ldots
d{\bf f}_{t_k}^{(i_k)},
$$

\vspace{2mm}
\begin{equation}
\label{yeee2}
J[\psi^{(k)}]_{T,t}^p=
\sum_{j_1=0}^{p}\ldots\sum_{j_k=0}^{p}
C_{j_k\ldots j_1}\left(
\prod_{l=1}^k\zeta_{j_l}^{(i_l)}-S_{j_1,\ldots,j_k}^{(i_1\ldots i_k)}
\right),
\end{equation}

\vspace{4mm}
\begin{equation}
\label{ppp1}
S_{j_1,\ldots,j_k}^{(i_1\ldots i_k)}=
\hbox{\vtop{\offinterlineskip\halign{
\hfil#\hfil\cr
{\rm l.i.m.}\cr
$\stackrel{}{{}_{N\to \infty}}$\cr
}} }\sum_{(l_1,\ldots,l_k)\in {\rm G}_k}
\phi_{j_{1}}(\tau_{l_1})
\Delta{\bf f}_{\tau_{l_1}}^{(i_1)}\ldots
\phi_{j_{k}}(\tau_{l_k})
\Delta{\bf f}_{\tau_{l_k}}^{(i_k)},
\end{equation}

\vspace{6mm}
\noindent
the Fourier coefficient $C_{j_k\ldots j_1}$ has the form {\rm (\ref{ppppa})},

\vspace{-2mm}
\begin{equation}
\label{rr232}
\zeta_{j}^{(i)}=
\int\limits_t^T \phi_{j}(s) d{\bf f}_s^{(i)}
\end{equation}

\vspace{2mm}
\noindent
are independent standard Gaussian random variables
for various
$i$ or $j$ $(i=1,\ldots,m),$

\vspace{-2mm}
$$
\sum\limits_{(j_1,\ldots,j_k)}
$$ 

\vspace{2mm}
\noindent
means the sum with respect to all
possible permutations 
$(j_1,\ldots,j_k)$. At the same time if 
$j_r$ swapped with $j_q$ in the permutation $(j_1,\ldots,j_k)$,
then $i_r$ swapped with $i_q$ in the permutation
$(i_1,\ldots,i_k);$
another notations are the same as in Theorem {\rm 1.}}

\vspace{2mm}

{\bf Remark 2.}\ {\it Note that

\vspace{1mm}
$$
{\sf M}\left\{J[\psi^{(k)}]_{T,t}
\int\limits_t^T \phi_{j_k}(t_k)
\ldots
\int\limits_t^{t_{2}}\phi_{j_{1}}(t_{1})
d{\bf f}_{t_1}^{(i_1)}\ldots
d{\bf f}_{t_k}^{(i_k)}\right\}=
$$

\vspace{2mm}
$$
={\sf M}\left\{\int\limits_t^T\psi_k(t_k) \ldots \int\limits_t^{t_{2}}
\psi_1(t_1) d{\bf f}_{t_1}^{(i_1)}\ldots
d{\bf f}_{t_k}^{(i_k)}
\int\limits_t^T \phi_{j_k}(t_k)
\ldots
\int\limits_t^{t_{2}}\phi_{j_{1}}(t_{1})
d{\bf f}_{t_1}^{(i_1)}\ldots
d{\bf f}_{t_k}^{(i_k)}\right\}=
$$

\vspace{2mm}
$$
=\int\limits_t^T\psi_k(t_k) \phi_{j_k}(t_k)\ldots \int\limits_t^{t_{2}}
\psi_1(t_1)\phi_{j_1}(t_1) dt_1\ldots dt_k=
C_{j_k\ldots j_1}.
$$

\vspace{5mm}

Therefore, from Theorem {\rm 2} for the case of pairwise different numbers
$i_1,\ldots,i_k$ we obtain 

\vspace{2mm}
$$
{\sf M}\left\{\left(J[\psi^{(k)}]_{T,t}-
J[\psi^{(k)}]_{T,t}^p\right)^2\right\}=
$$

\vspace{3mm}
$$
=
\int\limits_{[t,T]^k} K^2(t_1,\ldots,t_k)
dt_1\ldots dt_k - \sum_{j_1=0}^{p}\ldots\sum_{j_k=0}^{p}
C_{j_k\ldots j_1}^2.
$$

\vspace{6mm}

Moreover, if $i_1=\ldots=i_k,$ then from Theorem {\rm 2} we get

\vspace{2mm}
$$
{\sf M}\left\{\left(J[\psi^{(k)}]_{T,t}-
J[\psi^{(k)}]_{T,t}^p\right)^2\right\}=
$$

\vspace{3mm}
$$
=\int\limits_{[t,T]^k} K^2(t_1,\ldots,t_k)
dt_1\ldots dt_k - \sum_{j_1=0}^{p}\ldots\sum_{j_k=0}^{p}
C_{j_k\ldots j_1}\Biggl(\sum\limits_{(j_1,\ldots,j_k)}
C_{j_k\ldots j_1}\Biggr),
$$

\vspace{6mm}
\noindent
where
$$
\sum\limits_{(j_1,\ldots,j_k)}
$$ 

\vspace{2mm}
\noindent
means the sum with respect to all
possible permutations
$(j_1,\ldots,j_k).$

For example, for the case $k=3$ we have

\vspace{1mm}
$$
{\sf M}\left\{\left(J[\psi^{(3)}]_{T,t}-
J[\psi^{(3)}]_{T,t}^p\right)^2\right\}=
\int\limits_t^T\psi_3^2(t_3)\int\limits_t^{t_3}\psi_2^2(t_2)
\int\limits_t^{t_2}\psi_1^2(t_1)dt_1dt_2dt_3 -
$$

\vspace{3mm}
$$
- \sum_{j_1,j_2,j_3=0}^{p}
C_{j_3j_2j_1}\biggl(C_{j_3j_2j_1}+C_{j_3j_1j_2}+C_{j_2j_3j_1}+
C_{j_2j_1j_3}+C_{j_1j_2j_3}+C_{j_1j_3j_2}\biggr).
$$
}

\vspace{6mm}

{\bf Proof.}
Using Theorem 1 for the case $p_1=\ldots=p_k=p$, we obtain

\begin{equation}
\label{yyye1}
J[\psi^{(k)}]_{T,t}=\
\hbox{\vtop{\offinterlineskip\halign{
\hfil#\hfil\cr
{\rm l.i.m.}\cr
$\stackrel{}{{}_{p\to \infty}}$\cr
}} }\sum_{j_1=0}^{p}\ldots\sum_{j_k=0}^{p}
C_{j_k\ldots j_1}\left(
\prod_{l=1}^k\zeta_{j_l}^{(i_l)}-S_{j_1,\ldots,j_k}^{(i_1\ldots i_k)}
\right).
\end{equation}

\vspace{3mm}

For $n>p$ we can write 

\vspace{2mm}
$$
J[\psi^{(k)}]_{T,t}^n=
\left(\sum_{j_1=0}^{p}+\sum_{j_1=p+1}^n\right)\ldots
\left(\sum_{j_k=0}^{p}+\sum_{j_k=p+1}^n\right)
C_{j_k\ldots j_1}\left(
\prod_{l=1}^k\zeta_{j_l}^{(i_l)}-S_{j_1,\ldots,j_k}^{(i_1\ldots i_k)}
\right)=
$$

\vspace{4mm}
\begin{equation}
\label{yyye}
=J[\psi^{(k)}]_{T,t}^p + \xi[\psi^{(k)}]_{T,t}^{p+1,n}.
\end{equation}

\vspace{4mm}

Let us prove that due to the special structure of random variables 
$S_{j_1,\ldots,j_k}^{(i_1\ldots i_k)}$ (see (\ref{leto5001})--(\ref{a6}),
(\ref{ppp1})),
the following relations 
are correct

\vspace{-1mm}
\begin{equation}
\label{tyty}
{\sf M}\left\{
\prod_{l=1}^k\zeta_{j_l}^{(i_l)}-S_{j_1,\ldots,j_k}^{(i_1\ldots i_k)}
\right\}=0,
\end{equation}

\vspace{2mm}
\begin{equation}
\label{tyty1}
{\sf M}\left\{
\left(\prod_{l=1}^k\zeta_{j_l}^{(i_l)}-S_{j_1,\ldots,j_k}^{(i_1\ldots i_k)}
\right)
\left(\prod_{l=1}^k\zeta_{j_l'}^{(i_l)}-S_{j_1',\ldots,j_k'}^{(i_1\ldots i_k)}
\right)\right\}=0,
\end{equation}

\vspace{5mm}
\noindent
where
$$
(j_1,\ldots,j_k)\in{\rm K}_p,\ \ \ (j_1',\ldots,j_k')
\in{\rm K}_n\backslash {\rm K}_{p}
$$ 

\vspace{2mm}
\noindent
and

\vspace{-1mm}
$$
{\rm K}_n=\left\{(j_1,\ldots,j_k):\ 0\le j_1,\ldots,j_k\le n\right\},
$$

$$
{\rm K}_p=\left\{(j_1,\ldots,j_k):\ 0\le j_1,\ldots,j_k\le p\right\}.
$$

\vspace{6mm}

For the case $i_1,\ldots,i_k=0, 1,\ldots,m$ from the
proof of Theorem 1 in \cite{arxiv-1}
(also see \cite{1}-\cite{art-9}, \cite{arxiv-2}-\cite{new-2023a}) 
it follows that

\vspace{2mm}
$$
J[\psi^{(k)}]_{T,t}
=
\sum_{j_1=0}^{p_1}\ldots
\sum_{j_k=0}^{p_k}
C_{j_k\ldots j_1}\ \
\hbox{\vtop{\offinterlineskip\halign{
\hfil#\hfil\cr
{\rm l.i.m.}\cr
$\stackrel{}{{}_{N\to \infty}}$\cr
}} }
\sum\limits_{\stackrel{l_1,\ldots,l_k=0}{{}_{l_q\ne l_r;\ q\ne r;\ 
q, r=1,\ldots, k}}}^{N-1}
\phi_{j_1}(\tau_{l_1})\ldots
\phi_{j_k}(\tau_{l_k})
\Delta{\bf w}_{\tau_{l_1}}^{(i_1)}
\ldots
\Delta{\bf w}_{\tau_{l_k}}^{(i_k)}+
$$

\vspace{4mm}
$$
+R_{T,t}^{p_1,\ldots,p_k}=
$$

\vspace{6mm}
$$
=\sum_{j_1=0}^{p_1}\ldots
\sum_{j_k=0}^{p_k}
C_{j_k\ldots j_1}\Biggl(
\hbox{\vtop{\offinterlineskip\halign{
\hfil#\hfil\cr
{\rm l.i.m.}\cr
$\stackrel{}{{}_{N\to \infty}}$\cr
}} }\sum_{l_1,\ldots,l_k=0}^{N-1}
\phi_{j_1}(\tau_{l_1})
\ldots
\phi_{j_k}(\tau_{l_k})
\Delta{\bf w}_{\tau_{l_1}}^{(i_1)}
\ldots
\Delta{\bf w}_{\tau_{l_k}}^{(i_k)}
-\Biggr.
$$

\vspace{3mm}
$$
-\Biggl.
\hbox{\vtop{\offinterlineskip\halign{
\hfil#\hfil\cr
{\rm l.i.m.}\cr
$\stackrel{}{{}_{N\to \infty}}$\cr
}} }\sum_{(l_1,\ldots,l_k)\in {\rm G}_k}
\phi_{j_{1}}(\tau_{l_1})
\Delta{\bf w}_{\tau_{l_1}}^{(i_1)}\ldots
\phi_{j_{k}}(\tau_{l_k})
\Delta{\bf w}_{\tau_{l_k}}^{(i_k)}\Biggr)
+R_{T,t}^{p_1,\ldots,p_k}=
$$

\vspace{7mm}
$$
=\sum_{j_1=0}^{p_1}\ldots\sum_{j_k=0}^{p_k}
C_{j_k\ldots j_1}\Biggl(
\prod_{l=1}^k\zeta_{j_l}^{(i_l)}-
\hbox{\vtop{\offinterlineskip\halign{
\hfil#\hfil\cr
{\rm l.i.m.}\cr
$\stackrel{}{{}_{N\to \infty}}$\cr
}} }\sum_{(l_1,\ldots,l_k)\in {\rm G}_k}
\phi_{j_{1}}(\tau_{l_1})
\Delta{\bf w}_{\tau_{l_1}}^{(i_1)}\ldots
\phi_{j_{k}}(\tau_{l_k})
\Delta{\bf w}_{\tau_{l_k}}^{(i_k)}\Biggr)+
$$

\vspace{3mm}
\begin{equation}
\label{e1}
+R_{T,t}^{p_1,\ldots,p_k}\ \ \ \hbox{w.~p.~1},
\end{equation}

\vspace{5mm}
\noindent
where

$$
R_{T,t}^{p_1,\ldots,p_k}=\sum_{(t_1,\ldots,t_k)}
\int\limits_{t}^{T}
\ldots
\int\limits_{t}^{t_2}
\Biggl(K(t_1,\ldots,t_k)-
\sum_{j_1=0}^{p_1}\ldots
\sum_{j_k=0}^{p_k}
C_{j_k\ldots j_1}
\prod_{l=1}^k\phi_{j_l}(t_l)\Biggr)\times
$$

\vspace{1mm}

\begin{equation}
\label{e2}
\times
d{\bf w}_{t_1}^{(i_1)}
\ldots
d{\bf w}_{t_k}^{(i_k)},
\end{equation}

\vspace{2mm}
\noindent
where

\vspace{-1mm}
$$
\sum\limits_{(t_1,\ldots,t_k)}
$$ 

\vspace{1mm}
\noindent
means the sum with respect to all
possible permutations
$(t_1,\ldots,t_k),$ which are
performed only 
in the values $d{\bf w}_{t_1}^{(i_1)}
\ldots $
$d{\bf w}_{t_k}^{(i_k)}$. At the same time the indexes near 
upper limits of integration in the iterated stochastic integrals 
(see (\ref{e2}))
are changed correspondently and if $t_r$ swapped with $t_q$ in the  
permutation $(t_1,\ldots,t_k)$, then $i_r$ swapped with $i_q$ in the 
permutation $(i_1,\ldots,i_k)$.

For the case $i_1,\ldots,i_k=1,\ldots,m$ 
and $p_1=\ldots=p_k=p$ from (\ref{e1}) we obtain

\vspace{2mm}
$$
\prod_{l=1}^k\zeta_{j_l}^{(i_l)}-
S_{j_1,\ldots,j_k}^{(i_1\ldots i_k)}\ \ =\ \ 
\hbox{\vtop{\offinterlineskip\halign{
\hfil#\hfil\cr
{\rm l.i.m.}\cr
$\stackrel{}{{}_{N\to \infty}}$\cr
}} }
\sum_{\stackrel{l_1,\ldots,l_k=0}
{{}_{l_q\ne l_r;\ q\ne r;\
q,r=1,\ldots,k}}}^{N-1}
\phi_{j_1}(\tau_{l_1})\ldots
\phi_{j_k}(\tau_{l_k})
\Delta{\bf f}_{\tau_{l_1}}^{(i_1)}
\ldots
\Delta{\bf f}_{\tau_{l_k}}^{(i_k)}=
$$

\vspace{2mm}
\begin{equation}
\label{ttt2}
=\sum\limits_{(j_1,\ldots,j_k)}
\int\limits_t^T \phi_{j_k}(t_k)
\ldots
\int\limits_t^{t_{2}}\phi_{j_{1}}(t_{1})
d{\bf f}_{t_1}^{(i_1)}\ldots
d{\bf f}_{t_k}^{(i_k)}\ \ \ {\rm w.\ p.\ 1,}
\end{equation}

\vspace{3mm}
\noindent
where 
$$
\sum\limits_{(j_1,\ldots,j_k)}
$$ 

\vspace{1mm}
\noindent
means the sum with respect to all
possible permutations
$(j_1,\ldots,j_k).$ At the same time if 
$j_r$ swapped  with $j_q$ in the permutation $(j_1,\ldots,j_k)$,
then $i_r$ swapped  with $i_q$ in the permutation
$(i_1,\ldots,i_k);$
another notations are the same as in Theorem 1.

From (\ref{ttt2}) due to the moment property of the Ito
stochastic integral we obtain (\ref{tyty}).
Let us prove (\ref{tyty1}).
From (\ref{ttt2}) we have

\vspace{2mm}
$$
0\le \left\vert{\sf M}\left\{
\left(\prod_{l=1}^k\zeta_{j_l}^{(i_l)}-S_{j_1,\ldots,j_k}^{(i_1\ldots i_k)}
\right)
\left(\prod_{l=1}^k\zeta_{j_l'}^{(i_l)}-S_{j_1',\ldots,j_k'}^{(i_1\ldots i_k)}
\right)\right\}\right\vert=
$$

\vspace{5mm}
$$
=\left\vert
{\sf M}\left\{\sum\limits_{(j_1,\ldots,j_k)}\sum\limits_{(j_1',\ldots,j_k')}
\int\limits_t^T \phi_{j_k}(t_k)
\ldots
\int\limits_t^{t_{2}}\phi_{j_{1}}(t_{1})
d{\bf f}_{t_1}^{(i_1)}\ldots
d{\bf f}_{t_k}^{(i_k)}\times\right.\right.
$$

\vspace{3mm}
$$
\left.\left.\times
\int\limits_t^T \phi_{j_k'}(t_k)
\ldots
\int\limits_t^{t_{2}}\phi_{j_1'}(t_{1})
d{\bf f}_{t_1}^{(i_1)}\ldots
d{\bf f}_{t_k}^{(i_k)}\right\}\right\vert\le
$$

\vspace{3mm}

$$
\le\sum\limits_{(j_1',\ldots,j_k')}\ 
\int\limits_t^T \phi_{j_k}(t_k)\phi_{j_k'}(t_k)dt_k
\ldots
\int\limits_t^{T}\phi_{j_1}(t_{1})\phi_{j_1'}(t_{1})
dt_1=
$$

\vspace{3mm}
\begin{equation}
\label{hq11}
=
\sum\limits_{(j_1',\ldots,j_k')}
{\bf 1}_{\{j_1=j_1'\}}\ldots {\bf 1}_{\{j_k=j_k'\}},
\end{equation}

\vspace{5mm}
\noindent
where 
where ${\bf 1}_A$ is the indicator of the set $A$.
Using (\ref{hq11}), we obtain (\ref{tyty1}).

First, let us prove (\ref{hq11}) for the cases $k=2$ and $k=3$ in detail.
We have

\vspace{1mm}
$$
{\sf M}\left\{\sum\limits_{(j_1,j_2)}\sum\limits_{(j_1',j_2')}
\int\limits_t^T\phi_{j_{2}}(t_{2})
\int\limits_t^{t_{2}}\phi_{j_{1}}(t_{1})
d{\bf f}_{t_1}^{(i_1)}d{\bf f}_{t_2}^{(i_2)}
\int\limits_t^T\phi_{j_{2}'}(t_{2})
\int\limits_t^{t_{2}}\phi_{j_{1}'}(t_{1})
d{\bf f}_{t_1}^{(i_1)}d{\bf f}_{t_2}^{(i_2)}
\right\}=
$$

\vspace{3mm}
$$
=
\int\limits_t^{T}\phi_{j_2}(s)\phi_{j_2'}(s)ds
\int\limits_t^{T}\phi_{j_1}(s)\phi_{j_1'}(s)ds+
$$

\vspace{2mm}
$$
+{\bf 1}_{\{i_1=i_2\}}
\int\limits_t^{T}\phi_{j_2}(s)\phi_{j_1'}(s)ds
\int\limits_t^{T}\phi_{j_1}(s)\phi_{j_2'}(s)ds=
$$

\vspace{2mm}
\begin{equation}
\label{gt22}
=
{\bf 1}_{\{j_1=j_1'\}}{\bf 1}_{\{j_2=j_2'\}}+
{\bf 1}_{\{i_1=i_2\}} \cdot
{\bf 1}_{\{j_2=j_1'\}}{\bf 1}_{\{j_1=j_2'\}},
\end{equation}

\vspace{7mm}
$$
{\sf M}\left\{\sum\limits_{(j_1,j_2,j_3)}\sum\limits_{(j_1',j_2',j_3')}
\int\limits_t^T \phi_{j_3}(t_3)
\int\limits_t^{t_{3}}\phi_{j_{2}}(t_{2})
\int\limits_t^{t_{2}}\phi_{j_{1}}(t_{1})
d{\bf f}_{t_1}^{(i_1)}d{\bf f}_{t_2}^{(i_2)}
d{\bf f}_{t_3}^{(i_3)}\times\right.
$$

\vspace{3mm}
$$
\left.\times
\int\limits_t^T \phi_{j_3'}(t_3)
\int\limits_t^{t_{3}}\phi_{j_{2}'}(t_{2})
\int\limits_t^{t_{2}}\phi_{j_{1}'}(t_{1})
d{\bf f}_{t_1}^{(i_1)}d{\bf f}_{t_2}^{(i_2)}
d{\bf f}_{t_3}^{(i_3)}\right\}=
$$

\vspace{3mm}
$$
=\int\limits_t^T \phi_{j_3}(s)\phi_{j_3'}(s)ds
\int\limits_t^{T}\phi_{j_2}(s)\phi_{j_2'}(s)ds
\int\limits_t^{T}\phi_{j_1}(s)\phi_{j_1'}(s)ds+
$$

\vspace{3mm}
$$
+{\bf 1}_{\{i_1=i_2\}}\int\limits_t^T \phi_{j_3}(s)\phi_{j_3'}(s)ds
\int\limits_t^{T}\phi_{j_1}(s)\phi_{j_2'}(s)ds
\int\limits_t^{T}\phi_{j_2}(s)\phi_{j_1'}(s)ds+
$$

\vspace{3mm}
$$
+{\bf 1}_{\{i_2=i_3\}}\int\limits_t^T \phi_{j_1}(s)\phi_{j_1'}(s)ds
\int\limits_t^{T}\phi_{j_2}(s)\phi_{j_3'}(s)ds
\int\limits_t^{T}\phi_{j_3}(s)\phi_{j_2'}(s)ds+
$$

\vspace{3mm}
$$
+{\bf 1}_{\{i_1=i_3\}}\int\limits_t^T \phi_{j_1}(s)\phi_{j_3'}(s)ds
\int\limits_t^{T}\phi_{j_2}(s)\phi_{j_2'}(s)ds
\int\limits_t^{T}\phi_{j_3}(s)\phi_{j_1'}(s)ds+
$$

\vspace{3mm}
$$
+{\bf 1}_{\{i_1=i_2=i_3\}}\int\limits_t^T \phi_{j_2}(s)\phi_{j_3'}(s)ds
\int\limits_t^{T}\phi_{j_1}(s)\phi_{j_2'}(s)ds
\int\limits_t^{T}\phi_{j_3}(s)\phi_{j_1'}(s)ds+
$$

\vspace{3mm}
$$
+{\bf 1}_{\{i_1=i_2=i_3\}}\int\limits_t^T \phi_{j_1}(s)\phi_{j_3'}(s)ds
\int\limits_t^{T}\phi_{j_3}(s)\phi_{j_2'}(s)ds
\int\limits_t^{T}\phi_{j_2}(s)\phi_{j_1'}(s)ds=
$$

\vspace{2mm}
$$
={\bf 1}_{\{j_3=j_3'\}}{\bf 1}_{\{j_2=j_2'\}}{\bf 1}_{\{j_1=j_1'\}}+
{\bf 1}_{\{i_1=i_2\}} \cdot
{\bf 1}_{\{j_3=j_3'\}}{\bf 1}_{\{j_1=j_2'\}}{\bf 1}_{\{j_2=j_1'\}}+
$$

\vspace{2mm}
$$
+{\bf 1}_{\{i_2=i_3\}} \cdot
{\bf 1}_{\{j_1=j_1'\}}{\bf 1}_{\{j_2=j_3'\}}{\bf 1}_{\{j_3=j_2'\}}+
{\bf 1}_{\{i_1=i_3\}} \cdot
{\bf 1}_{\{j_1=j_3'\}}{\bf 1}_{\{j_2=j_2'\}}{\bf 1}_{\{j_3=j_1'\}}+
$$

\vspace{2mm}
$$
+{\bf 1}_{\{i_1=i_2=i_3\}}\cdot
{\bf 1}_{\{j_2=j_3'\}}{\bf 1}_{\{j_1=j_2'\}}{\bf 1}_{\{j_3=j_1'\}}+
$$

\vspace{2mm}
\begin{equation}
\label{gt23}
+
{\bf 1}_{\{i_1=i_2=i_3\}}\cdot
{\bf 1}_{\{j_1=j_3'\}}{\bf 1}_{\{j_3=j_2'\}}{\bf 1}_{\{j_2=j_1'\}}.
\end{equation}

\vspace{6mm}

From (\ref{gt22}) and (\ref{gt23}) we get

\vspace{1mm}
$$
\left\vert{\sf M}\left\{\sum\limits_{(j_1,j_2)}\sum\limits_{(j_1',j_2')}
\int\limits_t^T\phi_{j_{2}}(t_{2})
\int\limits_t^{t_{2}}\phi_{j_{1}}(t_{1})
d{\bf f}_{t_1}^{(i_1)}d{\bf f}_{t_2}^{(i_2)}
\times\right.\right.
$$

\vspace{3mm}
$$
\left.\left.\times
\int\limits_t^T\phi_{j_{2}'}(t_{2})
\int\limits_t^{t_{2}}\phi_{j_{1}'}(t_{1})
d{\bf f}_{t_1}^{(i_1)}d{\bf f}_{t_2}^{(i_2)}
\right\}\right\vert\le
$$

\vspace{3mm}
$$
\le
{\bf 1}_{\{j_1=j_1'\}}{\bf 1}_{\{j_2=j_2'\}}+
{\bf 1}_{\{j_2=j_1'\}}{\bf 1}_{\{j_1=j_2'\}}=
$$

\vspace{2mm}
$$
=
\sum\limits_{(j_1',j_2')}
{\bf 1}_{\{j_1=j_1'\}}{\bf 1}_{\{j_2=j_2'\}},
$$

\vspace{6mm}

$$
\left\vert{\sf M}\left\{\sum\limits_{(j_1,j_2,j_3)}\sum\limits_{(j_1',j_2',j_3')}
\int\limits_t^T \phi_{j_3}(t_3)
\int\limits_t^{t_{3}}\phi_{j_{2}}(t_{2})
\int\limits_t^{t_{2}}\phi_{j_{1}}(t_{1})
d{\bf f}_{t_1}^{(i_1)}d{\bf f}_{t_2}^{(i_2)}
d{\bf f}_{t_3}^{(i_3)}\times\right.\right.
$$

\vspace{3mm}
$$
\left.\left.\times
\int\limits_t^T \phi_{j_3'}(t_3)
\int\limits_t^{t_{3}}\phi_{j_{2}'}(t_{2})
\int\limits_t^{t_{2}}\phi_{j_{1}'}(t_{1})
d{\bf f}_{t_1}^{(i_1)}d{\bf f}_{t_2}^{(i_2)}
d{\bf f}_{t_3}^{(i_3)}\right\}\right\vert\le
$$

\vspace{3mm}
$$
\le{\bf 1}_{\{j_3=j_3'\}}{\bf 1}_{\{j_2=j_2'\}}{\bf 1}_{\{j_1=j_1'\}}+
{\bf 1}_{\{j_3=j_3'\}}{\bf 1}_{\{j_1=j_2'\}}{\bf 1}_{\{j_2=j_1'\}}+
$$

\vspace{2mm}
$$
+
{\bf 1}_{\{j_1=j_1'\}}{\bf 1}_{\{j_2=j_3'\}}{\bf 1}_{\{j_3=j_2'\}}+
{\bf 1}_{\{j_1=j_3'\}}{\bf 1}_{\{j_2=j_2'\}}{\bf 1}_{\{j_3=j_1'\}}+
$$

\vspace{2mm}
$$
+
{\bf 1}_{\{j_2=j_3'\}}{\bf 1}_{\{j_1=j_2'\}}{\bf 1}_{\{j_3=j_1'\}}+
{\bf 1}_{\{j_1=j_3'\}}{\bf 1}_{\{j_3=j_2'\}}{\bf 1}_{\{j_2=j_1'\}}
=
$$

\vspace{2mm}
$$
=
\sum\limits_{(j_1',j_2',j_3')}
{\bf 1}_{\{j_1=j_1'\}}{\bf 1}_{\{j_2=j_2'\}}{\bf 1}_{\{j_3=j_3'\}},
$$

\vspace{6mm}
\noindent
where we used the relation

\vspace{1mm}
$$
\int\limits_t^T
\phi_{i}(\tau)\phi_{j}(\tau)d\tau={\bf 1}_{\{i=j\}},\ \ \ i,j=0, 1, 2\ldots
$$

\vspace{5mm}

Now consider the case of an arbitrary $k\in \mathbb{N}.$ We have

\vspace{1mm}
$$
{\sf M}\left\{\sum\limits_{(j_1,\ldots,j_k)}\sum\limits_{(j_1',\ldots,j_k')}
\int\limits_t^T \phi_{j_k}(t_k)
\ldots
\int\limits_t^{t_{2}}\phi_{j_{1}}(t_{1})
d{\bf f}_{t_1}^{(i_1)}\ldots
d{\bf f}_{t_k}^{(i_k)}\times\right.
$$

\vspace{3.5mm}
$$
\left.\times
\int\limits_t^T \phi_{j_k'}(t_k)
\ldots
\int\limits_t^{t_{2}}\phi_{j_1'}(t_{1})
d{\bf f}_{t_1}^{(i_1)}\ldots
d{\bf f}_{t_k}^{(i_k)}\right\}=
$$

\vspace{3.5mm}
$$
={\sf M}\left\{\sum\limits_{(j_1,\ldots,j_k)}\sum\limits_{(j_1',\ldots,j_k')}
\int\limits_t^T \phi_{j_k}(t_k)
\ldots
\int\limits_t^{t_{2}}\phi_{j_{1}}(t_{1})
d{\bf f}_{t_1}^{(i_1)}\ldots
d{\bf f}_{t_k}^{(i_k)}\times\right.
$$

\vspace{3.5mm}
$$
\left.\times
\int\limits_t^T \phi_{j_k'}(t_k)
\ldots
\int\limits_t^{t_{2}}\phi_{j_1'}(t_{1})
d{\bf f}_{t_1}^{(i_1')}\ldots
d{\bf f}_{t_k}^{(i_k')}\right\}=
$$

\vspace{3.5mm}
$$
=\sum\limits_{(j_1,\ldots,j_k)}\sum\limits_{(j_1',\ldots,j_k')}
{\bf 1}_{\{i_k=i_k'\}}\ldots {\bf 1}_{\{i_1=i_1'\}}\times
$$

\vspace{3.5mm}
$$
\times
\int\limits_t^T \phi_{j_k}(t_k)\phi_{j_k'}(t_k)
\ldots
\int\limits_t^{t_{2}}\phi_{j_1}(t_1)\phi_{j_1'}(t_1)
dt_1\ldots dt_k=
$$

\vspace{3.5mm}
$$
=\sum\limits_{(j_1',\ldots,j_k')}
{\bf 1}_{\{i_k=i_k'\}}\ldots {\bf 1}_{\{i_1=i_1'\}}\times
$$

\vspace{3.5mm}
$$
\times
\int\limits_t^T \phi_{j_k}(t_k)\phi_{j_k'}(t_k)dt_k
\ldots
\int\limits_t^T\phi_{j_1}(t_1)\phi_{j_1'}(t_1)dt_1=
$$

\vspace{3.5mm}
\begin{equation}
\label{wen100}
=\sum\limits_{(j_1',\ldots,j_k')}
{\bf 1}_{\{i_k=i_k'\}}\ldots {\bf 1}_{\{i_1=i_1'\}}
{\bf 1}_{\{j_k=j_k'\}}\ldots {\bf 1}_{\{j_1=j_1'\}},
\end{equation}

\vspace{6mm}
\noindent
where $(i_1',\ldots,i_k')=(i_1,\ldots,i_k).$
However, if 
$j_r$ swapped  with $j_q$ in the permutation $(j_1,\ldots,j_k)$,
then $i_r$ swapped  with $i_q$ in the permutation
$(i_1,\ldots,i_k)$ and 
if $j_r'$ swapped  with $j_q'$ in the permutation $(j_1',\ldots,j_k')$,
then $i_r'$ swapped  with $i_q'$ in the permutation
$(i_1',\ldots,i_k')$. From (\ref{wen100}) we obtain (\ref{hq11}).
The equality (\ref{tyty1}) is proved.

Note that the formula (\ref{tyty1}) (in the light of the results of 
Sect.~5) can be 
interpreted as a consequence of the orthogonality of two random 
variables that are Hermite polynomials of vector random arguments.

\vspace{3mm}

Applying (\ref{tyty}) and (\ref{tyty1}), we obtain

\vspace{1mm}
$$
{\sf M}\left\{J[\psi^{(k)}]_{T,t}^p\xi[\psi^{(k)}]_{T,t}^{p+1,n}
\right\}=0.
$$

\vspace{5mm}

Due to (\ref{yeee2}), (\ref{yyye1}), and (\ref{yyye}) we can write 

\vspace{2mm}
$$
\xi[\psi^{(k)}]_{T,t}^{p+1,n}=J[\psi^{(k)}]_{T,t}^n-J[\psi^{(k)}]_{T,t}^p,
$$

\vspace{3mm}
$$
\hbox{\vtop{\offinterlineskip\halign{
\hfil#\hfil\cr
{\rm l.i.m.}\cr
$\stackrel{}{{}_{n\to \infty}}$\cr
}} }
\xi[\psi^{(k)}]_{T,t}^{p+1,n}=J[\psi^{(k)}]_{T,t}-J[\psi^{(k)}]_{T,t}^p
\stackrel{\rm def}{=}\xi[\psi^{(k)}]_{T,t}^{p+1}.
$$

\vspace{6mm}

We have

$$
0\le \left|{\sf M}\left\{
\xi[\psi^{(k)}]_{T,t}^{p+1}J[\psi^{(k)}]_{T,t}^p\right\}\right|=
$$

\vspace{4mm}
$$
=
\left|{\sf M}\left\{\left(
\xi[\psi^{(k)}]_{T,t}^{p+1}-
\xi[\psi^{(k)}]_{T,t}^{p+1,n}+\xi[\psi^{(k)}]_{T,t}^{p+1,n}\right)
J[\psi^{(k)}]_{T,t}^p\right\}\right|=
$$

\vspace{4mm}
$$
\le 
\left|{\sf M}\left\{\left(
\xi[\psi^{(k)}]_{T,t}^{p+1}-
\xi[\psi^{(k)}]_{T,t}^{p+1,n}\right)
J[\psi^{(k)}]_{T,t}^p\right\}\right|+
\left|{\sf M}\left\{
\xi[\psi^{(k)}]_{T,t}^{p+1,n}J[\psi^{(k)}]_{T,t}^p\right\}\right|=
$$

\vspace{4mm}
$$
=\left|{\sf M}\left\{\left(
J[\psi^{(k)}]_{T,t}-
J[\psi^{(k)}]_{T,t}^{n}\right)
J[\psi^{(k)}]_{T,t}^p\right\}\right|\le
$$

\vspace{4mm}
$$
\le \sqrt{{\sf M}\left\{\left(
J[\psi^{(k)}]_{T,t}-
J[\psi^{(k)}]_{T,t}^{n}\right)^2\right\}}
\sqrt{{\sf M}\left\{\left(
J[\psi^{(k)}]_{T,t}^p\right)^2\right\}}\le
$$

\vspace{4mm}
$$
\le \sqrt{{\sf M}\left\{\left(
J[\psi^{(k)}]_{T,t}-
J[\psi^{(k)}]_{T,t}^{n}\right)^2\right\}}\times
$$

\vspace{4mm}
$$
\times
\left(\sqrt{{\sf M}\left\{\left(
J[\psi^{(k)}]_{T,t}^p - J[\psi^{(k)}]_{T,t}\right)^2\right\}}
+ \sqrt{{\sf M}\left\{\left(
J[\psi^{(k)}]_{T,t}\right)^2\right\}}\right)\le
$$

\vspace{4mm}
\begin{equation}
\label{rrre}
\le K \sqrt{{\sf M}\left\{\left(
J[\psi^{(k)}]_{T,t}-
J[\psi^{(k)}]_{T,t}^{n}\right)^2\right\}} \to 0\ \ \ {\rm if}\ \ \ n\to\infty,
\end{equation}

\vspace{6mm}
\noindent
where $K$ is a constant.

From (\ref{rrre}) it follows that

$$
{\sf M}\left\{
\xi[\psi^{(k)}]_{T,t}^{p+1}J[\psi^{(k)}]_{T,t}^p\right\}=0
$$

\vspace{3mm}
\noindent
or

$$
{\sf M}\left\{
\left(J[\psi^{(k)}]_{T,t}-J[\psi^{(k)}]_{T,t}^p\right)
J[\psi^{(k)}]_{T,t}^p\right\}=0.
$$

\vspace{6mm}

The last equality means that

\vspace{1mm}
\begin{equation}
\label{yyyw}
{\sf M}\left\{
J[\psi^{(k)}]_{T,t}J[\psi^{(k)}]_{T,t}^p\right\}=
{\sf M}\left\{
\left(J[\psi^{(k)}]_{T,t}^p\right)^2\right\}.
\end{equation}

\vspace{5mm}

Taking into account (\ref{yyyw}), we obtain

\vspace{2mm}
$$
{\sf M}\left\{\left(J[\psi^{(k)}]_{T,t}-
J[\psi^{(k)}]_{T,t}^p\right)^2\right\}=
{\sf M}\left\{\left(J[\psi^{(k)}]_{T,t}\right)^2\right\}+
$$

\vspace{3mm}
$$
+
{\sf M}\left\{\left(J[\psi^{(k)}]_{T,t}^p\right)^2\right\}
-2{\sf M}\left\{J[\psi^{(k)}]_{T,t}J[\psi^{(k)}]_{T,t}^p\right\}=
{\sf M}\left\{\left(J[\psi^{(k)}]_{T,t}\right)^2\right\}-
$$

\vspace{4mm}
$$
-
{\sf M}\left\{J[\psi^{(k)}]_{T,t}J[\psi^{(k)}]_{T,t}^p\right\}=
$$

\vspace{2mm}
\begin{equation}
\label{tttr}
=\int\limits_{[t,T]^k} K^2(t_1,\ldots,t_k)
dt_1\ldots dt_k - 
{\sf M}\left\{J[\psi^{(k)}]_{T,t}J[\psi^{(k)}]_{T,t}^p\right\}.
\end{equation}

\vspace{6mm}

Let us consider the value

\vspace{1mm}
$$
{\sf M}\left\{J[\psi^{(k)}]_{T,t}J[\psi^{(k)}]_{T,t}^p\right\}.
$$

\vspace{5mm}

Using (\ref{yeee2}) and (\ref{ttt2}), we get

\vspace{1mm}
\begin{equation}
\label{z9}
J[\psi^{(k)}]_{T,t}^p=
\sum_{j_1=0}^{p}\ldots\sum_{j_k=0}^{p}
C_{j_k\ldots j_1}
\sum\limits_{(j_1,\ldots,j_k)}
\int\limits_t^T \phi_{j_k}(t_k)
\ldots
\int\limits_t^{t_{2}}\phi_{j_{1}}(t_{1})
d{\bf f}_{t_1}^{(i_1)}\ldots
d{\bf f}_{t_k}^{(i_k)}.
\end{equation}

\vspace{6mm}

After substituting (\ref{z9}) into (\ref{tttr}), we obtain

\vspace{2mm}
$$
{\sf M}\left\{\left(J[\psi^{(k)}]_{T,t}-
J[\psi^{(k)}]_{T,t}^p\right)^2\right\}
= \int\limits_{[t,T]^k} K^2(t_1,\ldots,t_k)
dt_1\ldots dt_k - 
$$

\vspace{3mm}
$$
- \sum_{j_1=0}^{p}\ldots\sum_{j_k=0}^{p}
C_{j_k\ldots j_1}
{\sf M}\left\{J[\psi^{(k)}]_{T,t}
\sum\limits_{(j_1,\ldots,j_k)}
\int\limits_t^T \phi_{j_k}(t_k)
\ldots
\int\limits_t^{t_{2}}\phi_{j_{1}}(t_{1})
d{\bf f}_{t_1}^{(i_1)}\ldots
d{\bf f}_{t_k}^{(i_k)}\right\}.
$$

\vspace{7mm}

Theorem 2 is proved.

\vspace{5mm}

\section{Exact Calculation of the Mean-Square Approximation Error
for the Cases $k=1,\ldots,5$}

\vspace{2mm}

Let us denote

$$
{\sf M}\left\{\left(J[\psi^{(k)}]_{T,t}-
J[\psi^{(k)}]_{T,t}^{p}\right)^2\right\}\stackrel{{\rm def}}
{=}E_k^{p},
$$

\vspace{3mm}
$$
\int\limits_{[t,T]^k} K^2(t_1,\ldots,t_k)
dt_1\ldots dt_k\stackrel{{\rm def}}{=}I_k.
$$

\vspace{3mm}

\subsection{The Case $k=1$}

\vspace{1mm}

For this case from Theorem 2 we obtain

\vspace{1mm}
$$
E_1^p
=I_1
-\sum_{j_1=0}^p
C_{j_1}^2.
$$

\vspace{2mm}

\subsection{The Case $k=2$}

For this case from Theorem 2 we have

\vspace{4mm}

(I).\ $i_1\ne i_2$:
\begin{equation}
\label{ee1}
E_2^p
=I_2
-\sum_{j_1,j_2=0}^p
C_{j_2j_1}^2,
\end{equation}

\vspace{2mm}

(II).\ $i_1=i_2:$
$$
E_2^p
=I_2
-\sum_{j_1,j_2=0}^p
C_{j_2j_1}^2-
\sum_{j_1,j_2=0}^p
C_{j_2j_1}C_{j_1j_2}.
$$

\vspace{5mm}

{\bf Example 1.} Let us consider the following
iterated Ito stochastic integral from the stochastic Taylor--Ito expansion 
\cite{Zapad-1}-\cite{Zapad-3} (also see \cite{arxiv-24} and 
\cite{1}-\cite{12})

\begin{equation}
\label{k1000}
I_{(00)T,t}^{(i_1i_2)}
=\int\limits_t^T\int\limits_t^{t_{2}}
d{\bf f}_{t_1}^{(i_1)}
d{\bf f}_{t_2}^{(i_2)},
\end{equation}

\vspace{3mm}
\noindent
where $i_1, i_2=1,\ldots,m.$

The approximation based on the expansion (\ref{leto5001}) for 
the integral (\ref{k1000}) (the case
of Legendre poly\-no\-mi\-als) has the following form \cite{1}-\cite{OK}

\vspace{1mm}
\begin{equation}
\label{4004}
I_{(00)T,t}^{(i_1 i_2)p}=
\frac{T-t}{2}\left(\zeta_0^{(i_1)}\zeta_0^{(i_2)}+\sum_{i=1}^{p}
\frac{1}{\sqrt{4i^2-1}}\left(
\zeta_{i-1}^{(i_1)}\zeta_{i}^{(i_2)}-
\zeta_i^{(i_1)}\zeta_{i-1}^{(i_2)}\right)-{\bf 1}_{\{i_1=i_2\}}\right).
\end{equation}

\vspace{5mm}

It should be noted that the formula (\ref{4004}) has been derived 
for the first time in \cite{old-art-1} (1997), \cite{old-art-2} (1998) 
with using the another approach, which
was developed in \cite{arxiv-6}.

Applying (\ref{ee1}), we obtain \cite{1}-\cite{OK} (also
see \cite{old-art-1} (1997), \cite{old-art-2} (1998))

\vspace{1mm}
\begin{equation}
\label{4007}
{\sf M}\biggl\{\left(I_{(00)T,t}^{(i_1 i_2)}-
I_{(00)T,t}^{(i_1 i_2)p}
\right)^2\biggr\}
=\frac{(T-t)^2}{2}\left(\frac{1}{2}-\sum_{i=1}^p
\frac{1}{4i^2-1}\right)\ \ \ (i_1\ne i_2).
\end{equation}

\vspace{2mm}

\subsection{The Case $k=3$}

For the case $k=3$ from Theorem 2 we obtain

\vspace{4mm}

(I).\ $i_1\ne i_2, i_1\ne i_3, i_2\ne i_3:$
\begin{equation}
\label{dest900}
E_3^p=I_3
-\sum_{j_1,j_2,j_3=0}^p C_{j_3j_2j_1}^2,
\end{equation}

\vspace{2mm}

(II).\ $i_1=i_2=i_3:$
$$
E_3^p = I_3 - \sum_{j_1,j_2,j_3=0}^{p}
C_{j_3j_2j_1}\Biggl(\sum\limits_{(j_1,j_2,j_3)}
C_{j_3j_2j_1}\Biggr),
$$

\vspace{2mm}

(III).1.\ $i_1=i_2\ne i_3:$
\begin{equation}
\label{qq1}
E_3^p=I_3
-\sum_{j_1,j_2,j_3=0}^p C_{j_3j_2j_1}^2-
\sum_{j_1,j_2,j_3=0}^p C_{j_3j_1j_2}C_{j_3j_2j_1},
\end{equation}

\vspace{2mm}

(III).2.\ $i_1\ne i_2=i_3:$
\begin{equation}
\label{dest1000}
E_3^p=I_3-
\sum_{j_1,j_2,j_3=0}^p C_{j_3j_2j_1}^2-
\sum_{j_1,j_2,j_3=0}^p C_{j_2j_3j_1}C_{j_3j_2j_1},
\end{equation}

\vspace{2mm}

(III).3.\ $i_1=i_3\ne i_2:$
\begin{equation}
\label{dest1001}
E_3^p=I_3
-\sum_{j_1,j_2,j_3=0}^p C_{j_3j_2j_1}^2-
\sum_{j_1,j_2,j_3=0}^p C_{j_3j_2j_1}C_{j_1j_2j_3}.
\end{equation}

\vspace{5mm}

It should be noted that the formulas from the above cases (I), (III).1,
(III).2, (III).3 have been derived in \cite{1}-\cite{12}
by direct calculation.

\vspace{2mm}

{\bf Example 2.} Let us consider the following
iterated Ito stochastic integral from the stochastic 
Taylor--Ito expansion \cite{Zapad-1}-\cite{Zapad-3} 
(also see \cite{arxiv-24} and \cite{1}-\cite{12})

\begin{equation}
\label{k1001}
I_{(000)T,t}^{(i_1i_2i_3)}
=\int\limits_t^T\int\limits_t^{t_{3}}\int\limits_t^{t_{2}}
d{\bf f}_{t_1}^{(i_1)}
d{\bf f}_{t_2}^{(i_2)}
d{\bf f}_{t_3}^{(i_3)},
\end{equation}

\vspace{3mm}
\noindent
where $i_1, i_2, i_3=1,\ldots,m.$

The approximation based on the expansion (\ref{leto5002}) for the 
integral (\ref{k1001}) (the case
of Legendre polynomials and $p_1=p_2=p_3=p$) has the following form 
\cite{1}-\cite{OK}

\vspace{2mm}
$$
I_{(000)T,t}^{(i_1i_2i_3)p}
=\sum_{j_1,j_2,j_3=0}^{p}
C_{j_3j_2j_1}\Biggl(
\zeta_{j_1}^{(i_1)}\zeta_{j_2}^{(i_2)}\zeta_{j_3}^{(i_3)}
-{\bf 1}_{\{i_1=i_2\}}
{\bf 1}_{\{j_1=j_2\}}
\zeta_{j_3}^{(i_3)}-
\Biggr.
$$

\vspace{2mm}
\begin{equation}
\label{sad001}
\Biggl.
-{\bf 1}_{\{i_2=i_3\}}
{\bf 1}_{\{j_2=j_3\}}
\zeta_{j_1}^{(i_1)}-
{\bf 1}_{\{i_1=i_3\}}
{\bf 1}_{\{j_1=j_3\}}
\zeta_{j_2}^{(i_2)}\Biggr),
\end{equation}

\vspace{5mm}
\noindent
where
\begin{equation}
\label{w1}
C_{j_3j_2j_1}
=\frac{\sqrt{(2j_1+1)(2j_2+1)(2j_3+1)}}{8}(T-t)^{3/2}\bar
C_{j_3j_2j_1},
\end{equation}

\vspace{1mm}
$$
\bar C_{j_3j_2j_1}=\int\limits_{-1}^{1}P_{j_3}(z)
\int\limits_{-1}^{z}P_{j_2}(y)
\int\limits_{-1}^{y}
P_{j_1}(x)dx dy dz,
$$

\vspace{5mm}
\noindent
where $P_i(x)$ is the Legendre polynomial $(i= 0, 1, 2,\ldots ).$

For example, using (\ref{qq1}) we obtain

\vspace{1mm}
$$
{\sf M}\left\{\left(
I_{(000)T,t}^{(i_1i_2 i_3)}-
I_{(000)T,t}^{(i_1i_2 i_3)p}\right)^2\right\}=
$$

\vspace{1mm}
$$
=
\frac{(T-t)^{3}}{6}-\sum_{j_1,j_2,j_3=0}^{p}
C_{j_3j_2j_1}^2
-\sum_{j_1,j_2,J_3=0}^{p}
C_{j_3j_1j_2}C_{j_3j_2j_1},
$$

\vspace{4mm}
\noindent
where $i_1=i_2\ne i_3$.

As mentioned in \cite{arxiv-4} (also see \cite{1}-\cite{12}), the exact 
values of coefficients 
$\bar C_{j_3j_2j_1}$ when $j_1,j_2,j_3=0, 1,\ldots,p$
can be calculated using DERIVE (computer
algebra system). In 
\cite{arxiv-4} (also see \cite{1}-\cite{12})
we can find several tables with exactly calculated
Fourier--Legendre coefficients
for approximations of iterated Ito stochastic integrals of 
multiplicities 1 to 5.
In addition, in \cite{Kuz-Kuz}, \cite{Mikh-1}, a database was obtained with 270,000 exactly
calculated Fourier--Legendre coefficients 
for approximations of iterated Ito and Stratonovich  stochastic integrals of 
multiplicities 1 to 6.

For the case $i_1=i_2=i_3$ we can use 
the following formula \cite{Zapad-1}-\cite{Zapad-3}

\vspace{1mm}
$$
I_{(000)T,t}^{(i_1 i_1 i_1)}=\frac{1}{6}(T-t)^{3/2}\left(
\left(\zeta_0^{(i_1)}\right)^3-3
\zeta_0^{(i_1)}\right)
$$

\vspace{4mm}
\noindent
w.~p.~1.

\vspace{5mm}

\subsection{The Case $k=4$}

For this case from Theorem 2 we obtain

\vspace{4mm}

(I).\ $i_1,\ldots,i_4$ are pairwise different:
$$
E_4^p= I_4 - \sum_{j_1,\ldots,j_4=0}^{p}C_{j_4\ldots j_1}^2,
$$

\vspace{4mm}

(II).\ $i_1=i_2=i_3=i_4$:
$$
E_4^p = I_4 -
 \sum_{j_1,\ldots,j_4=0}^{p}
C_{j_4\ldots j_1}\Biggl(\sum\limits_{(j_1,\ldots,j_4)}
C_{j_4\ldots j_1}\Biggr),
$$

\vspace{4mm}

(III).1.\ $i_1=i_2\ne i_3, i_4;\ i_3\ne i_4:$
$$
E^p_4 = I_4 - \sum_{j_1,\ldots,j_4=0}^{p}
C_{j_4\ldots j_1}\Biggl(\sum\limits_{(j_1,j_2)}
C_{j_4\ldots j_1}\Biggr),
$$

\vspace{4mm}

(III).2.\ $i_1=i_3\ne i_2, i_4;\ i_2\ne i_4:$
$$
E^p_4 = I_4 - \sum_{j_1,\ldots,j_4=0}^{p}
C_{j_4\ldots j_1}\Biggl(\sum\limits_{(j_1,j_3)}
C_{j_4\ldots j_1}\Biggr),
$$

\vspace{4mm}

(III).3.\ $i_1=i_4\ne i_2, i_3;\ i_2\ne i_3:$
$$
E^p_4 = I_4 - \sum_{j_1,\ldots,j_4=0}^{p}
C_{j_4\ldots j_1}\Biggl(\sum\limits_{(j_1,j_4)}
C_{j_4\ldots j_1}\Biggr),
$$

\vspace{4mm}

(III).4.\ $i_2=i_3\ne i_1, i_4;\ i_1\ne i_4:$
$$
E^p_4 = I_4 - \sum_{j_1,\ldots,j_4=0}^{p}
C_{j_4\ldots j_1}\Biggl(\sum\limits_{(j_2,j_3)}
C_{j_4\ldots j_1}\Biggr),
$$

\vspace{4mm}

(III).5.\ $i_2=i_4\ne i_1, i_3;\ i_1\ne i_3:$
$$
E^p_4 = I_4 - \sum_{j_1,\ldots,j_4=0}^{p}
C_{j_4\ldots j_1}\Biggl(\sum\limits_{(j_2,j_4)}
C_{j_4\ldots j_1}\Biggr),
$$

\vspace{4mm}

(III).6.\ $i_3=i_4\ne i_1, i_2;\ i_1\ne i_2:$
$$
E^p_4 = I_4 - \sum_{j_1,\ldots,j_4=0}^{p}
C_{j_4\ldots j_1}\Biggl(\sum\limits_{(j_3,j_4)}
C_{j_4\ldots j_1}\Biggr),
$$

\vspace{4mm}

(IV).1.\ $i_1=i_2=i_3\ne i_4$:
$$
E_4^p = I_4 -
 \sum_{j_1,\ldots,j_4=0}^{p}
C_{j_4\ldots j_1}\Biggl(\sum\limits_{(j_1,j_2,j_3)}
C_{j_4\ldots j_1}\Biggr),
$$

\vspace{4mm}

(IV).2.\ $i_2=i_3=i_4\ne i_1$:
$$
E_4^p = I_4 -
 \sum_{j_1,\ldots,j_4=0}^{p}
C_{j_4\ldots j_1}\Biggl(\sum\limits_{(j_2,j_3,j_4)}
C_{j_4\ldots j_1}\Biggr),
$$

\vspace{4mm}

(IV).3.\ $i_1=i_2=i_4\ne i_3$:
$$
E_4^p = I_4 -
 \sum_{j_1,\ldots,j_4=0}^{p}
C_{j_4\ldots j_1}\Biggl(\sum\limits_{(j_1,j_2,j_4)}
C_{j_4\ldots j_1}\Biggr),
$$

\vspace{4mm}

(IV).4.\ $i_1=i_3=i_4\ne i_2$:
$$
E_4^p = I_4 -
 \sum_{j_1,\ldots,j_4=0}^{p}
C_{j_4\ldots j_1}\Biggl(\sum\limits_{(j_1,j_3,j_4)}
C_{j_4\ldots j_1}\Biggr),
$$

\vspace{4mm}

(V).1.\ $i_1=i_2\ne i_3=i_4$:
$$
E^p_4 = I_4 - \sum_{j_1,\ldots,j_4=0}^{p}
C_{j_4\ldots j_1}\Biggl(\sum\limits_{(j_1,j_2)}\Biggl(
\sum\limits_{(j_3,j_4)}
C_{j_4\ldots j_1}\Biggr)\Biggr),
$$

\vspace{4mm}

(V).2.\ $i_1=i_3\ne i_2=i_4$:
$$
E^p_4 = I_4 - \sum_{j_1,\ldots,j_4=0}^{p}
C_{j_4\ldots j_1}\Biggl(\sum\limits_{(j_1,j_3)}\Biggl(
\sum\limits_{(j_2,j_4)}
C_{j_4\ldots j_1}\Biggr)\Biggr),
$$

\vspace{4mm}

(V).3.\ $i_1=i_4\ne i_2=i_3$:
$$
E^p_4 = I_4 - \sum_{j_1,\ldots,j_4=0}^{p}
C_{j_4\ldots j_1}\Biggl(\sum\limits_{(j_1,j_4)}\Biggl(
\sum\limits_{(j_2,j_3)}
C_{j_4\ldots j_1}\Biggr)\Biggr).
$$

\vspace{5mm}

\subsection{The Case $k=5$}

For the case $k=5$ from Theorem 2 we obtain

\vspace{4mm}

(I).\ $i_1,\ldots,i_5$ are pairwise different:
$$
E_5^p= I_5 - \sum_{j_1,\ldots,j_5=0}^{p}C_{j_5\ldots j_1}^2,
$$

\vspace{4mm}

(II).\ $i_1=i_2=i_3=i_4=i_5$:
$$
E_5^p = I_5 - \sum_{j_1,\ldots,j_5=0}^{p}
C_{j_5\ldots j_1}\Biggl(\sum\limits_{(j_1,\ldots,j_5)}
C_{j_5\ldots j_1}\Biggr),
$$

\vspace{7mm}

(III).1.\ $i_1=i_2\ne i_3,i_4,i_5$\ ($i_3,i_4,i_5$ are pairwise different):
$$
E^p_5 = I_5 - \sum_{j_1,\ldots,j_5=0}^{p}
C_{j_5\ldots j_1}\Biggl(\sum\limits_{(j_1,j_2)}
C_{j_5\ldots j_1}\Biggr),
$$

\vspace{7mm}

(III).2.\ $i_1=i_3\ne i_2,i_4,i_5$\ ($i_2,i_4,i_5$ are pairwise different):
$$
E^p_5 = I_5 - \sum_{j_1,\ldots,j_5=0}^{p}
C_{j_5\ldots j_1}\Biggl(\sum\limits_{(j_1,j_3)}
C_{j_5\ldots j_1}\Biggr),
$$

\vspace{7mm}

(III).3.\ $i_1=i_4\ne i_2,i_3,i_5$\ ($i_2,i_3,i_5$ are pairwise different):
$$
E^p_5 = I_5 - \sum_{j_1,\ldots,j_5=0}^{p}
C_{j_5\ldots j_1}\Biggl(\sum\limits_{(j_1,j_4)}
C_{j_5\ldots j_1}\Biggr),
$$

\vspace{7mm}

(III).4.\ $i_1=i_5\ne i_2,i_3,i_4$\ ($i_2,i_3,i_4$  are pairwise different):
$$
E^p_5 = I_5 - \sum_{j_1,\ldots,j_5=0}^{p}
C_{j_5\ldots j_1}\Biggl(\sum\limits_{(j_1,j_5)}
C_{j_5\ldots j_1}\Biggr),
$$

\vspace{7mm}

(III).5.\ $i_2=i_3\ne i_1,i_4,i_5$\ ($i_1,i_4,i_5$ are pairwise different):
$$
E^p_5 = I_5 - \sum_{j_1,\ldots,j_5=0}^{p}
C_{j_5\ldots j_1}\Biggl(\sum\limits_{(j_2,j_3)}
C_{j_5\ldots j_1}\Biggr),
$$

\vspace{7mm}

(III).6.\ $i_2=i_4\ne i_1,i_3,i_5$\ ($i_1,i_3,i_5$ are pairwise different):
$$
E^p_5 = I_5 - \sum_{j_1,\ldots,j_5=0}^{p}
C_{j_5\ldots j_1}\Biggl(\sum\limits_{(j_2,j_4)}
C_{j_5\ldots j_1}\Biggr),
$$

\vspace{7mm}

(III).7.\ $i_2=i_5\ne i_1,i_3,i_4$\ ($i_1,i_3,i_4$ are pairwise different):
$$
E^p_5 = I_5 - \sum_{j_1,\ldots,j_5=0}^{p}
C_{j_5\ldots j_1}\Biggl(\sum\limits_{(j_2,j_5)}
C_{j_5\ldots j_1}\Biggr),
$$

\vspace{7mm}

(III).8.\ $i_3=i_4\ne i_1,i_2,i_5$\ ($i_1,i_2,i_5$  are pairwise different):
$$
E^p_5 = I_5 - \sum_{j_1,\ldots,j_5=0}^{p}
C_{j_5\ldots j_1}\Biggl(\sum\limits_{(j_3,j_4)}
C_{j_5\ldots j_1}\Biggr),
$$

\vspace{7mm}

(III).9.\ $i_3=i_5\ne i_1,i_2,i_4$\ ($i_1,i_2,i_4$  are pairwise different):
$$
E^p_5 = I_5 - \sum_{j_1,\ldots,j_5=0}^{p}
C_{j_5\ldots j_1}\Biggl(\sum\limits_{(j_3,j_5)}
C_{j_5\ldots j_1}\Biggr),
$$

\vspace{7mm}

(III).10.\ $i_4=i_5\ne i_1,i_2,i_3$\ ($i_1,i_2,i_3$  are pairwise different):
$$
E^p_5 = I_5 - \sum_{j_1,\ldots,j_5=0}^{p}
C_{j_5\ldots j_1}\Biggl(\sum\limits_{(j_4,j_5)}
C_{j_5\ldots j_1}\Biggr),
$$

\vspace{7mm}

(IV).1.\ $i_1=i_2=i_3\ne i_4, i_5$\ $(i_4\ne i_5$):
$$
E^p_5 = I_5 - \sum_{j_1,\ldots,j_5=0}^{p}
C_{j_5\ldots j_1}\Biggl(\sum\limits_{(j_1,j_2,j_3)}
C_{j_5\ldots j_1}\Biggr),
$$

\vspace{7mm}

(IV).2.\ $i_1=i_2=i_4\ne i_3, i_5$\  $(i_3\ne i_5$):
$$
E^p_5 = I_5 - \sum_{j_1,\ldots,j_5=0}^{p}
C_{j_5\ldots j_1}\Biggl(\sum\limits_{(j_1,j_2,j_4)}
C_{j_5\ldots j_1}\Biggr),
$$

\vspace{7mm}

(IV).3.\ $i_1=i_2=i_5\ne i_3, i_4$\  $(i_3\ne i_4$):
$$
E^p_5 = I_5 - \sum_{j_1,\ldots,j_5=0}^{p}
C_{j_5\ldots j_1}\Biggl(\sum\limits_{(j_1,j_2,j_5)}
C_{j_5\ldots j_1}\Biggr),
$$

\vspace{7mm}

(IV).4.\ $i_2=i_3=i_4\ne i_1, i_5$\  $(i_1\ne i_5$):
$$
E^p_5 = I_5 - \sum_{j_1,\ldots,j_5=0}^{p}
C_{j_5\ldots j_1}\Biggl(\sum\limits_{(j_2,j_3,j_4)}
C_{j_5\ldots j_1}\Biggr),
$$

\vspace{7mm}
(IV).5.\ $i_2=i_3=i_5\ne i_1, i_4$\  $(i_1\ne i_4$):
$$
E^p_5 = I_5 - \sum_{j_1,\ldots,j_5=0}^{p}
C_{j_5\ldots j_1}\Biggl(\sum\limits_{(j_2,j_3,j_5)}
C_{j_5\ldots j_1}\Biggr),
$$

\vspace{7mm}

(IV).6.\ $i_2=i_4=i_5\ne i_1, i_3$\  $(i_1\ne i_3$):
$$
E^p_5 = I_5 - \sum_{j_1,\ldots,j_5=0}^{p}
C_{j_5\ldots j_1}\Biggl(\sum\limits_{(j_2,j_4,j_5)}
C_{j_5\ldots j_1}\Biggr),
$$

\vspace{7mm}

(IV).7.\ $i_3=i_4=i_5\ne i_1, i_2$\  $(i_1\ne i_2$):
$$
E^p_5 = I_5 - \sum_{j_1,\ldots,j_5=0}^{p}
C_{j_5\ldots j_1}\Biggl(\sum\limits_{(j_3,j_4,j_5)}
C_{j_5\ldots j_1}\Biggr),
$$

\vspace{7mm}

(IV).8.\ $i_1=i_3=i_5\ne i_2, i_4$\  $(i_2\ne i_4$):
$$
E^p_5 = I_5 - \sum_{j_1,\ldots,j_5=0}^{p}
C_{j_5\ldots j_1}\Biggl(\sum\limits_{(j_1,j_3,j_5)}
C_{j_5\ldots j_1}\Biggr),
$$

\vspace{7mm}

(IV).9.\ $i_1=i_3=i_4\ne i_2, i_5$\  $(i_2\ne i_5$):
$$
E^p_5 = I_5 - \sum_{j_1,\ldots,j_5=0}^{p}
C_{j_5\ldots j_1}\Biggl(\sum\limits_{(j_1,j_3,j_4)}
C_{j_5\ldots j_1}\Biggr),
$$

\vspace{7mm}

(IV).10.\ $i_1=i_4=i_5\ne i_2, i_3$\  $(i_2\ne i_3$):
$$
E^p_5 = I_5 - \sum_{j_1,\ldots,j_5=0}^{p}
C_{j_5\ldots j_1}\Biggl(\sum\limits_{(j_1,j_4,j_5)}
C_{j_5\ldots j_1}\Biggr),
$$

\vspace{7mm}

(V).1.\ $i_1=i_2=i_3=i_4\ne i_5$:
$$
E^p_5 = I_5 - \sum_{j_1,\ldots,j_5=0}^{p}
C_{j_5\ldots j_1}\Biggl(\sum\limits_{(j_1,j_2,j_3,j_4)}
C_{j_5\ldots j_1}\Biggr),
$$

\vspace{7mm}

(V).2.\ $i_1=i_2=i_3=i_5\ne i_4$:
$$
E^p_5 = I_5 - \sum_{j_1,\ldots,j_5=0}^{p}
C_{j_5\ldots j_1}\Biggl(\sum\limits_{(j_1,j_2,j_3,j_5)}
C_{j_5\ldots j_1}\Biggr),
$$

\vspace{7mm}

(V).3.\ $i_1=i_2=i_4=i_5\ne i_3$:
$$
E^p_5 = I_5 - \sum_{j_1,\ldots,j_5=0}^{p}
C_{j_5\ldots j_1}\Biggl(\sum\limits_{(j_1,j_2,j_4,j_5)}
C_{j_5\ldots j_1}\Biggr),
$$

\vspace{7mm}

(V).4.\ $i_1=i_3=i_4=i_5\ne i_2$:
$$
E^p_5 = I_5 - \sum_{j_1,\ldots,j_5=0}^{p}
C_{j_5\ldots j_1}\Biggl(\sum\limits_{(j_1,j_3,j_4,j_5)}
C_{j_5\ldots j_1}\Biggr),
$$

\vspace{7mm}

(V).5.\ $i_2=i_3=i_4=i_5\ne i_1$:
$$
E^p_5 = I_5 - \sum_{j_1,\ldots,j_5=0}^{p}
C_{j_5\ldots j_1}\Biggl(\sum\limits_{(j_2,j_3,j_4,j_5)}
C_{j_5\ldots j_1}\Biggr),
$$

\vspace{7mm}

(VI).1.\ $i_5\ne i_1=i_2\ne i_3=i_4\ne i_5$:
$$
E^p_5 = I_5 - \sum_{j_1,\ldots,j_5=0}^{p}
C_{j_5\ldots j_1}\Biggl(\sum\limits_{(j_1,j_2)}\Biggl(
\sum\limits_{(j_3,j_4)}
C_{j_5\ldots j_1}\Biggr)\Biggr),
$$

\vspace{7mm}

(VI).2.\ $i_5\ne i_1=i_3\ne i_2=i_4\ne i_5$:
$$
E^p_5 = I_5 - \sum_{j_1,\ldots,j_5=0}^{p}
C_{j_5\ldots j_1}\Biggl(\sum\limits_{(j_1,j_3)}\Biggl(
\sum\limits_{(j_2,j_4)}
C_{j_5\ldots j_1}\Biggr)\Biggr),
$$

\vspace{7mm}

(VI).3.\ $i_5\ne i_1=i_4\ne i_2=i_3\ne i_5$:
$$
E^p_5 = I_5 - \sum_{j_1,\ldots,j_5=0}^{p}
C_{j_5\ldots j_1}\Biggl(\sum\limits_{(j_1,j_4)}\Biggl(
\sum\limits_{(j_2,j_3)}
C_{j_5\ldots j_1}\Biggr)\Biggr),
$$

\vspace{7mm}

(VI).4.\ $i_4\ne i_1=i_2\ne i_3=i_5\ne i_4$:
$$
E^p_5 = I_5 - \sum_{j_1,\ldots,j_5=0}^{p}
C_{j_5\ldots j_1}\Biggl(\sum\limits_{(j_1,j_2)}\Biggl(
\sum\limits_{(j_3,j_5)}
C_{j_5\ldots j_1}\Biggr)\Biggr),
$$

\vspace{7mm}

(VI).5.\ $i_4\ne i_1=i_5\ne i_2=i_3\ne i_4$:
$$
E^p_5 = I_5 - \sum_{j_1,\ldots,j_5=0}^{p}
C_{j_5\ldots j_1}\Biggl(\sum\limits_{(j_1,j_5)}\Biggl(
\sum\limits_{(j_2,j_3)}
C_{j_5\ldots j_1}\Biggr)\Biggr),
$$

\vspace{7mm}

(VI).6.\ $i_4\ne i_2=i_5\ne i_1=i_3\ne i_4$:
$$
E^p_5 = I_5 - \sum_{j_1,\ldots,j_5=0}^{p}
C_{j_5\ldots j_1}\Biggl(\sum\limits_{(j_2,j_5)}\Biggl(
\sum\limits_{(j_1,j_3)}
C_{j_5\ldots j_1}\Biggr)\Biggr),
$$

\vspace{7mm}

(VI).7.\ $i_3\ne i_2=i_5\ne i_1=i_4\ne i_3$:
$$
E^p_5 = I_5 - \sum_{j_1,\ldots,j_5=0}^{p}
C_{j_5\ldots j_1}\Biggl(\sum\limits_{(j_2,j_5)}\Biggl(
\sum\limits_{(j_1,j_4)}
C_{j_5\ldots j_1}\Biggr)\Biggr),
$$

\vspace{7mm}

(VI).8.\ $i_3\ne i_1=i_2\ne i_4=i_5\ne i_3$:
$$
E^p_5 = I_5 - \sum_{j_1,\ldots,j_5=0}^{p}
C_{j_5\ldots j_1}\Biggl(\sum\limits_{(j_1,j_2)}\Biggl(
\sum\limits_{(j_4,j_5)}
C_{j_5\ldots j_1}\Biggr)\Biggr),
$$

\vspace{7mm}

(VI).9.\ $i_3\ne i_2=i_4\ne i_1=i_5\ne i_3$:
$$
E^p_5 = I_5 - \sum_{j_1,\ldots,j_5=0}^{p}
C_{j_5\ldots j_1}\Biggl(\sum\limits_{(j_2,j_4)}\Biggl(
\sum\limits_{(j_1,j_5)}
C_{j_5\ldots j_1}\Biggr)\Biggr),
$$

\vspace{7mm}

(VI).10.\ $i_2\ne i_1=i_4\ne i_3=i_5\ne i_2$:
$$
E^p_5 = I_5 - \sum_{j_1,\ldots,j_5=0}^{p}
C_{j_5\ldots j_1}\Biggl(\sum\limits_{(j_1,j_4)}\Biggl(
\sum\limits_{(j_3,j_5)}
C_{j_5\ldots j_1}\Biggr)\Biggr),
$$

\vspace{7mm}

(VI).11.\ $i_2\ne i_1=i_3\ne i_4=i_5\ne i_2$:
$$
E^p_5 = I_5 - \sum_{j_1,\ldots,j_5=0}^{p}
C_{j_5\ldots j_1}\Biggl(\sum\limits_{(j_1,j_3)}\Biggl(
\sum\limits_{(j_4,j_5)}
C_{j_5\ldots j_1}\Biggr)\Biggr),
$$

\vspace{7mm}

(VI).12.\ $i_2\ne i_1=i_5\ne i_3=i_4\ne i_2$:
$$
E^p_5 = I_5 - \sum_{j_1,\ldots,j_5=0}^{p}
C_{j_5\ldots j_1}\Biggl(\sum\limits_{(j_1,j_5)}\Biggl(
\sum\limits_{(j_3,j_4)}
C_{j_5\ldots j_1}\Biggr)\Biggr),
$$

\vspace{7mm}

(VI).13.\ $i_1\ne i_2=i_3\ne i_4=i_5\ne i_1$:
$$
E^p_5 = I_5 - \sum_{j_1,\ldots,j_5=0}^{p}
C_{j_5\ldots j_1}\Biggl(\sum\limits_{(j_2,j_3)}\Biggl(
\sum\limits_{(j_4,j_5)}
C_{j_5\ldots j_1}\Biggr)\Biggr),
$$

\vspace{7mm}

(VI).14.\ $i_1\ne i_2=i_4\ne i_3=i_5\ne i_1$:
$$
E^p_5 = I_5 - \sum_{j_1,\ldots,j_5=0}^{p}
C_{j_5\ldots j_1}\Biggl(\sum\limits_{(j_2,j_4)}\Biggl(
\sum\limits_{(j_3,j_5)}
C_{j_5\ldots j_1}\Biggr)\Biggr),
$$

\vspace{7mm}

(VI).15.\ $i_1\ne i_2=i_5\ne i_3=i_4\ne i_1$:
$$
E^p_5 = I_5 - \sum_{j_1,\ldots,j_5=0}^{p}
C_{j_5\ldots j_1}\Biggl(\sum\limits_{(j_2,j_5)}\Biggl(
\sum\limits_{(j_3,j_4)}
C_{j_5\ldots j_1}\Biggr)\Biggr),
$$

\vspace{7mm}

(VII).1.\ $i_1=i_2=i_3\ne i_4=i_5$:
$$
E^p_5 = I_5 - \sum_{j_1,\ldots,j_5=0}^{p}
C_{j_5\ldots j_1}\Biggl(\sum\limits_{(j_4,j_5)}\Biggl(
\sum\limits_{(j_1,j_2,j_3)}
C_{j_5\ldots j_1}\Biggr)\Biggr),
$$

\vspace{7mm}

(VII).2.\ $i_1=i_2=i_4\ne i_3=i_5$:
$$
E^p_5 = I_5 - \sum_{j_1,\ldots,j_5=0}^{p}
C_{j_5\ldots j_1}\Biggl(\sum\limits_{(j_3,j_5)}\Biggl(
\sum\limits_{(j_1,j_2,j_4)}
C_{j_5\ldots j_1}\Biggr)\Biggr),
$$

\vspace{7mm}

(VII).3.\ $i_1=i_2=i_5\ne i_3=i_4$:
$$
E_p = I - \sum_{j_1,\ldots,j_5=0}^{p}
C_{j_5\ldots j_1}\Biggl(\sum\limits_{(j_3,j_4)}\Biggl(
\sum\limits_{(j_1,j_2,j_5)}
C_{j_5\ldots j_1}\Biggr)\Biggr),
$$

\vspace{8mm}

(VII).4.\ $i_2=i_3=i_4\ne i_1=i_5$:
$$
E^p_5 = I_5 - \sum_{j_1,\ldots,j_5=0}^{p}
C_{j_5\ldots j_1}\Biggl(\sum\limits_{(j_1,j_5)}\Biggl(
\sum\limits_{(j_2,j_3,j_4)}
C_{j_5\ldots j_1}\Biggr)\Biggr),
$$

\vspace{8mm}

(VII).5.\ $i_2=i_3=i_5\ne i_1=i_4$:
$$
E^p_5 = I_5 - \sum_{j_1,\ldots,j_5=0}^{p}
C_{j_5\ldots j_1}\Biggl(\sum\limits_{(j_1,j_4)}\Biggl(
\sum\limits_{(j_2,j_3,j_5)}
C_{j_5\ldots j_1}\Biggr)\Biggr),
$$

\vspace{8mm}

(VII).6.\ $i_2=i_4=i_5\ne i_1=i_3$:
$$
E^p_5 = I_5 - \sum_{j_1,\ldots,j_5=0}^{p}
C_{j_5\ldots j_1}\Biggl(\sum\limits_{(j_1,j_3)}\Biggl(
\sum\limits_{(j_2,j_4,j_5)}
C_{j_5\ldots j_1}\Biggr)\Biggr),
$$

\vspace{8mm}

(VII).7.\ $i_3=i_4=i_5\ne i_1=i_2$:
$$
E^p_5 = I_5 - \sum_{j_1,\ldots,j_5=0}^{p}
C_{j_5\ldots j_1}\Biggl(\sum\limits_{(j_1,j_2)}\Biggl(
\sum\limits_{(j_3,j_4,j_5)}
C_{j_5\ldots j_1}\Biggr)\Biggr),
$$

\vspace{8mm}

(VI).8.\ $i_1=i_3=i_5\ne i_2=i_4$:
$$
E^p_5 = I_5 - \sum_{j_1,\ldots,j_5=0}^{p}
C_{j_5\ldots j_1}\Biggl(\sum\limits_{(j_2,j_4)}\Biggl(
\sum\limits_{(j_1,j_3,j_5)}
C_{j_5\ldots j_1}\Biggr)\Biggr),
$$

\vspace{8mm}

(VII).9.\ $i_1=i_3=i_4\ne i_2=i_5$:
$$
E^p_5 = I_5 - \sum_{j_1,\ldots,j_5=0}^{p}
C_{j_5\ldots j_1}\Biggl(\sum\limits_{(j_2,j_5)}\Biggl(
\sum\limits_{(j_1,j_3,j_4)}
C_{j_5\ldots j_1}\Biggr)\Biggr),
$$

\vspace{8mm}

(VII).10.\ $i_1=i_4=i_5\ne i_2=i_3$:
$$
E^p_5 = I_5 - \sum_{j_1,\ldots,j_5=0}^{p}
C_{j_5\ldots j_1}\Biggl(\sum\limits_{(j_2,j_3)}\Biggl(
\sum\limits_{(j_1,j_4,j_5)}
C_{j_5\ldots j_1}\Biggr)\Biggr).
$$

\vspace{5mm}

\section{Method of Generalized Multiple Fourier Series. The Case 
of an Arbitrary Complete Orthonormal Systems  
of Functions in the Space $L_2([t,T])$ 
and Weight Functions $\psi_1(\tau),\ldots,\psi_k(\tau)\in L_2([t, T])$
}

\vspace{5mm}

For further consideration, let us 
consider the generalization of formulas (\ref{a1})--(\ref{a6}) 
for the case of an arbitrary multiplicity $k$ $(k\in\mathbb{N})$ of 
the iterated Ito stochastic integral $J[\psi^{(k)}]_{T,t}$ defined by (\ref{sodom20}).
In order to do this, let us
introduce some notations. 
Let us consider the unordered
set $\{1, 2, \ldots, k\}$ 
and separate it into two parts:
the first part consists of $r$ unordered 
pairs (sequence order of these pairs is also unimportant) and the 
second one consists of the 
remaining $k-2r$ numbers.
So, we have

\vspace{1mm}
\begin{equation}
\label{leto5007}
(\{
\underbrace{\{g_1, g_2\}, \ldots, 
\{g_{2r-1}, g_{2r}\}}_{\small{\hbox{part 1}}}
\},
\{\underbrace{q_1, \ldots, q_{k-2r}}_{\small{\hbox{part 2}}}
\}),
\end{equation}

\vspace{4mm}
\noindent
where 

\vspace{-2mm}
$$
\{g_1, g_2, \ldots, 
g_{2r-1}, g_{2r}, q_1, \ldots, q_{k-2r}\}=\{1, 2, \ldots, k\},
$$

\vspace{4mm}
\noindent
braces   
mean an unordered 
set, and pa\-ren\-the\-ses mean an ordered set.

We will say that (\ref{leto5007}) is a partition 
and consider the sum with respect to all possible
partitions

\vspace{1mm}
\begin{equation}
\label{leto5008}
\sum_{\stackrel{(\{\{g_1, g_2\}, \ldots, 
\{g_{2r-1}, g_{2r}\}\}, \{q_1, \ldots, q_{k-2r}\})}
{{}_{\{g_1, g_2, \ldots, 
g_{2r-1}, g_{2r}, q_1, \ldots, q_{k-2r}\}=\{1, 2, \ldots, k\}}}}
a_{g_1 g_2, \ldots, 
g_{2r-1} g_{2r}, q_1 \ldots q_{k-2r}},
\end{equation}

\vspace{5mm}
\noindent
where $a_{g_1 g_2, \ldots, 
g_{2r-1} g_{2r}, q_1 \ldots q_{k-2r}}\in \mathbb{R}.$

Below there are several examples of sums in the form (\ref{leto5008})

\vspace{2mm}
$$
\sum_{\stackrel{(\{g_1, g_2\})}{{}_{\{g_1, g_2\}=\{1, 2\}}}}
a_{g_1 g_2}=a_{12},
$$

\vspace{3mm}
$$
\sum_{\stackrel{(\{\{g_1, g_2\}, \{g_3, g_4\}\})}
{{}_{\{g_1, g_2, g_3, g_4\}=\{1, 2, 3, 4\}}}}
a_{g_1 g_2, g_3 g_4}=a_{12,34} + a_{13,24} + a_{23,14},
$$

\vspace{3mm}
$$
\sum_{\stackrel{(\{g_1, g_2\}, \{q_1, q_{2}\})}
{{}_{\{g_1, g_2, q_1, q_{2}\}=\{1, 2, 3, 4\}}}}
a_{g_1 g_2, q_1 q_{2}}=
$$

$$
=a_{12,34}+a_{13,24}+a_{14,23}
+a_{23,14}+a_{24,13}+a_{34,12},
$$

\vspace{3mm}
$$
\sum_{\stackrel{(\{g_1, g_2\}, \{q_1, q_{2}, q_3\})}
{{}_{\{g_1, g_2, q_1, q_{2}, q_3\}=\{1, 2, 3, 4, 5\}}}}
a_{g_1 g_2, q_1 q_{2}q_3}
=
$$

$$
=a_{12,345}+a_{13,245}+a_{14,235}
+a_{15,234}+a_{23,145}+a_{24,135}+
$$
$$
+a_{25,134}+a_{34,125}+a_{35,124}+a_{45,123},
$$

\vspace{3mm}
$$
\sum_{\stackrel{(\{\{g_1, g_2\}, \{g_3, g_{4}\}\}, \{q_1\})}
{{}_{\{g_1, g_2, g_3, g_{4}, q_1\}=\{1, 2, 3, 4, 5\}}}}
a_{g_1 g_2, g_3 g_{4},q_1}
=
$$

$$
=
a_{12,34,5}+a_{13,24,5}+a_{14,23,5}+
a_{12,35,4}+a_{13,25,4}+a_{15,23,4}+
$$
$$
+a_{12,54,3}+a_{15,24,3}+a_{14,25,3}+a_{15,34,2}+a_{13,54,2}+a_{14,53,2}+
$$
$$
+
a_{52,34,1}+a_{53,24,1}+a_{54,23,1}.
$$

\vspace{5mm}

Now we can write (\ref{tyyyxxx}) as

\vspace{1mm}

$$
J[\psi^{(k)}]_{T,t}=
\hbox{\vtop{\offinterlineskip\halign{
\hfil#\hfil\cr
{\rm l.i.m.}\cr
$\stackrel{}{{}_{p_1,\ldots,p_k\to \infty}}$\cr
}} }
\sum\limits_{j_1=0}^{p_1}\ldots
\sum\limits_{j_k=0}^{p_k}
C_{j_k\ldots j_1}\Biggl(
\prod_{l=1}^k\zeta_{j_l}^{(i_l)}+\sum\limits_{r=1}^{[k/2]}
(-1)^r \times
\Biggr.
$$

\vspace{3mm}
\begin{equation}
\label{leto6000}
\times
\sum_{\stackrel{(\{\{g_1, g_2\}, \ldots, 
\{g_{2r-1}, g_{2r}\}\}, \{q_1, \ldots, q_{k-2r}\})}
{{}_{\{g_1, g_2, \ldots, 
g_{2r-1}, g_{2r}, q_1, \ldots, q_{k-2r}\}=\{1, 2, \ldots, k\}}}}
\prod\limits_{s=1}^r
{\bf 1}_{\{i_{g_{{}_{2s-1}}}=~i_{g_{{}_{2s}}}\ne 0\}}
\Biggl.{\bf 1}_{\{j_{g_{{}_{2s-1}}}=~j_{g_{{}_{2s}}}\}}
\prod_{l=1}^{k-2r}\zeta_{j_{q_l}}^{(i_{q_l})}\Biggr),
\end{equation}

\vspace{5mm}
\noindent
where $[x]$ is an integer part of a real number $x,$
$\prod\limits_{\emptyset}
\stackrel{\sf def}{=}1,$ $\sum\limits_{\emptyset}
\stackrel{\sf def}{=}0;$ 
another notations are the same as in Theorem~1.

another notations are the same as in Theorem {\bf 1}.

\vspace{2mm}

In particular, from (\ref{leto6000}) for $k=5$ we obtain

\vspace{3mm}

$$
J[\psi^{(5)}]_{T,t}=
\hbox{\vtop{\offinterlineskip\halign{
\hfil#\hfil\cr
{\rm l.i.m.}\cr
$\stackrel{}{{}_{p_1,\ldots,p_5\to \infty}}$\cr
}} }\sum_{j_1=0}^{p_1}\ldots\sum_{j_5=0}^{p_5}
C_{j_5\ldots j_1}\Biggl(
\prod_{l=1}^5\zeta_{j_l}^{(i_l)}-\Biggr.
$$

\vspace{2mm}
$$
-
\sum\limits_{\stackrel{(\{g_1, g_2\}, \{q_1, q_{2}, q_3\})}
{{}_{\{g_1, g_2, q_{1}, q_{2}, q_3\}=\{1, 2, 3, 4, 5\}}}}
{\bf 1}_{\{i_{g_{{}_{1}}}=~i_{g_{{}_{2}}}\ne 0\}}
{\bf 1}_{\{j_{g_{{}_{1}}}=~j_{g_{{}_{2}}}\}}
\prod_{l=1}^{3}\zeta_{j_{q_l}}^{(i_{q_l})}+
$$

\vspace{2mm}
$$
+
\sum_{\stackrel{(\{\{g_1, g_2\}, 
\{g_{3}, g_{4}\}\}, \{q_1\})}
{{}_{\{g_1, g_2, g_{3}, g_{4}, q_1\}=\{1, 2, 3, 4, 5\}}}}
{\bf 1}_{\{i_{g_{{}_{1}}}=~i_{g_{{}_{2}}}\ne 0\}}
{\bf 1}_{\{j_{g_{{}_{1}}}=~j_{g_{{}_{2}}}\}}
\Biggl.{\bf 1}_{\{i_{g_{{}_{3}}}=~i_{g_{{}_{4}}}\ne 0\}}
{\bf 1}_{\{j_{g_{{}_{3}}}=~j_{g_{{}_{4}}}\}}
\zeta_{j_{q_1}}^{(i_{q_1})}\Biggr).
$$

\vspace{6mm}
\noindent
The last equality obviously agrees with
(\ref{a5}).

Further, we will use the definition of the multiple Wiener 
stochastic integral from \cite{ito1951}, \cite{Kuo} to generalize Theorem
1 to the case of an arbitrary 
complete orthonormal system of functions in the space $L_2([t, T])$
and $\psi_1(\tau),$ $\ldots,\psi_k(\tau)\in L_2([t, T]).$

Consider the following step function on the hypercube $[t, T]^k$

\vspace{-2mm}
\begin{equation}
\label{chain3}
\Phi_N(t_1,\ldots,t_k)=\sum\limits_{l_1,\ldots,l_k=0}^{N-1}
a_{l_1 \ldots l_k} {\bf 1}_{[\tau_{l_1},\tau_{l_1+1})}(t_1) \ldots
{\bf 1}_{[\tau_{l_k},\tau_{l_k+1})}(t_k),
\end{equation}

\vspace{3mm}
\noindent
where $a_{l_1 \ldots l_k}\in\mathbb{R}$ and such that 
$a_{l_1 \ldots l_k}=0$ if $l_p=l_q$ for some $p\ne q,$

\vspace{1mm}
$$
{\bf 1}_A (\tau)=\left\{
\begin{matrix}
1\ &{\rm if}\ \tau\in A \cr\cr
0\ &\hbox{\rm otherwise}
\end{matrix}\right.,
$$

\vspace{4mm}
\noindent
$N\in\mathbb{N},$ $\left\{\tau_{j}\right\}_{j=0}^{N}$ is a partition of
$[t,T],$ which satisfies the condition (\ref{1111}):

\vspace{1mm}
\begin{equation}
\label{1111xxx1}
t=\tau_0<\ldots <\tau_N=T,\ \ \
\Delta_N=
\hbox{\vtop{\offinterlineskip\halign{
\hfil#\hfil\cr
{\rm max}\cr
$\stackrel{}{{}_{0\le j\le N-1}}$\cr
}} }\Delta\tau_j\to 0\ \ \hbox{if}\ \ N\to \infty,\ \ \ 
\Delta\tau_j=\tau_{j+1}-\tau_j.
\end{equation}

\vspace{4mm}

Let us define the multiple Wiener stochastic integral for $\Phi_N(t_1,\ldots,t_k)$ 
\cite{ito1951}, \cite{Kuo}

\vspace{1mm}
\begin{equation}
\label{chain9}
J'[\Phi_N]_{T,t}^{(i_1\ldots i_k)}\stackrel{\sf def}{=}
\sum\limits_{l_1,\ldots,l_k=0}^{N-1}
a_{l_1 \ldots l_k}
\Delta{\bf w}_{\tau_{l_1}}^{(i_1)}\ldots \Delta{\bf w}_{\tau_{l_k}}^{(i_k)},
\end{equation}

\vspace{4mm}
\noindent
where $\Delta{\bf w}_{\tau_{j}}^{(i)}=
{\bf w}_{\tau_{j+1}}^{(i)}-{\bf w}_{\tau_{j}}^{(i)},$\
$i=0, 1,\ldots,m,$\ ${\bf w}_{\tau}^{(0)}=\tau.$

It is known (see \cite{Kuo}, Lemma~9.6.4)
that for any $\Phi(t_1,\ldots,t_k)\in L_2([t, T]^k)$ 
there exists a sequence of step functions $\Phi_N(t_1,\ldots,t_k)$ of the form (\ref{chain3})
such that

\vspace{1mm}
\begin{equation}
\label{chain15}
\lim\limits_{N\to\infty} \int\limits_{[t,T]^k}
\left(\Phi(t_1,\ldots,t_k)-\Phi_N(t_1,\ldots,t_k)\right)^2 dt_1\ldots dt_k=0.
\end{equation}

\vspace{4mm}

We have

\vspace{-1mm}
$$
\Phi_N(t_1,\ldots,t_k)=\sum\limits_{l_1,\ldots,l_k=0}^{N-1}
a_{l_1 \ldots l_k} {\bf 1}_{[\tau_{l_1},\tau_{l_1+1})}(t_1) \ldots
{\bf 1}_{[\tau_{l_k},\tau_{l_k+1})}(t_k)=
$$

\vspace{1mm}
\begin{equation}
\label{chain5}
=\sum\limits_{(l_1,\ldots,l_k)}
\sum_{\stackrel{l_1,\ldots,l_k=0}{{}_{l_1<l_2<\ldots < l_k}}}^{N-1}
a_{l_1 \ldots l_k} {\bf 1}_{[\tau_{l_1},\tau_{l_1+1})}(t_1) \ldots
{\bf 1}_{[\tau_{l_k},\tau_{l_k+1})}(t_k),
\end{equation}

\vspace{4mm}
\noindent
where permutations $(l_1,\ldots,l_k)$ when summing are 
performed only in the expression $l_1<l_2<\ldots < l_k$
(recall that $a_{l_1 \ldots l_k}=0$ if $l_p=l_q$ for some $p\ne q$).

Using (\ref{chain5}), we get

\vspace{-1mm}
\begin{equation}
\label{chain30}
\sum_{(t_1,\ldots,t_k)}
\int\limits_{t}^{T}
\ldots
\int\limits_{t}^{t_2}
\Phi_N(t_1,\ldots,t_k)d{\bf w}_{t_1}^{(i_1)}
\ldots
d{\bf w}_{t_k}^{(i_k)}=
\end{equation}

\vspace{1mm}
$$
=\sum\limits_{(l_1,\ldots,l_k)}
\sum_{\stackrel{l_1,\ldots,l_k=0}{{}_{l_1<l_2<\ldots < l_k}}}^{N-1}
a_{l_1 \ldots l_k} 
\Delta{\bf w}_{\tau_{l_1}}^{(i_1)} \ldots \Delta{\bf w}_{\tau_{l_k}}^{(i_k)}=
$$

\vspace{1mm}
$$
=\sum\limits_{\stackrel{l_1,\ldots,l_k=0}{{}_{l_q\ne l_r;\ q\ne r;\ 
q, r=1,\ldots, k}}}^{N-1}
a_{l_1 \ldots l_k} 
\Delta{\bf w}_{\tau_{l_1}}^{(i_1)} \ldots \Delta{\bf w}_{\tau_{l_k}}^{(i_k)}=
$$

\vspace{1mm}
\begin{equation}
\label{chain10}
=J'[\Phi_N]_{T,t}^{(i_1\ldots i_k)}\ \ \ \hbox{w.\ p.\ 1},
\end{equation}

\vspace{3mm}
\noindent
where permutations $(t_1,\ldots,t_k)$ when summing are 
performed only in the values
$d{\bf w}_{t_1}^{(i_1)}
\ldots $
$d{\bf w}_{t_k}^{(i_k)}$ 
and permutations $(l_1,\ldots,l_k)$ when summing are 
performed only in the expression $l_1<l_2<\ldots < l_k.$
At the same time the indices near 
upper 
limits of integration in the iterated stochastic integrals in (\ref{chain30}) are changed 
correspondently and if $t_r$ swapped with $t_q$ in the  
permutation $(t_1,\ldots,t_k)$, then $i_r$ swapped with $i_q$ in 
the permutation $(i_1,\ldots,i_k)$ (see (\ref{chain30})).
In addition, the multiple Wiener stochastic integral 
$J'[\Phi_N]_{T,t}^{(i_1\ldots i_k)}$ is defined by (\ref{chain9})
and

$$
\int\limits_{t}^{T}
\ldots
\int\limits_{t}^{t_2}
\Phi_N(t_1,\ldots,t_k)d{\bf w}_{t_1}^{(i_1)}
\ldots
d{\bf w}_{t_k}^{(i_k)}
$$

\vspace{4mm}
\noindent
is the iterated Ito stochastic integral.

Since the integration of a bounded function with respect to the
set of measure zero for Lebesgue integrals gives zero result, then the 
following formula is correct for these integrals

\begin{equation}
\label{riemann}
\int\limits_{[t, T]^k}|G(t_1,\ldots,t_k)|dt_1\ldots dt_k=
\sum_{(t_1,\ldots,t_k)}
\int\limits_{t}^{T}
\ldots
\int\limits_{t}^{t_2}
|G(t_1,\ldots,t_k)|dt_1\ldots dt_k,
\end{equation}

\vspace{4mm}
\noindent
where permutations $(t_1,\ldots,t_k)$ when summing
are performed only 
in the 
va\-lues $dt_1,\ldots, dt_k$. At the same time the indexes near upper 
limits of integration are changed correspondently
and the function $|G(t_1,\ldots,t_k)|$ is assumed to be integrable in 
the hypercube $[t, T]^k.$

Using (\ref{chain15}), (\ref{chain10}), (\ref{riemann}), and
standard moment properties of the Ito stochastic integral, we have

\vspace{1mm}
$$
{\sf M}\left\{\left(J'[\Phi_N]_{T,t}^{(i_1\ldots i_k)}-
J'[\Phi_M]_{T,t}^{(i_1\ldots i_k)}\right)^2\right\}\le
$$

\vspace{2mm}
$$
\le C_k 
\sum_{(t_1,\ldots,t_k)}
\int\limits_{t}^{T}
\ldots
\int\limits_{t}^{t_2}
\left(\Phi_N(t_1,\ldots,t_k)-\Phi_M(t_1,\ldots,t_k)\right)^2 dt_1
\ldots dt_k=
$$

\vspace{2mm}
$$
=C_k 
\int\limits_{[t,T]^k}
\left(\Phi_N(t_1,\ldots,t_k)-\Phi_M(t_1,\ldots,t_k)\right)^2 dt_1
\ldots dt_k=
$$

\vspace{2mm}
$$
=C_k\left\Vert \Phi_N-\Phi_M\right\Vert_{L_2([t, T]^k)}^2\le
$$

\vspace{2mm}
$$
\le 2 C_k \left(\left\Vert \Phi_N-\Phi\right\Vert_{L_2([t, T]^k)}^2+
\left\Vert \Phi-\Phi_M\right\Vert_{L_2([t, T]^k)}^2\right)^2\ \to 0
$$

\vspace{5mm}
\noindent
if $N,M\to\infty,$ 
where constant $C_k$ 
depends only
on the multiplicity $k$ of the multiple Wiener stochastic integral.

Thus, there exists the limit 

$$
\hbox{\vtop{\offinterlineskip\halign{
\hfil#\hfil\cr
{\rm l.i.m.}\cr
$\stackrel{}{{}_{N\to \infty}}$\cr
}} }J'[\Phi_N]_{T,t}^{(i_1\ldots i_k)}.
$$

\vspace{4mm}

We will define the multiple Wiener stochastic integral for $\Phi(t_1,\ldots,t_k)\in L_2([t, T]^k)$ 
by the formula 

\begin{equation}
\label{WiI}
J'[\Phi]_{T,t}^{(i_1\ldots i_k)}\stackrel{\sf def}{=}
\hbox{\vtop{\offinterlineskip\halign{
\hfil#\hfil\cr
{\rm l.i.m.}\cr
$\stackrel{}{{}_{N\to \infty}}$\cr
}} }J'[\Phi_N]_{T,t}^{(i_1\ldots i_k)}=
\hbox{\vtop{\offinterlineskip\halign{
\hfil#\hfil\cr
{\rm l.i.m.}\cr
$\stackrel{}{{}_{N\to \infty}}$\cr
}} }
\sum\limits_{l_1,\ldots,l_k=0}^{N-1}
a_{l_1 \ldots l_k}
\Delta{\bf w}_{\tau_{l_1}}^{(i_1)}\ldots \Delta{\bf w}_{\tau_{l_k}}^{(i_k)},
\end{equation}

\vspace{4mm}
\noindent
where $\Phi_N(t_1,\ldots,t_k)$ is defined by 
(\ref{chain3}),
$\Delta{\bf w}_{\tau_{j}}^{(i)}=
{\bf w}_{\tau_{j+1}}^{(i)}-{\bf w}_{\tau_{j}}^{(i)},$\
$i=0, 1,\ldots,m,$\ ${\bf w}_{\tau}^{(0)}=\tau.$

Let us prove the following equality 

\begin{equation}
\label{Wi110}
J'[\Phi]_{T,t}^{(i_1\ldots i_k)}=\sum_{(t_1,\ldots,t_k)}
\int\limits_{t}^{T}
\ldots
\int\limits_{t}^{t_2}
\Phi(t_1,\ldots,t_k)d{\bf w}_{t_1}^{(i_1)}
\ldots
d{\bf w}_{t_k}^{(i_k)}\ \ \ \hbox{w.\ p.\ 1},
\end{equation}

\vspace{3mm}
\noindent
where permutations $(t_1,\ldots,t_k)$ when summing are 
performed only in the values
$d{\bf w}_{t_1}^{(i_1)}
\ldots $
$d{\bf w}_{t_k}^{(i_k)}.$ At the same time the indices near 
upper 
limits of integration in the iterated stochastic integrals are changed 
correspondently and if $t_r$ swapped with $t_q$ in the  
permutation $(t_1,\ldots,t_k)$, then $i_r$ swapped with $i_q$ in 
the permutation $(i_1,\ldots,i_k).$ 
In addition, the multiple Wiener stochastic integral 
$J'[\Phi]_{T,t}^{(i_1\ldots i_k)}$ is defined by (\ref{WiI})
and 

\vspace{-1mm}
$$
\int\limits_{t}^{T}
\ldots
\int\limits_{t}^{t_2}
\Phi(t_1,\ldots,t_k)d{\bf w}_{t_1}^{(i_1)}
\ldots
d{\bf w}_{t_k}^{(i_k)}
$$

\vspace{4mm}
\noindent
is the iterated Ito stochastic integral.

The equality (\ref{Wi110}) has already been proved for the case 
$\Phi(t_1,\ldots,t_k)=\Phi_N(t_1,\ldots,t_k)$ (see (\ref{chain10})).
From (\ref{chain10}) we have
$$
J'[\Phi_N]_{T,t}^{(i_1\ldots i_k)}=
\sum_{(t_1,\ldots,t_k)}
\int\limits_{t}^{T}
\ldots
\int\limits_{t}^{t_2}
\Phi_N(t_1,\ldots,t_k)d{\bf w}_{t_1}^{(i_1)}
\ldots
d{\bf w}_{t_k}^{(i_k)}=
$$

\vspace{2mm}
$$
=\sum_{(t_1,\ldots,t_k)}
\int\limits_{t}^{T}
\ldots
\int\limits_{t}^{t_2}
\Phi(t_1,\ldots,t_k)d{\bf w}_{t_1}^{(i_1)}
\ldots
d{\bf w}_{t_k}^{(i_k)}+
$$

\vspace{2mm}
\begin{equation}
\label{chain11}
+\sum_{(t_1,\ldots,t_k)}
\int\limits_{t}^{T}
\ldots
\int\limits_{t}^{t_2}
\left(\Phi_N(t_1,\ldots,t_k)-\Phi(t_1,\ldots,t_k)\right)d{\bf w}_{t_1}^{(i_1)}
\ldots
d{\bf w}_{t_k}^{(i_k)}\ \ \ \hbox{w.~p.~1.}
\end{equation}

\vspace{5mm}

Passing to the limit $\hbox{\vtop{\offinterlineskip\halign{
\hfil#\hfil\cr
{\rm l.i.m.}\cr
$\stackrel{}{{}_{N\to \infty}}$\cr
}} }$ in the equality (\ref{chain11}), we obtain

\vspace{1mm}
$$
J'[\Phi]_{T,t}^{(i_1\ldots i_k)}=
\sum_{(t_1,\ldots,t_k)}
\int\limits_{t}^{T}
\ldots
\int\limits_{t}^{t_2}
\Phi(t_1,\ldots,t_k)d{\bf w}_{t_1}^{(i_1)}
\ldots
d{\bf w}_{t_k}^{(i_k)}+
$$

\vspace{2mm}
\begin{equation}
\label{chain12}
+\hbox{\vtop{\offinterlineskip\halign{
\hfil#\hfil\cr
{\rm l.i.m.}\cr
$\stackrel{}{{}_{N\to \infty}}$\cr
}} }\sum_{(t_1,\ldots,t_k)}
\int\limits_{t}^{T}
\ldots
\int\limits_{t}^{t_2}
\left(\Phi_N(t_1,\ldots,t_k)-\Phi(t_1,\ldots,t_k)\right)d{\bf w}_{t_1}^{(i_1)}
\ldots
d{\bf w}_{t_k}^{(i_k)}\ \ \ \hbox{w.~p.~1.}
\end{equation}

\vspace{5mm}

Using (\ref{chain15}), (\ref{riemann}),  and
standard moment properties of the Ito stochastic integral,
we get

\vspace{1mm}
$$
{\sf M}\left\{\left(
\sum_{(t_1,\ldots,t_k)}
\int\limits_{t}^{T}
\ldots
\int\limits_{t}^{t_2}
\left(\Phi_N(t_1,\ldots,t_k)-\Phi(t_1,\ldots,t_k)\right)d{\bf w}_{t_1}^{(i_1)}
\ldots
d{\bf w}_{t_k}^{(i_k)}\right)^2\right\}\le
$$

\vspace{2mm}
$$
\le C_k 
\sum_{(t_1,\ldots,t_k)}
\int\limits_{t}^{T}
\ldots
\int\limits_{t}^{t_2}
\left(\Phi_N(t_1,\ldots,t_k)-\Phi(t_1,\ldots,t_k)\right)^2 dt_1
\ldots dt_k=
$$

\vspace{2mm}
\begin{equation}
\label{chain20}
=C_k 
\int\limits_{[t,T]^k}
\left(\Phi_N(t_1,\ldots,t_k)-\Phi(t_1,\ldots,t_k)\right)^2 dt_1
\ldots dt_k\ \to 0
\end{equation}

\vspace{4mm}
\noindent
if $N\to\infty,$ 
where constant $C_k$ 
depends only
on the multiplicity $k$ of the multiple Wiener stochastic integral.

The relations (\ref{chain12}) and (\ref{chain20}) prove the equality 
(\ref{Wi110}).
From (\ref{Wi110}) we have

\begin{equation}
\label{wi1001}
J[\psi^{(k)}]_{T,t}^{(i_1\ldots i_k)}=\int\limits_t^T\psi_k(t_k) \ldots \int\limits_t^{t_{2}}
\psi_1(t_1) d{\bf w}_{t_1}^{(i_1)}\ldots
d{\bf w}_{t_k}^{(i_k)}=J'[K]_{T,t}^{(i_1\ldots i_k)}\ \ \ \hbox{w.\ p.\ 1},
\end{equation}

\vspace{3mm}
\noindent
where $J[\psi^{(k)}]_{T,t}^{(i_1\ldots i_k)}$ is the iterated Ito stochastic integral
(\ref{sodom20}), $K=K(t_1,\ldots,t_k)$ is defined by (\ref{ppp}), i.e.

\begin{equation}
\label{chain200}
K(t_1,\ldots,t_k)=
\left\{\begin{matrix}
\psi_1(t_1)\ldots \psi_k(t_k),\ &t_1<\ldots<t_k\cr\cr\cr
0,\ &\hbox{\rm otherwise}
\end{matrix}
\right.,
\end{equation}

\vspace{4mm}
\noindent
where $\psi_1(\tau),\ldots,\psi_k(\tau)\in L_2([t,T])$,\ $t_1,\ldots,t_k\in [t, T]$ $(k\ge 2)$ and 
$K(t_1)\equiv\psi_1(t_1)$ for $t_1\in[t, T].$

Applying (\ref{wi1001}) and the linearity property of the Ito stochastic integral, we obtain

\vspace{1mm}
$$
J[\psi^{(k)}]_{T,t}^{(i_1\ldots i_k)}=J'[K]_{T,t}^{(i_1\ldots i_k)}=
$$

\vspace{2mm}
\begin{equation}
\label{chain102}
=\sum_{j_1=0}^{p_1}\ldots
\sum_{j_k=0}^{p_k}
C_{j_k\ldots j_1}
J'[\phi_{j_1}\ldots \phi_{j_k}]_{T,t}^{(i_1\ldots i_k)}+
J'[R_{p_1\ldots p_k}]_{T,t}^{(i_1\ldots i_k)}\ \ \ \hbox{w.~p.~1,}
\end{equation}

\vspace{4mm}
\noindent
where
$$
R_{p_1\ldots p_k}(t_1,\ldots,t_k)\stackrel{{\rm def}}{=}
K(t_1,\ldots,t_k)-
\sum_{j_1=0}^{p_1}\ldots
\sum_{j_k=0}^{p_k}
C_{j_k\ldots j_1}
\prod_{l=1}^k\phi_{j_l}(t_l)
$$

\vspace{3mm}
\noindent
and
\begin{equation}
\label{chain300}
C_{j_k\ldots j_1}=\int\limits_{[t,T]^k}
K(t_1,\ldots,t_k)\prod_{l=1}^{k}\phi_{j_l}(t_l)dt_1\ldots dt_k
\end{equation}

\vspace{4mm}
\noindent
is the Fourier coefficient corresponding to $K(t_1,\ldots,t_k).$

Again applying (\ref{Wi110}), we have

\vspace{1mm}
$$
J'[R_{p_1\ldots p_k}]_{T,t}^{(i_1\ldots i_k)}
=
$$

\vspace{2mm}
\begin{equation}
\label{wi2005}
=
\sum_{(t_1,\ldots,t_k)}
\int\limits_{t}^{T}
\ldots
\int\limits_{t}^{t_2}
\Biggl(K(t_1,\ldots,t_k)-
\sum_{j_1=0}^{p_1}\ldots
\sum_{j_k=0}^{p_k}
C_{j_k\ldots j_1}
\prod_{l=1}^k\phi_{j_l}(t_l)\Biggr)
d{\bf w}_{t_1}^{(i_1)}
\ldots
d{\bf w}_{t_k}^{(i_k)},
\end{equation}

\vspace{3mm}
\noindent
where permutations $(t_1,\ldots,t_k)$ when summing are performed only 
in the values $d{\bf w}_{t_1}^{(i_1)}
\ldots $
$d{\bf w}_{t_k}^{(i_k)}$. At the same time the indices near 
upper limits of integration in the iterated stochastic integrals 
are changed correspondently and if $t_r$ swapped with $t_q$ in the  
permutation $(t_1,\ldots,t_k)$, then $i_r$ swapped with $i_q$ in the 
permutation $(i_1,\ldots,i_k).$
In addition, the multiple Wiener stochastic integral
$J'[R_{p_1\ldots p_k}]_{T,t}^{(i_1\ldots i_k)}$ is defined by 
(\ref{WiI}).

According to (\ref{sos1z}), (\ref{riemann}), and
the standard moment properties of the Ito stochastic integral, we have

\vspace{1mm}
$$
{\sf M}\left\{\left(J'[R_{p_1\ldots p_k}]_{T,t}^{(i_1\ldots i_k)}\right)^2\right\}
\le 
$$

\vspace{2mm}
$$
\le C_k
\sum_{(t_1,\ldots,t_k)}
\int\limits_{t}^{T}
\ldots
\int\limits_{t}^{t_2}
\left(K(t_1,\ldots,t_k)-
\sum_{j_1=0}^{p_1}\ldots
\sum_{j_k=0}^{p_k}
C_{j_k\ldots j_1}
\prod_{l=1}^k\phi_{j_l}(t_l)\right)^2
dt_1
\ldots
dt_k=
$$

\vspace{2mm}
\begin{equation}
\label{chain8810}
=C_k\int\limits_{[t,T]^k}
\left(K(t_1,\ldots,t_k)-
\sum_{j_1=0}^{p_1}\ldots
\sum_{j_k=0}^{p_k}
C_{j_k\ldots j_1}
\prod_{l=1}^k\phi_{j_l}(t_l)\right)^2
dt_1
\ldots
dt_k\to 0
\end{equation}

\vspace{4mm}
\noindent
if $p_1,\ldots,p_k\to\infty,$ where constant $C_k$ 
depends only
on the multiplicity $k$ of the 
iterated  Ito stochastic integral
$J[\psi^{(k)}]_{T,t}^{(i_1\ldots i_k)}$.

Applying (\ref{chain102}) and (\ref{chain8810}), we obtain the following expansion

\vspace{1mm}
\begin{equation}
\label{chain102xx}
J[\psi^{(k)}]_{T,t}^{(i_1\ldots i_k)}=
\hbox{\vtop{\offinterlineskip\halign{
\hfil#\hfil\cr
{\rm l.i.m.}\cr
$\stackrel{}{{}_{p_1,\ldots,p_k\to \infty}}$\cr
}} }\sum_{j_1=0}^{p_1}\ldots
\sum_{j_k=0}^{p_k}
C_{j_k\ldots j_1}
J'[\phi_{j_1}\ldots \phi_{j_k}]_{T,t}^{(i_1\ldots i_k)}.
\end{equation}

\vspace{5.5mm}

In \cite{11a} (Sect.~1.14, Theorem~1.23), \cite{new-2023a} (Theorem~5) it is shown that

$$
J'[\phi_{j_1}\ldots\phi_{j_k}]^{(i_1 \ldots i_k)}_{T,t}=
\prod_{l=1}^k\zeta_{j_l}^{(i_l)}+\sum\limits_{r=1}^{[k/2]}
(-1)^r \times \Biggr.
$$

\vspace{2mm}
\begin{equation}
\label{leto6000xxa}
\times
\sum_{\stackrel{(\{\{g_1, g_2\}, \ldots, 
\{g_{2r-1}, g_{2r}\}\}, \{q_1, \ldots, q_{k-2r}\})}
{{}_{\{g_1, g_2, \ldots, 
g_{2r-1}, g_{2r}, q_1, \ldots, q_{k-2r}\}=\{1,2, \ldots, k\}}}}
\prod\limits_{s=1}^r
{\bf 1}_{\{i_{g_{{}_{2s-1}}}=~i_{g_{{}_{2s}}}\ne 0\}}
\Biggl.{\bf 1}_{\{j_{g_{{}_{2s-1}}}=~j_{g_{{}_{2s}}}\}}
\prod_{l=1}^{k-2r}\zeta_{j_{q_l}}^{(i_{q_l})}
\end{equation}

\vspace{5mm}
\noindent 
w.~p.~{\rm 1,} where $\{\phi_j(x)\}_{j=0}^{\infty}$ is an arbitrary complete orthonormal system  
of functions in the space $L_2([t,T]),$
$i_1,\ldots,i_k=0,1,\ldots,m,$
$J'[\phi_{j_1}\ldots\phi_{j_k}]^{(i_1 \ldots i_k)}_{T,t}$
is defined by {\rm (\ref{WiI}),}
$[x]$ is an integer part of a real number $x,$
$\prod\limits_{\emptyset}
\stackrel{\sf def}{=}1,$ $\sum\limits_{\emptyset}
\stackrel{\sf def}{=}0;$ 
another notations are the same as in Theorem~1.

\vspace{2mm}

Note that the right-hand side of 
(\ref{leto6000xxa}) is nothing but the Hermite polynomial 
of random vector argument with components 
$\zeta_{j_1}^{(i_1)},\ldots,\zeta_{j_k}^{(i_k)}.$

The equalities (\ref{chain102xx}) and (\ref{leto6000xxa})
prove Theorem~1 for the case 
of an arbitrary complete orthonormal systems  
of functions
$\{\phi_j(x)\}_{j=0}^{\infty}$ is  in the space $L_2([t,T])$
and $\psi_1(\tau),\ldots,\psi_k(\tau)\in L_2([t,T])$.

Let us find the representation of the right-hand side
of (\ref{leto6000xxa}) through the product of Hermite
polynomials of scalar arguments.

\vspace{2mm}

{\it We will say that the condition {\rm ($\star\star$)} is fulfilled
for the multi-index $(i_1\ldots i_k)$ $(i_1,\ldots,i_k=0, 1,\ldots, m)$ if
$m_1,\ldots,m_k$ are multiplicities of the elements $i_1,\ldots,i_k,$ respectively$,$ i.e.

\vspace{-1mm}
$$
\{i_1,\ldots, i_k\}\hspace{-0.4mm}=\hspace{-0.4mm}\{\overbrace{{i_1, \ldots, i_1}}^{m_1},
\overbrace{{i_2, \ldots, i_2}}^{m_2},
\ldots, \overbrace{{i_r, \ldots, i_r}}^{m_r}\}\ \ \ (m_{r+1}=\ldots=m_k=0),
$$

\vspace{3mm}
\noindent
where $r=1,\ldots, k,$ braces   
mean an unordered 
set, and pa\-ren\-the\-ses mean an ordered set. At that, 
$m_1+\ldots+m_k=k,$\ $m_1,\ldots, m_k=0,1,\ldots,k,$\ 
and all elements with nonzero multiplicities are pairwise different.}

\vspace{2mm}

Let the condition {\rm ($\star\star$)} is fulfilled
for the mul\-ti-\-in\-dex $(i_1 \ldots i_k).$ Then

\vspace{2mm}
$$
J'\left[\phi_{j_1}\ldots \phi_{j_k}\right]_{T,t}^{(i_1\ldots i_k)}
=
J'\biggl[\underbrace{\phi_{j_{g_1}}
\ldots \phi_{j_{g_{{}_{m_1}}}}}_{m_1}
\underbrace{\phi_{j_{g_{m_1+1}}}
\ldots \phi_{j_{g_{m_1+m_2}}}}_{m_2}\ldots \biggr.
$$

\vspace{2mm}
\begin{equation}
\label{newe11}
\biggl.\ldots
\underbrace{\phi_{j_{g_{m_1+m_2+\ldots+m_{k-1}+1}}}\ldots
\phi_{j_{g_{m_1+m_2+\ldots+m_k}}}}_{m_k}\biggr]_{T,t}^
{(\overbrace{{}_{i_1 \ldots i_1}}^{m_1}
\overbrace{{}_{i_2 \ldots i_2}}^{m_2}
\ldots \overbrace{{}_{i_k \ldots i_k}}^{m_k})}
\end{equation}

\vspace{4mm}
\noindent
w.~p.~1, where 
$J'\left[\phi_{j_1}\ldots \phi_{j_k}\right]_{T,t}^{(i_1\ldots i_k)}$
is defined by (\ref{WiI}),
$\Phi(t_1,\ldots,t_k)=\phi_{j_1}(t_1)\ldots \phi_{j_k}(t_k),$
$\{\phi_j(x)\}_{j=0}^{\infty}$ is an arbitrary complete orthonormal system  
of functions in the space $L_2([t,T]),$ 
$\{j_{g_1},\ldots,j_{g_{m_1+m_2+\ldots+m_k}}\}$ 
$=\{j_1,\ldots,j_k\}$,
braces   
mean an unordered 
set, and pa\-ren\-the\-ses mean an ordered set.

From (\ref{newe11}) we have 

$$
J'\left[\phi_{j_{1}}\ldots \phi_{j_{k}}\right]_{T,t}^{(i_1\ldots i_k)}=
$$

\vspace{2mm}
$$
=
J'\left[\phi_{j_{g_1}}
\ldots \phi_{j_{g_{{}_{m_1}}}}\right]_{T,t}^{
(\hspace{0.5mm}\overbrace{{}_{i_1 \ldots i_1}}^{m_1}\hspace{0.5mm})}
\cdot 
J'\left[\phi_{j_{g_{m_1+1}}}
\ldots \phi_{j_{g_{m_1+m_2}}}\right]_{T,t}^{
(\hspace{0.5mm}\overbrace{{}_{i_2 \ldots i_2}}^{m_2}\hspace{0.5mm})}
\cdot \ldots 
$$

\vspace{2mm}
\begin{equation}
\label{ziko30}
\ldots \cdot 
J'\left[\phi_{j_{g_{m_1+m_2+\ldots+m_{k-1}+1}}}\ldots
\phi_{j_{g_{m_1+m_2+\ldots+m_k}}}\right]_{T,t}^{
(\hspace{0.5mm}\overbrace{{}_{i_k \ldots i_k}}^{m_k}\hspace{0.5mm})}
\end{equation}

\vspace{2mm}
\noindent
w.~p.~1, where
\begin{equation}
\label{ziko10}
J'\left[\phi_{j_{g_{m_1+m_2+\ldots+m_{l-1}+1}}}\ldots
\phi_{j_{g_{m_1+m_2+\ldots+m_l}}}\right]_{T,t}^{
(\hspace{0.5mm}\overbrace{{}_{i_l \ldots i_l}}^{m_l}\hspace{0.5mm})}
\stackrel{\sf def}{=}1\ \ \ \hbox{for}\ \ \ m_l=0.
\end{equation}

\vspace{5mm}

The detailed proof of the equality (\ref{ziko30}) 
is given in \cite{11a} (Sect.~1.14), \cite{new-2023a}, Sect.~2.2).

Let us consider the following multiple Wiener stochastic integral 

$$
J'\left[\phi_{j_{g_{m_1+m_2+\ldots+m_{l-1}+1}}}\ldots
\phi_{j_{g_{m_1+m_2+\ldots+m_l}}}\right]_{T,t}^{
(\hspace{0.5mm}\overbrace{{}_{i_l \ldots i_l}}^{m_l}\hspace{0.5mm})}\ \ \ (m_l>0),
$$

\vspace{4mm}
\noindent
where we suppose that 

$$
\bigl\{j_{g_{m_1+m_2+\ldots+m_{l-1}+1}}, \ldots, j_{g_{m_1+m_2+\ldots+m_{l}}}
\bigr\}=
$$
\begin{equation}
\label{ziko999}
=\bigl\{\underbrace{j_{h_{1,l}}, \ldots, j_{h_{1,l}}}_{n_{1,l}}\ \hspace{-1mm},
\underbrace{j_{h_{2,l}}, \ldots, j_{h_{2,l}}}_{n_{2,l}}\ \hspace{-1mm}, \ldots,
\underbrace{j_{h_{d_l,l}}, \ldots, j_{h_{d_l,l}}}_{n_{d_l,l}}\bigr\},
\end{equation}

\vspace{4mm}
\noindent
where
$n_{1,l}+n_{2,l}+\ldots+n_{d_l,l}=m_l;$\ \ $n_{1,l}, n_{2,l}, \ldots, n_{d_l,l}=1,\ldots, m_l;$\ \ 
$d_l=1,\ldots,m_l;$\ \ $l=1,\ldots,k.$ Note that the numbers $m_1,\ldots,m_k,$\ $g_1,\ldots,g_k$
depend on $(i_1,\ldots ,i_k)$ and the numbers
$n_{1,l},\ldots,n_{d_l,l},$\ $h_{1,l},\ldots,h_{d_l,l},$\ $d_l$
depend on $\{j_1,\ldots,j_k\}$. Moreover, 
$\left\{j_{g_1},\ldots,j_{g_k}\right\}
=\{j_1,\ldots,j_k\}.$

Using Theorem~9.6.9 \cite{Kuo} (also see 
\cite{ito1951}, Theorem~3.1), we get w.~p.~1

\vspace{2mm}
$$
J'\left[\phi_{j_{g_{m_1+m_2+\ldots+m_{l-1}+1}}}\ldots
\phi_{j_{g_{m_1+m_2+\ldots+m_l}}}
\right]_{T,t}^{(\hspace{0.5mm}\overbrace{{}_{i_l \ldots i_l}}^{m_l}\hspace{0.5mm})}
=
$$

\vspace{5mm}
\begin{equation}
\label{ziko20}
=\left\{
\begin{matrix}
H_{n_{1,l}}\left(\zeta_{j_{h_{1,l}}}^{(i_l)}\right)\ldots 
H_{n_{d_l,l}}\left(\zeta_{j_{h_{d_l,l}}}^{(i_l)}\right),\ 
&\hbox{\rm if}\ \ \ 
i_l\ne 0\cr\cr
\left(\zeta_{j_{h_{1,l}}}^{(0)}\right)^{n_{1,l}}\ldots
\left(\zeta_{j_{h_{d_l,l}}}^{(0)}\right)^{n_{d_l,l}},\  &\hbox{\rm if}\ \ \ 
i_l=0
\end{matrix}\right.\ \ \ (m_l>0),
\end{equation}

\vspace{6mm}
\noindent
where $H_n(x)$ is the Hermite polynomial of degree $n$

\vspace{1mm}
$$
H_n(x)=(-1)^n e^{x^2/2} \frac{d^n}{dx^n}\left(e^{-x^2/2}\right)
$$

\vspace{2mm}
\noindent
or

\vspace{-5mm}
\begin{equation}
\label{ziko500}
H_n(x)=n!\sum\limits_{m=0}^{[n/2]}\frac{(-1)^m x^{n-2m}}{m!(n-2m)! 2^m}\ \ \ (n\in\mathbb{N}),
\end{equation}

\vspace{4mm}
\noindent
and $\zeta_j^{(i)}$ $(i=0,1,\ldots,m,\ j=0,1,\ldots)$ is defined by (\ref{rr23}).

Note that the equality (\ref{ziko20}) is proved in \cite{11a}, \cite{new-2023a} using the Ito formula
(see the detailed proof in \cite{11a} (Sect.~1.14) or \cite{new-2023a}, Sect.~2.2).

From (\ref{ziko10}) and (\ref{ziko20}) we obtain w.~p.~1

\vspace{2mm}
$$
J'\left[\phi_{j_{g_{m_1+m_2+\ldots+m_{l-1}+1}}}\ldots
\phi_{j_{g_{m_1+m_2+\ldots+m_l}}}
\right]_{T,t}^{(\hspace{0.5mm}\overbrace{{}_{i_l \ldots i_l}}^{m_l}\hspace{0.5mm})}
=
$$

\vspace{5mm}
\begin{equation}
\label{ziko40}
={\bf 1}_{\{m_l=0\}}+{\bf 1}_{\{m_l>0\}}\left\{
\begin{matrix}
H_{n_{1,l}}\left(\zeta_{j_{h_{1,l}}}^{(i_l)}\right)\ldots 
H_{n_{d_l,l}}\left(\zeta_{j_{h_{d_l,l}}}^{(i_l)}\right),\ 
&\hbox{\rm if}\ \ \ 
i_l\ne 0\cr\cr
\left(\zeta_{j_{h_{1,l}}}^{(0)}\right)^{n_{1,l}}\ldots
\left(\zeta_{j_{h_{d_l,l}}}^{(0)}\right)^{n_{d_l,l}},\  &\hbox{\rm if}\ \ \ 
i_l=0
\end{matrix}\right.,
\end{equation}

\vspace{8mm}
\noindent
where ${\bf 1}_A$ denotes the indicator of the set $A$.

Using (\ref{leto6000xxa}), (\ref{ziko30}), and (\ref{ziko40}), we get w.~p.~1

\vspace{1mm}
$$
J'\left[\phi_{j_1}\ldots \phi_{j_k}\right]_{T,t}^{(i_1\ldots i_k)}=
$$

\vspace{2mm}
\begin{equation}
\label{ziko50}
=\prod_{l=1}^k\left({\bf 1}_{\{m_l=0\}}+{\bf 1}_{\{m_l>0\}}\left\{
\begin{matrix}
H_{n_{1,l}}\left(\zeta_{j_{h_{1,l}}}^{(i_l)}\right)\ldots 
H_{n_{d_l,l}}\left(\zeta_{j_{h_{d_l,l}}}^{(i_l)}\right),\ 
&\hbox{\rm if}\ \ \ 
i_l\ne 0\cr\cr
\left(\zeta_{j_{h_{1,l}}}^{(0)}\right)^{n_{1,l}}\ldots
\left(\zeta_{j_{h_{d_l,l}}}^{(0)}\right)^{n_{d_l,l}},\  &\hbox{\rm if}\ \ \ 
i_l=0
\end{matrix}\right.\ \right)=
\end{equation}

\vspace{7mm}
$$
=\prod_{l=1}^k\zeta_{j_l}^{(i_l)}
+\sum\limits_{r=1}^{[k/2]}
(-1)^r \times
$$

\vspace{2mm}
\begin{equation}
\label{chain401}
\times\sum_{\stackrel{(\{\{g_1, g_2\}, \ldots, 
\{g_{2r-1}, g_{2r}\}\}, \{q_1, \ldots, q_{k-2r}\})}
{{}_{\{g_1, g_2, \ldots, 
g_{2r-1}, g_{2r}, q_1, \ldots, q_{k-2r}\}=\{1, 2, \ldots, k\}}}}
\prod\limits_{s=1}^r
{\bf 1}_{\{i_{g_{{}_{2s-1}}}=~i_{g_{{}_{2s}}}\ne 0\}}
\Biggl.{\bf 1}_{\{j_{g_{{}_{2s-1}}}=~j_{g_{{}_{2s}}}\}}
\prod_{l=1}^{k-2r}\zeta_{j_{q_l}}^{(i_{q_l})}
\end{equation}

\vspace{5mm}
\noindent
w.~p.~1, where 
the multiple Wiener stochastic integral 
$J'[\phi_{j_1}\ldots \phi_{j_k}]_{T,t}^{(i_1\ldots i_k)}$ is defined by 
(\ref{WiI}); another
notations are the same as in (\ref{leto6000xxa}).

Thus, the following theorem is proved.

\vspace{2mm}
         
{\bf Theorem 3}\ \cite{11a} (Sect.~1.11), \cite{arxiv-1} (Sect.~15).
{\it Suppose that
the condition {\rm ($\star\star$)} is fulfilled
for the multi-index $(i_1 \ldots i_k)$
and the condition {\rm (\ref{ziko999})} is also 
fulfilled.
Furthermore$,$ let 
$\psi_l(\tau)\in L_2([t, T])$ $(l=$ $1,\ldots, k)$ and
$\{\phi_j(x)\}_{j=0}^{\infty}$ is an arbitrary complete orthonormal system  
of functions in the space $L_2([t,T]).$
Then the following expansions

\vspace{1mm}
$$
J[\psi^{(k)}]_{T,t}^{(i_1\ldots i_k)}=
\hbox{\vtop{\offinterlineskip\halign{
\hfil#\hfil\cr
{\rm l.i.m.}\cr
$\stackrel{}{{}_{p_1,\ldots,p_k\to \infty}}$\cr
}} }
\sum\limits_{j_1=0}^{p_1}\ldots
\sum\limits_{j_k=0}^{p_k}
C_{j_k\ldots j_1}\times
$$

\vspace{5mm}
\begin{equation}
\label{ddddxxxx1111}
\times
\prod_{l=1}^k\left({\bf 1}_{\{m_l=0\}}+{\bf 1}_{\{m_l>0\}}\left\{
\begin{matrix}
H_{n_{1,l}}\left(\zeta_{j_{h_{1,l}}}^{(i_l)}\right)\ldots 
H_{n_{d_l,l}}\left(\zeta_{j_{h_{d_l,l}}}^{(i_l)}\right),\ 
&\hbox{\rm if}\ \ \ 
i_l\ne 0\cr\cr
\left(\zeta_{j_{h_{1,l}}}^{(0)}\right)^{n_{1,l}}\ldots
\left(\zeta_{j_{h_{d_l,l}}}^{(0)}\right)^{n_{d_l,l}},\  &\hbox{\rm if}\ \ \ 
i_l=0
\end{matrix}\right.\ \right),
\end{equation}

\vspace{8mm}
$$
J[\psi^{(k)}]_{T,t}^{(i_1\ldots i_k)}=
\hbox{\vtop{\offinterlineskip\halign{
\hfil#\hfil\cr
{\rm l.i.m.}\cr
$\stackrel{}{{}_{p_1,\ldots,p_k\to \infty}}$\cr
}} }
\sum\limits_{j_1=0}^{p_1}\ldots
\sum\limits_{j_k=0}^{p_k}
C_{j_k\ldots j_1}\Biggl(
\prod_{l=1}^k\zeta_{j_l}^{(i_l)}+\sum\limits_{r=1}^{[k/2]}
(-1)^r \times
\Biggr.
$$

\vspace{3mm}
\begin{equation}
\label{chainxx90}
\times
\sum_{\stackrel{(\{\{g_1, g_2\}, \ldots, 
\{g_{2r-1}, g_{2r}\}\}, \{q_1, \ldots, q_{k-2r}\})}
{{}_{\{g_1, g_2, \ldots, 
g_{2r-1}, g_{2r}, q_1, \ldots, q_{k-2r}\}=\{1, 2, \ldots, k\}}}}
\prod\limits_{s=1}^r
{\bf 1}_{\{i_{g_{{}_{2s-1}}}=~i_{g_{{}_{2s}}}\ne 0\}}
\Biggl.{\bf 1}_{\{j_{g_{{}_{2s-1}}}=~j_{g_{{}_{2s}}}\}}
\prod_{l=1}^{k-2r}\zeta_{j_{q_l}}^{(i_{q_l})}\Biggr)
\end{equation}

\vspace{5mm}
\noindent
con\-verg\-ing in the mean-square sense are valid$,$
where $[x]$ is an integer part of a real number $x,$\ \
$n_{1,l}+n_{2,l}+\ldots+n_{d_l,l}=m_l,$\ \ $n_{1,l}, n_{2,l}, \ldots, n_{d_l,l}=1,\ldots, m_l,$\ \ 
$d_l=1,\ldots,m_l,$\ \ $l=1,\ldots,k;$\ \ $m_1+\ldots+m_k=k,$\ 
the numbers $m_1,\ldots,m_k,$\ $g_1,\ldots,g_k$
depend on $(i_1,\ldots,i_k)$ and 
the numbers $n_{1,l},\ldots,n_{d_l,l},$\ $h_{1,l},\ldots,h_{d_l,l},$\ $d_l$
depend on $\{j_1,\ldots,j_k\};$ moreover$,$ $\left\{j_{g_1},\ldots,j_{g_k}\right\}
=\{j_1,\ldots,j_k\};$
$H_n(x)$ is the Hermite polynomial {\rm (\ref{ziko500});}
another
notations are the same as in {\rm (\ref{leto6000xxa})} and in Theorem {\rm 1}.}

\vspace{2mm}

It should be noted that an analogue of the expansion
(\ref{ddddxxxx1111}) was considered 
in \cite{Rybakov1000}. 
However, the proof of an analogue of the expansion (\ref{ddddxxxx1111})
from \cite{Rybakov1000} is different from the proof given in this section
(see \cite{new-2023a}, Sect.~4 for details).

\vspace{5mm}

\section{Exact Calculation of the Mean-Square Approximation Error in The Method of 
Generalized Multiple Fourier Series. The Case 
of an Arbitrary Complete Orthonormal Systems  
of Functions in the Space $L_2([t,T])$ 
and Weight Functions $\psi_1(\tau),\ldots,\psi_k(\tau)\in L_2([t, T])$}

\vspace{5mm}

In this section, we generalize Theorem~2 to the case 
of an arbitrary complete orthonormal systems  
of functions in the Space $L_2([t,T])$ 
and $\psi_1(\tau),\ldots,\psi_k(\tau)\in L_2([t, T])$.

\vspace{2mm}

{\bf Theorem 4}\ \cite{arxiv-1}.
{\it Suppose that
$\psi_1(\tau),\ldots,\psi_k(\tau)\in L_2([t, T])$ 
and
$\{\phi_j(x)\}_{j=0}^{\infty}$ is an arbitrary complete orthonormal system  
of functions in the space $L_2([t,T]).$ 
Then

\vspace{1mm}
$$
{\sf M}\left\{\left(J[\psi^{(k)}]_{T,t}-
J[\psi^{(k)}]_{T,t}^p\right)^2\right\}
= \int\limits_{[t,T]^k} K^2(t_1,\ldots,t_k)
dt_1\ldots dt_k - 
$$

\vspace{1mm}
\begin{equation}
\label{chain100}
- \sum_{j_1=0}^{p}\ldots\sum_{j_k=0}^{p}
C_{j_k\ldots j_1}
{\sf M}\left\{J[\psi^{(k)}]_{T,t}
\sum\limits_{(j_1,\ldots,j_k)}
\int\limits_t^T \phi_{j_k}(t_k)
\ldots
\int\limits_t^{t_{2}}\phi_{j_{1}}(t_{1})
d{\bf f}_{t_1}^{(i_1)}\ldots
d{\bf f}_{t_k}^{(i_k)}\right\},
\end{equation}

\vspace{5mm}
\noindent
where
$$
J[\psi^{(k)}]_{T,t}=\int\limits_t^T\psi_k(t_k) \ldots \int\limits_t^{t_{2}}
\psi_1(t_1) d{\bf f}_{t_1}^{(i_1)}\ldots
d{\bf f}_{t_k}^{(i_k)},
$$

\vspace{2mm}
\begin{equation}
\label{chain101}
J[\psi^{(k)}]_{T,t}^p=
\sum_{j_1=0}^{p}\ldots\sum_{j_k=0}^{p}
C_{j_k\ldots j_1} J'[\phi_{j_1}\ldots \phi_{j_k}]_{T,t}^{(i_1\ldots i_k)},
\end{equation}

\vspace{5mm}

\noindent
$J'[\phi_{j_1}\ldots \phi_{j_k}]_{T,t}^{(i_1\ldots i_k)}$
is the multiple Wiener stochastic integral 
defined by {\rm (\ref{WiI}),}
the Fourier coefficient $C_{j_k\ldots j_1}$ has the form {\rm (\ref{chain300}),}
$K(t_1,\ldots,t_k)$ is defined by {\rm (\ref{chain200}),}

\vspace{-1mm}
$$
\zeta_{j}^{(i)}=
\int\limits_t^T \phi_{j}(s) d{\bf f}_s^{(i)}
$$

\vspace{2mm}
\noindent
are independent standard Gaussian random variables
for various
$i$ or $j$ $(i=1,\ldots,m),$

\vspace{-1mm}
$$
\sum\limits_{(j_1,\ldots,j_k)}
$$ 

\vspace{2mm}
\noindent
means the sum with respect to all
possible permutations 
$(j_1,\ldots,j_k).$ At the same time if 
$j_r$ swapped with $j_q$ in the permutation $(j_1,\ldots,j_k)$,
then $i_r$ swapped with $i_q$ in the permutation
$(i_1,\ldots,i_k)$ {\rm (}see {\rm (\ref{chain100})).}}

\vspace{2mm}

{\bf Proof.}\ First, note that the formula (\ref{chain101}) 
appears due to (\ref{chain102xx}).
Using the equality (\ref{Wi110}), we get

\vspace{-2mm}
\begin{equation}
\label{chain103}
J'[\phi_{j_1}\ldots \phi_{j_k}]_{T,t}^{(i_1\ldots i_k)}=\sum_{(t_1,\ldots,t_k)}
\int\limits_t^T \phi_{j_k}(t_k)
\ldots
\int\limits_t^{t_{2}}\phi_{j_{1}}(t_{1})
d{\bf f}_{t_1}^{(i_1)}\ldots
d{\bf f}_{t_k}^{(i_k)}\ \ \ \hbox{w.\ p.\ 1},
\end{equation}

\vspace{3mm}
\noindent
where permutations $(t_1,\ldots,t_k)$ when summing are 
performed only in the values
$d{\bf f}_{t_1}^{(i_1)}
\ldots $
$d{\bf f}_{t_k}^{(i_k)}.$ At the same time the indices near 
upper 
limits of integration in the iterated stochastic integrals are changed 
correspondently and if $t_r$ swapped with $t_q$ in the  
permutation $(t_1,\ldots,t_k)$, then $i_r$ swapped with $i_q$ in 
the permutation $(i_1,\ldots,i_k).$

It is easy to see that the equality (\ref{chain103}) can be written in the form

\begin{equation}
\label{chain104}
J'[\phi_{j_1}\ldots \phi_{j_k}]_{T,t}^{(i_1\ldots i_k)}=
\sum\limits_{(j_1,\ldots,j_k)}
\int\limits_t^T \phi_{j_k}(t_k)
\ldots
\int\limits_t^{t_{2}}\phi_{j_{1}}(t_{1})
d{\bf f}_{t_1}^{(i_1)}\ldots
d{\bf f}_{t_k}^{(i_k)}\ \ \ \hbox{w.\ p.\ 1},
\end{equation}

\vspace{3mm}
\noindent
where 

\vspace{-1mm}
$$
\sum\limits_{(j_1,\ldots,j_k)}
$$ 

\vspace{2mm}
\noindent
means the sum with respect to all
possible permutations 
$(j_1,\ldots,j_k).$ At the same time if 
$j_r$ swapped with $j_q$ in the permutation $(j_1,\ldots,j_k)$,
then $i_r$ swapped with $i_q$ in the permutation.

Further proof of Theorem~4 is based on the equality 
(\ref{chain104}) and is similar to the proof of Theorem 2.
Theorem~4 is proved.

\vspace{5mm}

\section{Estimate for the Mean-Square Approximation Error in the Method
of Generalized Multiple Fourier Series}

\vspace{5mm}

In this section, we prove the useful estimate for the
mean-square error of approximation based on Theorem 3.

\vspace{2mm}

{\bf Theorem 5.}\ {\it Suppose that
$\psi_1(\tau),\ldots,\psi_k(\tau)\in L_2([t, T])$ 
and
$\{\phi_j(x)\}_{j=0}^{\infty}$ is an arbitrary complete orthonormal system  
of functions in the space $L_2([t,T]).$ 
Then the estimate

\vspace{1mm}
$$
{\sf M}\left\{\left(
J[\psi^{(k)}]_{T,t}-J[\psi^{(k)}]_{T,t}^{p_1,\ldots,p_k}
\right)^2\right\}
\le 
$$

\begin{equation}
\label{z1}
\le k!\left(\int\limits_{[t,T]^k}
K^2(t_1,\ldots,t_k)
dt_1\ldots dt_k -\sum_{j_1=0}^{p_1}\ldots
\sum_{j_k=0}^{p_k}C^2_{j_k\ldots j_1}\right)
\end{equation}

\vspace{4mm}
\noindent
is valid for the following cases{\rm :}

\vspace{3mm}

{\rm 1.}\ $i_1,\ldots,i_k=1,\ldots,m$\ \ and\ \ $0<T-t<\infty,$

\vspace{2mm}

{\rm 2.}\ $i_1,\ldots,i_k=0, 1,\ldots,m,$\ \ $i_1^2+\ldots+i_k^2>0,$\ \
and\ \ $0<T-t<1,$

\vspace{4mm}
\noindent
where $J[\psi^{(k)}]_{T,t}$ is the iterated Ito stochastic integral {\rm (\ref{sodom20}),}
$J[\psi^{(k)}]_{T,t}^{p_1,\ldots,p_k}$ is the expression
on the right-hand side of {\rm (\ref{chainxx90})} before passing to the limit
$\hbox{\vtop{\offinterlineskip\halign{
\hfil#\hfil\cr
{\rm l.i.m.}\cr
$\stackrel{}{{}_{p_1,\ldots,p_k\to \infty}}$\cr
}} };$ another
notations are the same as in Theorem {\rm 3}.
}

\vspace{2mm}

{\bf Proof.}\ Using (\ref{chain102}), (\ref{chain8810}), (\ref{chain102xx}), Theorem~3, 
orthonormality of the system $\{\phi_j(x)\}_{j=0}^{\infty}$, 
and the elementary inequality
\begin{equation}
\label{y5}
\left(a_1+a_2+\ldots+a_n\right)^2 \le
n\left(a_1^2+a_2^2+\ldots+a_n^2\right),
\end{equation}

\vspace{3mm}
\noindent
we obtain for the case $i_1,\ldots,i_k=1,\dots,m$
$(0<T-t<\infty)$
the following estimate 

\vspace{1mm}
$$
{\sf M}\left\{\left(
J[\psi^{(k)}]_{T,t}-J[\psi^{(k)}]_{T,t}^{p_1,\ldots,p_k}
\right)^2\right\}
\le 
$$

\vspace{2mm}
$$
\le k!
\sum_{(t_1,\ldots,t_k)}
\int\limits_{t}^{T}
\ldots
\int\limits_{t}^{t_2}
\left(K(t_1,\ldots,t_k)-
\sum_{j_1=0}^{p_1}\ldots
\sum_{j_k=0}^{p_k}
C_{j_k\ldots j_1}
\prod_{l=1}^k\phi_{j_l}(t_l)\right)^2
dt_1
\ldots
dt_k=
$$

\vspace{3mm}
$$
=k!\int\limits_{[t,T]^k}
\left(K(t_1,\ldots,t_k)-
\sum_{j_1=0}^{p_1}\ldots
\sum_{j_k=0}^{p_k}
C_{j_k\ldots j_1}
\prod_{l=1}^k\phi_{j_l}(t_l)\right)^2
dt_1
\ldots
dt_k
$$

\vspace{2mm}

$$
= k!\left(\int\limits_{[t,T]^k}
K^2(t_1,\ldots,t_k)
dt_1\ldots dt_k -\sum_{j_1=0}^{p_1}\ldots
\sum_{j_k=0}^{p_k}C^2_{j_k\ldots j_1}\right).
$$

\vspace{5mm}

Similarly using standard moment properties of stochastic integrals,
we get 

\vspace{1mm}
$$
{\sf M}\left\{\left(
J[\psi^{(k)}]_{T,t}-J[\psi^{(k)}]_{T,t}^{p_1,\ldots,p_k}
\right)^2\right\}
\le 
$$

\vspace{2mm}

$$
= C_k\left(\int\limits_{[t,T]^k}
K^2(t_1,\ldots,t_k)
dt_1\ldots dt_k -\sum_{j_1=0}^{p_1}\ldots
\sum_{j_k=0}^{p_k}C^2_{j_k\ldots j_1}\right),
$$

\vspace{5mm}
\noindent
where $i_1,\ldots,i_k=0, 1,\ldots,m$
$(i_1^2+\ldots+i_k^2>0),$ 
and $C_k$ is a constant.

It is not difficult to see that $C_k$ depends on 
$k$ ($k$ is the multiplicity 
of the iterated Ito stochastic integral) and $T-t$ ($T-t$ is the length
of integration interval for the iterated Ito stochastic integral).
Moreover, $C_k$ has the following form

$$
C_k=k!\cdot{\rm max}\biggl\{
(T-t)^{\alpha_1},\ (T-t)^{\alpha_2},\ \ldots,\ (T-t)^{\alpha_{k!}}
\biggr\},
$$

\vspace{3mm}
\noindent
where $\alpha_1, \alpha_2, \ldots, \alpha_{k!}=0,\ 1,\ldots,\ k-1.$

Then for the case
$i_1,\ldots,i_k=0, 1,\ldots,m,$ $i_1^2+\ldots+i_k^2>0$
$(0<T-t<1)$
we obtain (\ref{z1}). Theorem 5 is proved.

\vspace{2mm}

{\bf Example 3.}\ Let us consider the estimate (\ref{z1})
for the iterated Ito stochastic integral $I_{(000)T,t}^{(i_1i_2 i_3)}$ defined by
(\ref{k1001})

$$
{\sf M}\left\{\left(
I_{(000)T,t}^{(i_1i_2 i_3)}-
I_{(000)T,t}^{(i_1i_2 i_3)p}\right)^2\right\}\le
6\left(\frac{(T-t)^{3}}{6}-\sum_{j_1,j_2,j_3=0}^{p}
C_{j_3j_2j_1}^2\right)\ \ \ (i_1, i_2, i_3=1,\ldots,m),
$$

\vspace{5mm}
\noindent
where $C_{j_3j_2j_1}$ has the form (\ref{w1}).

\vspace{5mm}

\section{Proof of Convergence With Probability 1 in the Method of Generalized 
Multiple Fourier Series. The Cases of Complete Orthonormal Systems of Legendre
Polynomials and Trigonometric Functions in the Space $L_2([t, T])$}

\vspace{5mm}

This section is written on the base of 
\cite{11a} (Sect.~1.7.2), \cite{new-new-2}.
Remind that in a lot of author's publications \cite{1}-\cite{new-2023a}
the convergence in Theorem 1 has been considered in different
probabilistic
senses. For example, the mean-square convergence \cite{1} (2006) (also see
\cite{2}-\cite{11a}) and convergence in the mean of degree
$2n$ $(n\in\mathbb{N})$ 
\cite{11a} (Sect.~1.1.9, 1.11, 1.12), \cite{arxiv-1} (Sect.~6, 15, 16)
have been proved. On the examples
of specific iterated Ito stochastic integrals of mutiplicities 1 and 2
the convergence with probability 1 has been considered in 
\cite{3} (2007) (also see \cite{4}-\cite{12},
\cite{arxiv-4}). However, these examples are narrow particular 
cases of the iterated Ito stochastic integrals (\ref{sodom20}).

In this section, we formulate and prove the theorem 
\cite{11a} (Sect.~1.7.2), \cite{new-new-2} on 
convergence with probability 1 of the expansions 
of iterated Ito stochastic integrals 
of multiplicity $k$ $(k\in\mathbb{N})$
from Theorems 1, 3.

Let us remind the well-known fact from the mathematical analysis
which is connected to existence
of iterated limits.

\vspace{2mm}

{\bf Proposition 1.}\ {\it Let $\bigl\{x_{n,m}\bigr\}_{n,m=1}^{\infty}$
be a double sequence and let there exists the limit

\vspace{-1mm}
$$
\lim\limits_{n,m\to\infty}x_{n,m}=a<\infty.
$$

\vspace{2mm}

Moreover, let there exist the limits

\vspace{-1mm}
$$
\lim\limits_{n\to\infty}x_{n,m}<\infty\ \ \ \hbox{for all}\ \ \ m,\ \ \ \
\lim\limits_{m\to\infty}x_{n,m}<\infty\ \ \ \hbox{for all}\ \ \ n.
$$

\vspace{2mm}

Then there exist the iterated limits

\vspace{-1mm}
$$
\lim\limits_{n\to\infty}\lim\limits_{m\to\infty}x_{n,m},\ \ \ 
\lim\limits_{m\to\infty}\lim\limits_{n\to\infty}x_{n,m}
$$

\vspace{2mm}
and moreover,

$$
\lim\limits_{n\to\infty}\lim\limits_{m\to\infty}x_{n,m}=
\lim\limits_{m\to\infty}\lim\limits_{n\to\infty}x_{n,m}=a.
$$
}

\vspace{3mm}

{\bf Theorem 6} \cite{11a}-\cite{11aaa-a}, \cite{arxiv-1}, \cite{arxiv-3}, \cite{arxiv-4}. {\it Let 
$\psi_l(\tau)$ $(l=1,\ldots, k)$ are 
continuously differentiable nonrandom functions at the interval
$[t, T]$ and $\{\phi_j(x)\}_{j=0}^{\infty}$ is a complete
orthonormal system of Legendre polynomials or 
trigonometric functions in the space $L_2([t, T]).$
Then 

$$
J[\psi^{(k)}]_{T,t}^{p,\ldots,p}\ \to \ J[\psi^{(k)}]_{T,t}\ \ \ 
\hbox{if}\ \ \ p\to \infty
$$

\vspace{3mm}
\noindent
w.~p.~{\rm 1,} where $J[\psi^{(k)}]_{T,t}^{p,\ldots,p}$
is the expression on the right-hand sides of {\rm (\ref{tyyyxxx})}
and {\rm (\ref{leto6000})}
before passing to the limit 
$\hbox{\vtop{\offinterlineskip\halign{
\hfil#\hfil\cr
{\rm l.i.m.}\cr
$\stackrel{}{{}_{p_1,\ldots,p_k\to \infty}}$\cr
}} }$  for the case $p_1=\ldots=p_k=p,$ i.e.

\vspace{2mm}
$$
J[\psi^{(k)}]_{T,t}^{p,\ldots,p}=
\sum_{j_1=0}^{p}\ldots\sum_{j_k=0}^{p}
C_{j_k\ldots j_1}\Biggl(
\prod_{l=1}^k\zeta_{j_l}^{(i_l)}\ -
\Biggr.
$$

\vspace{4mm}
$$
-\ \Biggl.
\hbox{\vtop{\offinterlineskip\halign{
\hfil#\hfil\cr
{\rm l.i.m.}\cr
$\stackrel{}{{}_{N\to \infty}}$\cr
}} }\sum_{(l_1,\ldots,l_k)\in {\rm G}_k}
\phi_{j_{1}}(\tau_{l_1})
\Delta{\bf w}_{\tau_{l_1}}^{(i_1)}\ldots
\phi_{j_{k}}(\tau_{l_k})
\Delta{\bf w}_{\tau_{l_k}}^{(i_k)}\Biggr)
$$

\vspace{5mm}
\noindent
or

$$
J[\psi^{(k)}]_{T,t}^{p,\ldots,p}=
\sum\limits_{j_1=0}^{p}\ldots
\sum\limits_{j_k=0}^{p}
C_{j_k\ldots j_1}\Biggl(
\prod_{l=1}^k\zeta_{j_l}^{(i_l)}+\sum\limits_{r=1}^{[k/2]}
(-1)^r \times
\Biggr.
$$

\vspace{5mm}
$$
\times
\sum_{\stackrel{(\{\{g_1, g_2\}, \ldots, 
\{g_{2r-1}, g_{2r}\}\}, \{q_1, \ldots, q_{k-2r}\})}
{{}_{\{g_1, g_2, \ldots, 
g_{2r-1}, g_{2r}, q_1, \ldots, q_{k-2r}\}=\{1, 2, \ldots, k\}}}}
\prod\limits_{s=1}^r
{\bf 1}_{\{i_{g_{{}_{2s-1}}}=~i_{g_{{}_{2s}}}\ne 0\}}
\Biggl.{\bf 1}_{\{j_{g_{{}_{2s-1}}}=~j_{g_{{}_{2s}}}\}}
\prod_{l=1}^{k-2r}\zeta_{j_{q_l}}^{(i_{q_l})}\Biggr),
$$

\vspace{7mm}
\noindent
where $i_1,\ldots,i_k=1,\ldots,m$.}

\vspace{2mm}

{\bf Proof.} Let us consider the Parseval equality

\vspace{1mm}
\begin{equation}
\label{par1}
\int\limits_{[t,T]^k}K^2(t_1,\ldots,t_k)dt_1\ldots dt_k=
\lim\limits_{p_1,\ldots,p_k\to\infty}
\sum_{j_1=0}^{p_1}\ldots \sum_{j_k=0}^{p_k}
C_{j_k\ldots j_1}^2,
\end{equation}

\vspace{4mm}
\noindent
where

\begin{equation}
\label{pppx}
K(t_1,\ldots,t_k)=
\begin{cases}
\psi_1(t_1)\ldots \psi_k(t_k),\ &t_1<\ldots<t_k\\
~\\
~\\
0,\ &\hbox{\rm otherwise}
\end{cases}\ \ \ \ 
=\ \ \ \ 
\prod\limits_{l=1}^k
\psi_l(t_l)\ \prod\limits_{l=1}^{k-1}{\bf 1}_{\{t_l<t_{l+1}\}},\ 
\end{equation}

\vspace{5mm}
\noindent
where $t_1,\ldots,t_k\in [t, T]$ for $k\ge 2$ and 
$K(t_1)\equiv\psi_1(t_1)$ for $t_1\in[t, T],$ 
${\bf 1}_A$ denotes the indicator of the set $A$,

\vspace{-2mm}
\begin{equation}
\label{ppppax}
C_{j_k\ldots j_1}=\int\limits_{[t,T]^k}
K(t_1,\ldots,t_k)\prod_{l=1}^{k}\phi_{j_l}(t_l)dt_1\ldots dt_k
\end{equation}

\vspace{4mm}
\noindent
is the Fourier coefficient.

Using (\ref{pppx}), we obtain

$$
C_{j_k\ldots j_1}=
\int\limits_t^T
\phi_{j_k}(t_k)\psi_k(t_k)\ldots \int\limits_t^{t_2}
\phi_{j_1}(t_1)\psi_1(t_1)dt_1\ldots dt_k.
$$

\vspace{3mm}

Further, we denote

\vspace{-1mm}
$$
\lim\limits_{p_1,\ldots,p_k\to\infty}
\sum_{j_1=0}^{p_1}\ldots \sum_{j_k=0}^{p_k}
C_{j_k\ldots j_1}^2\stackrel{\sf def}{=}
\sum_{j_1,\ldots,j_k=0}^{\infty}
C_{j_k\ldots j_1}^2.
$$

\vspace{4mm}

If $p_1=\ldots=p_k=p,$ then we also write

\vspace{1mm}
$$
\lim\limits_{p\to\infty}
\sum_{j_1=0}^{p}\ldots \sum_{j_k=0}^{p}
C_{j_k\ldots j_1}^2\stackrel{\sf def}{=}
\sum_{j_1,\ldots,j_k=0}^{\infty}
C_{j_k\ldots j_1}^2.
$$

\vspace{4mm}

From the other hand, for iterated limits we write

\vspace{1mm}
$$
\lim\limits_{p_1\to\infty}\ldots \lim\limits_{p_k\to\infty}
\sum_{j_1=0}^{p_1}\ldots \sum_{j_k=0}^{p_k}
C_{j_k\ldots j_1}^2\stackrel{\sf def}{=}
\sum_{j_1=0}^{\infty}\ldots
\sum_{j_k=0}^{\infty}
C_{j_k\ldots j_1}^2,
$$

\vspace{3mm}
$$
\lim\limits_{p_1\to\infty}\lim\limits_{p_2,\ldots,p_k\to\infty}
\sum_{j_1=0}^{p_1}\ldots \sum_{j_k=0}^{p_k}
C_{j_k\ldots j_1}^2\stackrel{\sf def}{=}
\sum_{j_1=0}^{\infty}
\sum_{j_2,\ldots,j_k=0}^{\infty}
C_{j_k\ldots j_1}^2
$$

\vspace{5mm}
\noindent
and so on.

\vspace{2mm}

{\bf Lemma 1.}\ {\it The following equalities are fulfilled

\vspace{1mm}
$$
\sum_{j_1,\ldots,j_k=0}^{\infty}
C_{j_k\ldots j_1}^2=
\sum_{j_1=0}^{\infty}\ldots
\sum_{j_k=0}^{\infty}
C_{j_k\ldots j_1}^2=
$$

\vspace{2mm}
\begin{equation}
\label{lem1}
=\sum_{j_k=0}^{\infty}\ldots
\sum_{j_1=0}^{\infty}
C_{j_k\ldots j_1}^2=
\sum_{j_{q_1}=0}^{\infty}\ldots
\sum_{j_{q_k}=0}^{\infty}
C_{j_k\ldots j_1}^2
\end{equation}

\vspace{5mm}
\noindent
for any permutation $(q_1,\ldots,q_k)$ such that
$\{q_1,\ldots,q_k\}=\{1,\ldots,k\}.$}

\vspace{2mm}

{\bf Proof.} Let us consider the value

\vspace{-3mm}
\begin{equation}
\label{21}
\sum_{j_{q_l}=0}^{p}\ldots
\sum_{j_{q_k}=0}^{p}
C_{j_k\ldots j_1}^2
\end{equation}

\vspace{3mm}
\noindent
for any permutation $(q_l,\ldots,q_k)$, where $l=1,2,\ldots,k$,
$\{q_1,\ldots,q_k\}=\{1,\ldots,k\}.$

Obviously, (\ref{21}) 
is the non-decreasing sequence with respect to $p$.
Moreover,

\vspace{2mm}
$$
\sum_{j_{q_l}=0}^{p}\ldots
\sum_{j_{q_k}=0}^{p}
C_{j_k\ldots j_1}^2\le 
\sum_{j_{q_1}=0}^{p}\sum_{j_{q_2}=0}^{p}\ldots
\sum_{j_{q_k}=0}^{p}
C_{j_k\ldots j_1}^2\le 
$$

\vspace{3mm}
$$
\le
\sum_{j_1,\ldots,j_k=0}^{\infty}
C_{j_k\ldots j_1}^2<\infty.
$$

\vspace{3mm}

Then the following limit

\vspace{1mm}
$$
\lim\limits_{p\to\infty}\sum\limits_{j_{q_l}=0}^p \ldots 
\sum\limits_{j_{q_k}=0}^{p}
C_{j_k\ldots j_1}^2=
\sum_{j_{q_l},\ldots,j_{q_k}=0}^{\infty}
C_{j_k\ldots j_1}^2
$$

\vspace{3mm}
\noindent
exists.

Let $p_l,\ldots,p_k$ simultaneously tend to infinity.
Then $g, r\to \infty$, where $g=\min\{p_l,\ldots,p_k\}$ and
$r=\max\{p_l,\ldots,p_k\}$. Moreover,

\vspace{2mm}
$$
\sum_{j_{q_l}=0}^{g}\ldots
\sum_{j_{q_k}=0}^{g}
C_{j_k\ldots j_1}^2\le 
\sum_{j_{q_l}=0}^{p_l}\ldots
\sum_{j_{q_k}=0}^{p_k}
C_{j_k\ldots j_1}^2\le
\sum_{j_{q_l}=0}^{r}\ldots
\sum_{j_{q_k}=0}^{r}
C_{j_k\ldots j_1}^2.
$$

\vspace{5mm}

This means that the existence of the limit 

\vspace{1mm}
\begin{equation}
\label{1c1c}
\lim\limits_{p\to\infty}\sum_{j_{q_l}=0}^{p}\ldots
\sum_{j_{q_k}=0}^{p}
C_{j_k\ldots j_1}^2
\end{equation}

\vspace{3mm}
\noindent
implies the existence of the limit 

\vspace{-1mm}

\begin{equation}
\label{1d1d}
\lim\limits_{p_l,\ldots,p_k\to\infty}\sum_{j_{q_l}=0}^{p_l}\ldots
\sum_{j_{q_k}=0}^{p_k}
C_{j_k\ldots j_1}^2
\end{equation}

\vspace{3mm}
\noindent
and equality of limits (\ref{1c1c}) and (\ref{1d1d}).
 
Taking into account the above reasoning, we have 

\vspace{1mm}
$$
\lim\limits_{p,q\to\infty}\sum_{j_{q_l}=0}^{q}\sum_{j_{q_{l+1}}=0}^{p}\ldots
\sum_{j_{q_k}=0}^{p}
C_{j_k\ldots j_1}^2=
\lim\limits_{p\to\infty}\sum_{j_{q_l}=0}^{p}\ldots
\sum_{j_{q_k}=0}^{p}
C_{j_k\ldots j_1}^2=
$$

\vspace{2mm}
\begin{equation}
\label{1h1h}
=\lim\limits_{p_l,\ldots,p_k\to\infty}\sum_{j_{q_l}=0}^{p_l}\ldots
\sum_{j_{q_k}=0}^{p_k}
C_{j_k\ldots j_1}^2.
\end{equation}

\vspace{5mm}

Since the limit
$$
\sum_{j_1,\ldots,j_k=0}^{\infty}
C_{j_k\ldots j_1}^2
$$

\vspace{3mm}
\noindent
exists (see the Parseval equality (\ref{par1})), then from Proposition 1
we have

\vspace{1mm}
$$
\sum_{j_{q_1}=0}^{\infty}\sum_{j_{q_2},\ldots,j_{q_k}=0}^{\infty}
C_{j_k\ldots j_1}^2=
\lim\limits_{q\to\infty}
\lim\limits_{p\to\infty}
\sum_{j_{q_1}=0}^{q}\sum_{j_{q_2}=0}^p \ldots \sum_{j_{q_k}=0}^{p}
C_{j_k\ldots j_1}^2=
$$

\vspace{3mm}
\begin{equation}
\label{1b1b}
=\lim\limits_{q,p\to\infty}
\sum_{j_{q_1}=0}^{q}\sum_{j_{q_2}=0}^p \ldots \sum_{j_{q_k}=0}^{p}
C_{j_k\ldots j_1}^2=
\sum_{j_1,\ldots,j_k=0}^{\infty}
C_{j_k\ldots j_1}^2.
\end{equation}

\vspace{5mm}

Using (\ref{1h1h}) and Proposition 1, we get

\vspace{1mm}
$$
\sum_{j_{q_2}=0}^{\infty}\sum_{j_{q_3},\ldots,j_{q_k}=0}^{\infty}
C_{j_k\ldots j_1}^2=
\lim\limits_{q\to\infty}
\lim\limits_{p\to\infty}
\sum_{j_{q_2}=0}^{q}\sum_{j_{q_3}=0}^p \ldots \sum_{j_{q_k}=0}^{p}
C_{j_k\ldots j_1}^2=
$$

\vspace{3mm}
\begin{equation}
\label{1a1a}
=\lim\limits_{q,p\to\infty}
\sum_{j_{q_2}=0}^{q}\sum_{j_{q_3}=0}^p \ldots \sum_{j_{q_k}=0}^{p}
C_{j_k\ldots j_1}^2=
\sum_{j_{q_2},\ldots,j_{q_k}=0}^{\infty}
C_{j_k\ldots j_1}^2.
\end{equation}

\vspace{5mm}

Combining (\ref{1a1a}) and (\ref{1b1b}), we obtain

\vspace{1mm}
$$
\sum_{j_{q_1}=0}^{\infty}\sum_{j_{q_2}=0}^{\infty}
\sum_{j_{q_3},\ldots,j_{q_k}=0}^{\infty}
C_{j_k\ldots j_1}^2=
\sum_{j_{1},\ldots,j_{k}=0}^{\infty}
C_{j_k\ldots j_1}^2.
$$

\vspace{4mm}

Repeating the above reasoning, we complete the proof of Lemma
1. 

Further, let us show that for $s=1,\ldots,k$

\vspace{1mm}
$$
\sum_{j_1=0}^{\infty}\ldots
\sum_{j_{s-1}=0}^{\infty}
\sum_{j_s=p+1}^{\infty}\sum_{j_{s+1}=0}^{\infty}\ldots \sum_{j_k=0}^{\infty}
C_{j_k\ldots j_1}^2=
$$

\vspace{3mm}
\begin{equation}
\label{d11}
=
\sum_{j_s=p+1}^{\infty}\sum_{j_{s-1}=0}^{\infty}\ldots
\sum_{j_{1}=0}^{\infty}
\sum_{j_{s+1}=0}^{\infty}\ldots \sum_{j_k=0}^{\infty}
C_{j_k\ldots j_1}^2.
\end{equation}

\vspace{5mm}

Using the arguments which we used when proving Lemma 1, we have

\vspace{2mm}
$$
\lim\limits_{n\to\infty}
\sum_{j_1=0}^{n}\ldots
\sum_{j_{s-1}=0}^{n}
\sum_{j_s=0}^{p}\sum_{j_{s+1}=0}^{n}\ldots \sum_{j_k=0}^{n}
C_{j_k\ldots j_1}^2=
$$

\vspace{3mm}
\begin{equation}
\label{ura0}
=\sum_{j_s=0}^{p}\ \sum_{j_{1},\ldots, j_{s-1}, j_{s+1},\ldots,j_k=0}^{\infty}
C_{j_k\ldots j_1}^2
=\sum_{j_s=0}^{p}\sum_{j_{q_1}=0}^{\infty}\ldots
\sum_{j_{q_{k-1}}=0}^{\infty}
C_{j_k\ldots j_1}^2
\end{equation}

\vspace{5mm}
\noindent
for any permutation $(q_1,\ldots,q_{k-1})$ such that
$\{q_1,\ldots,q_{k-1}\}=\{1,\ldots,s-1,s+1,\ldots,k\}$,
where $p$ is a fixed natural number.

Obviously, we have

\vspace{2mm}
$$
\sum_{j_s=0}^{p}\sum_{j_{q_1}=0}^{\infty}\ldots
\sum_{j_{q_{k-1}}=0}^{\infty}
C_{j_k\ldots j_1}^2=
\sum_{j_{q_1}=0}^{\infty}\ldots \sum_{j_s=0}^{p} \ldots
\sum_{j_{q_{k-1}}=0}^{\infty}C_{j_k\ldots j_1}^2 = \ldots =
$$

\vspace{3mm}
\begin{equation}
\label{ura1}
=
\sum_{j_{q_1}=0}^{\infty}\ldots 
\sum_{j_{q_{k-1}}=0}^{\infty}
\sum_{j_s=0}^{p}
C_{j_k\ldots j_1}^2.
\end{equation}

\vspace{5mm}

Using (\ref{ura0}), (\ref{ura1}), and Lemma 1, we obtain

\vspace{2mm}
$$
\sum_{j_1=0}^{\infty}\ldots
\sum_{j_{s-1}=0}^{\infty}
\sum_{j_s=p+1}^{\infty}\sum_{j_{s+1}=0}^{\infty}\ldots \sum_{j_k=0}^{\infty}
C_{j_k\ldots j_1}^2=
\sum_{j_1=0}^{\infty}\ldots
\sum_{j_{s-1}=0}^{\infty}
\sum_{j_s=0}^{\infty}\sum_{j_{s+1}=0}^{\infty}\ldots \sum_{j_k=0}^{\infty}
C_{j_k\ldots j_1}^2-
$$

\vspace{3mm}
$$
-\sum_{j_1=0}^{\infty}\ldots
\sum_{j_{s-1}=0}^{\infty}
\sum_{j_s=0}^{p}\sum_{j_{s+1}=0}^{\infty}\ldots \sum_{j_k=0}^{\infty}
C_{j_k\ldots j_1}^2=
$$

\vspace{3mm}
$$
=
\sum_{j_s=0}^{\infty}
\sum_{j_{s-1}=0}^{\infty}\ldots
\sum_{j_1=0}^{\infty}\sum_{j_{s+1}=0}^{\infty}\ldots \sum_{j_k=0}^{\infty}
C_{j_k\ldots j_1}^2-
\sum_{j_s=0}^{p}
\sum_{j_{s-1}=0}^{\infty}\ldots
\sum_{j_1=0}^{\infty}\sum_{j_{s+1}=0}^{\infty}\ldots \sum_{j_k=0}^{\infty}
C_{j_k\ldots j_1}^2=
$$

\vspace{3mm}
$$
=
\sum_{j_s=p+1}^{\infty}
\sum_{j_{s-1}=0}^{\infty}\ldots
\sum_{j_1=0}^{\infty}\sum_{j_{s+1}=0}^{\infty}\ldots \sum_{j_k=0}^{\infty}
C_{j_k\ldots j_1}^2.
$$

\vspace{5mm}

The equality (\ref{d11}) is proved.

Applying the Parseval equality and Lemma 1, we obtain

\vspace{2mm}
$$
\int\limits_{[t,T]^k}K^2(t_1,\ldots,t_k)dt_1\ldots dt_k-
\sum_{j_1=0}^{p}\ldots \sum_{j_k=0}^{p}
C_{j_k\ldots j_1}^2=
$$

\vspace{3mm}

$$
=\sum_{j_1,\ldots,j_k=0}^{\infty}
C_{j_k\ldots j_1}^2-
\sum_{j_1=0}^{p}\ldots \sum_{j_k=0}^{p}
C_{j_k\ldots j_1}^2=
$$

\vspace{5mm}

$$
=
\sum_{j_1=0}^{\infty}\ldots \sum_{j_k=0}^{\infty}
C_{j_k\ldots j_1}^2-
\sum_{j_1=0}^{p}\ldots \sum_{j_k=0}^{p}
C_{j_k\ldots j_1}^2=
$$

\vspace{5mm}

$$
=\sum_{j_1=0}^{p}\sum_{j_2=0}^{\infty}\ldots \sum_{j_k=0}^{\infty}
C_{j_k\ldots j_1}^2+
\sum_{j_1=p+1}^{\infty}\sum_{j_2=0}^{\infty}\ldots \sum_{j_k=0}^{\infty}
C_{j_k\ldots j_1}^2-
\sum_{j_1=0}^{p}\ldots \sum_{j_k=0}^{p}
C_{j_k\ldots j_1}^2=
$$

\vspace{6mm}

$$
=\sum_{j_1=0}^{p}\sum_{j_2=0}^{p}\sum_{j_3=0}^{\infty}
\ldots \sum_{j_k=0}^{\infty}
C_{j_k\ldots j_1}^2+
\sum_{j_1=0}^{p}\sum_{j_2=p+1}^{\infty}
\sum_{j_3=0}^{\infty}
\ldots \sum_{j_k=0}^{\infty}+
$$

\vspace{4mm}

$$
+\sum_{j_1=p+1}^{\infty}\sum_{j_2=0}^{\infty}\ldots \sum_{j_k=0}^{\infty}
C_{j_k\ldots j_1}^2-
\sum_{j_1=0}^{p}\ldots \sum_{j_k=0}^{p}
C_{j_k\ldots j_1}^2=\ldots =
$$

\vspace{6mm}

$$
=\sum_{j_1=p+1}^{\infty}\sum_{j_2=0}^{\infty}\ldots \sum_{j_k=0}^{\infty}
C_{j_k\ldots j_1}^2+
\sum_{j_1=0}^p
\sum_{j_2=p+1}^{\infty}\sum_{j_2=0}^{\infty}\ldots \sum_{j_k=0}^{\infty}
C_{j_k\ldots j_1}^2+
$$

\vspace{6mm}

$$
+\sum_{j_1=0}^p\sum_{j_2=0}^p
\sum_{j_3=p+1}^{\infty}\sum_{j_4=0}^{\infty}\ldots \sum_{j_k=0}^{\infty}
C_{j_k\ldots j_1}^2+ \ldots +
\sum_{j_1=0}^p\ldots \sum_{j_{k-1}=0}^p
\sum_{j_k=p+1}^{\infty}C_{j_k\ldots j_1}^2\le
$$

\vspace{6mm}

$$
\le\sum_{j_1=p+1}^{\infty}\sum_{j_2=0}^{\infty}\ldots \sum_{j_k=0}^{\infty}
C_{j_k\ldots j_1}^2+
\sum_{j_1=0}^{\infty}
\sum_{j_2=p+1}^{\infty}\sum_{j_2=0}^{\infty}\ldots \sum_{j_k=0}^{\infty}
C_{j_k\ldots j_1}^2+
$$

\vspace{4mm}

$$
+\sum_{j_1=0}^{\infty}\sum_{j_2=0}^{\infty}
\sum_{j_3=p+1}^{\infty}\sum_{j_4=0}^{\infty}\ldots \sum_{j_k=0}^{\infty}
C_{j_k\ldots j_1}^2+ \ldots +
\sum_{j_1=0}^{\infty}\ldots \sum_{j_{k-1}=0}^{\infty}
\sum_{j_k=p+1}^{\infty}C_{j_k\ldots j_1}^2=
$$

\vspace{4mm}

\begin{equation}
\label{aaap}
=\sum\limits_{s=1}^k \left(\sum_{j_1=0}^{\infty}\ldots
\sum_{j_{s-1}=0}^{\infty}
\sum_{j_s=p+1}^{\infty}\sum_{j_{s+1}=0}^{\infty}\ldots \sum_{j_k=0}^{\infty}
C_{j_k\ldots j_1}^2\right).
\end{equation}

\vspace{6mm}

Note that deriving (\ref{aaap}), we used the following

\vspace{1mm}
$$
\sum_{j_1=0}^{p}\ldots
\sum_{j_{s-1}=0}^{p}
\sum_{j_s=p+1}^{\infty}\sum_{j_{s+1}=0}^{\infty}\ldots \sum_{j_k=0}^{\infty}
C_{j_k\ldots j_1}^2\le
$$

\vspace{3mm}
$$
\le
\sum_{j_1=0}^{m_1}\ldots
\sum_{j_{s-1}=0}^{m_{s-1}}
\sum_{j_s=p+1}^{\infty}\sum_{j_{s+1}=0}^{\infty}\ldots \sum_{j_k=0}^{\infty}
C_{j_k\ldots j_1}^2\le
$$

\vspace{3mm}
$$
\le
\lim\limits_{m_{s-1}\to\infty}
\sum_{j_1=0}^{m_1}\ldots
\sum_{j_{s-1}=0}^{m_{s-1}}
\sum_{j_s=p+1}^{\infty}\sum_{j_{s+1}=0}^{\infty}\ldots \sum_{j_k=0}^{\infty}
C_{j_k\ldots j_1}^2=
$$

\vspace{3mm}
$$
=
\sum_{j_1=0}^{m_1}\ldots
\sum_{j_{s-2}=0}^{m_{s-2}}\sum_{j_{s-1}=0}^{\infty}
\sum_{j_s=p+1}^{\infty}\sum_{j_{s+1}=0}^{\infty}\ldots \sum_{j_k=0}^{\infty}
C_{j_k\ldots j_1}^2\le
$$

\vspace{2mm}
$$
\le\ldots\le
$$

\vspace{1mm}
$$
\le\sum_{j_1=0}^{\infty}\ldots
\sum_{j_{s-1}=0}^{\infty}
\sum_{j_s=p+1}^{\infty}\sum_{j_{s+1}=0}^{\infty}\ldots \sum_{j_k=0}^{\infty}
C_{j_k\ldots j_1}^2,
$$

\vspace{5mm}
\noindent
where $m_1,\ldots,m_{s-1}>p.$

Denote

$$
C_{j_s\ldots j_1}(\tau)=
\int\limits_t^{\tau}
\phi_{j_s}(t_s)\psi_s(t_s)\ldots \int\limits_t^{t_2}
\phi_{j_1}(t_1)\psi_1(t_1)dt_1\ldots dt_s,
$$

\vspace{3mm}
\noindent
where
$s=1,\ldots,k-1.$

Let us remind the Dini Theorem, which we will use further.

\vspace{2mm}

{\bf Theorem (Dini).} {\it 
Let the functional sequence $u_n(x)$ 
be non-decreasing at each point of the interval $[a, b]$.
In addition, all the functions $u_n(x)$
of this sequence and the limit function $u(x)$ are continuous on the interval
$[a, b].$ Then the convergence $u_n(x)$ to 
$u(x)$ is uniform on the interval $[a,b].$}

\vspace{2mm}

For $s<k$ due to the Parseval equality and Dini Theorem as well as
(\ref{d11}) we obtain

\vspace{3mm}
$$
\sum_{j_1=0}^{\infty}\ldots
\sum_{j_{s-1}=0}^{\infty}
\sum_{j_s=p+1}^{\infty}\sum_{j_{s+1}=0}^{\infty}\ldots \sum_{j_k=0}^{\infty}
C_{j_k\ldots j_1}^2=
$$

\vspace{6mm}

$$
\stackrel{\hbox{(\ref{d11})}}{=}\ \ \ 
\sum_{j_s=p+1}^{\infty}
\sum_{j_{s-1}=0}^{\infty}\ldots
\sum_{j_{1}=0}^{\infty}
\sum_{j_{s+1}=0}^{\infty}\ldots \sum_{j_k=0}^{\infty}
C_{j_k\ldots j_1}^2=
$$

\vspace{6mm}

$$
\stackrel{\hbox{(Parseval Eq.)}}{=}\ \ \
\sum_{j_s=p+1}^{\infty}
\sum_{j_{s-1}=0}^{\infty}\ldots
\sum_{j_{1}=0}^{\infty}
\sum_{j_{s+1}=0}^{\infty}\ldots 
\sum_{j_{k-1}=0}^{\infty}
\int\limits_t^T \psi_k^2(t_k) \left(C_{j_{k-1}\ldots j_1}(t_k)\right)^2 dt_k=
$$

\vspace{6mm}

$$
\stackrel{\hbox{(Dini Th.)}}{=}\ \ \
\sum_{j_s=p+1}^{\infty}
\sum_{j_{s-1}=0}^{\infty}\ldots
\sum_{j_{1}=0}^{\infty}
\sum_{j_{s+1}=0}^{\infty}\ldots 
\sum_{j_{k-2}=0}^{\infty}
\int\limits_t^T \psi_k^2(t_k) 
\sum_{j_{k-1}=0}^{\infty}\left(C_{j_{k-1}\ldots j_1}(t_k)\right)^2 dt_k=
$$

\vspace{6mm}

$$
\stackrel{\hbox{(Parseval Eq.)}}{=}\ \ \
\sum_{j_s=p+1}^{\infty}
\sum_{j_{s-1}=0}^{\infty}\ldots
\sum_{j_{1}=0}^{\infty}
\sum_{j_{s+1}=0}^{\infty}\ldots 
\sum_{j_{k-2}=0}^{\infty}
\int\limits_t^T \psi_k^2(t_k) \int\limits_t^{t_k} \psi_{k-1}^2(t_{k-1}) 
\left(C_{j_{k-2}\ldots j_1}(t_{k-1})\right)^2\times
$$

\vspace{2mm}
$$
\times dt_{k-1}dt_k\le
$$

\vspace{4mm}
$$
\le C\sum_{j_s=p+1}^{\infty}
\sum_{j_{s-1}=0}^{\infty}\ldots
\sum_{j_{1}=0}^{\infty}
\sum_{j_{s+1}=0}^{\infty}\ldots 
\sum_{j_{k-2}=0}^{\infty}
\int\limits_t^T 
\left(C_{j_{k-2}\ldots j_1}(\tau)\right)^2 d\tau =
$$

\vspace{5mm}

$$
\stackrel{\hbox{(Dini Th.)}}{=}\ \ \
C\sum_{j_s=p+1}^{\infty}
\sum_{j_{s-1}=0}^{\infty}\ldots
\sum_{j_{1}=0}^{\infty}
\sum_{j_{s+1}=0}^{\infty}\ldots 
\sum_{j_{k-3}=0}^{\infty}
\int\limits_t^T 
\sum_{j_{k-2}=0}^{\infty}
\left(C_{j_{k-2}\ldots j_1}(\tau)\right)^2 d\tau =
$$

\vspace{5mm}

$$
\stackrel{\hbox{(Parseval Eq.)}}{=}\ \ \ C\sum_{j_s=p+1}^{\infty}
\sum_{j_{s-1}=0}^{\infty}\ldots
\sum_{j_{1}=0}^{\infty}
\sum_{j_{s+1}=0}^{\infty}\ldots 
\sum_{j_{k-3}=0}^{\infty}
\int\limits_t^T \int\limits_t^{\tau}
\psi_{k-2}^2(\theta)
\left(C_{j_{k-3}\ldots j_1}(\theta)\right)^2
d\theta d\tau\le 
$$

\vspace{5mm}

$$
\le K
\sum_{j_s=p+1}^{\infty}
\sum_{j_{s-1}=0}^{\infty}\ldots
\sum_{j_{1}=0}^{\infty}
\sum_{j_{s+1}=0}^{\infty}\ldots 
\sum_{j_{k-3}=0}^{\infty}
\int\limits_t^T
\left(C_{j_{k-3}\ldots j_1}(\tau)\right)^2
d\tau\le 
$$

\vspace{4mm}
$$
\le \ldots \le
$$

$$
\le C_k
\sum_{j_s=p+1}^{\infty}
\sum_{j_{s-1}=0}^{\infty}\ldots
\sum_{j_{1}=0}^{\infty}
\int\limits_t^T 
\left(C_{j_{s}\ldots j_1}(\tau)\right)^2 d\tau=
$$

\vspace{5mm}

\begin{equation}
\label{d14}
\stackrel{\hbox{(Dini Th.)}}{=}\ \ \ C_k
\sum_{j_s=p+1}^{\infty}
\sum_{j_{s-1}=0}^{\infty}\ldots
\sum_{j_{2}=0}^{\infty}
\int\limits_t^T  \sum_{j_{1}=0}^{\infty}
\left(C_{j_{s}\ldots j_1}(\tau)\right)^2 d\tau,
\end{equation}

\vspace{6mm}
\noindent
where constants $C,$ $K$ depend on $T-t$ and
constant $C_k$ depends on $k$ and $T-t.$

Let us explane more precisely how we obtain (\ref{d14}).
For any function $g(s)\in L_2([t,T])$ we have the following
Parseval equality

$$
\sum\limits_{j=0}^{\infty}\left(\int\limits_t^{\tau}
\phi_j(s)g(s)ds\right)^2=
\sum\limits_{j=0}^{\infty}\left(\int\limits_t^T
{\bf 1}_{\{s<\tau\}}\phi_j(s)g(s)ds\right)^2=
$$

\begin{equation}
\label{d15}
=\int\limits_t^T
\left({\bf 1}_{\{s<\tau\}}\right)^2 g^2(s)ds=
\int\limits_t^{\tau}
g^2(s)ds.
\end{equation}

\vspace{3mm}

The equality (\ref{d15}) has been applied repeatedly when we obtaining
(\ref{d14}).
Using the integration order replacement in Riemann integrals, we have

\vspace{1mm}
$$
C_{j_s\ldots j_1}(\tau)=
\int\limits_t^{\tau}
\phi_{j_s}(t_s)\psi_s(t_s)\ldots \int\limits_t^{t_2}
\phi_{j_1}(t_1)\psi_1(t_1)dt_1\ldots dt_s=
$$

\vspace{2mm}
$$
=\int\limits_t^{\tau}
\phi_{j_1}(t_1)\psi_1(t_1)\int\limits_{t_1}^{\tau}
\phi_{j_2}(t_2)\psi_2(t_2)
\ldots
\int\limits_{t_{s-1}}^{\tau}
\phi_{j_s}(t_s)\psi_s(t_s)dt_s\ldots dt_2dt_1
\stackrel{\sf def}{=}
$$

\vspace{4mm}
$$
\stackrel{\sf def}{=}
{\tilde C}_{j_s\ldots j_1}(\tau).
$$

\vspace{7mm}

For $l=1,\ldots,s$ we will use the following notation

\vspace{2mm}
$$
{\tilde C}_{j_s\ldots j_l}(\tau,\theta)=
\int\limits_{\theta}^{\tau}
\phi_{j_l}(t_l)\psi_l(t_l)\int\limits_{t_l}^{\tau}
\phi_{j_{l+1}}(t_{l+1})\psi_{l+1}(t_{l+1})
\ldots
\int\limits_{t_{s-1}}^{\tau}
\phi_{j_s}(t_s)\psi_s(t_s)dt_s\ldots dt_{l+1}dt_l.
$$

\vspace{5mm}

Applying the Parseval equality and Dini Theorem, from (\ref{d14}) we obtain

\vspace{2mm}
$$
\sum_{j_1=0}^{\infty}\ldots
\sum_{j_{s-1}=0}^{\infty}
\sum_{j_s=p+1}^{\infty}\sum_{j_{s+1}=0}^{\infty}\ldots \sum_{j_k=0}^{\infty}
C_{j_k\ldots j_1}^2\le
$$

\vspace{5mm}

$$
\le
C_k
\sum_{j_s=p+1}^{\infty}
\sum_{j_{s-1}=0}^{\infty}\ldots
\sum_{j_{2}=0}^{\infty}
\int\limits_t^T  \sum_{j_{1}=0}^{\infty}
\left(C_{j_{s}\ldots j_1}(\tau)\right)^2 d\tau=
$$

\vspace{5mm}

$$
=C_k
\sum_{j_s=p+1}^{\infty}
\sum_{j_{s-1}=0}^{\infty}\ldots
\sum_{j_{2}=0}^{\infty}
\int\limits_t^T  \sum_{j_{1}=0}^{\infty}
\left({\tilde C}_{j_{s}\ldots j_1}(\tau)\right)^2 d\tau=
$$

\vspace{5mm}

\begin{equation}
\label{molod1}
\stackrel{\hbox{(Parseval Eq.)}}{=}\ \ \ C_k
\sum_{j_s=p+1}^{\infty}
\sum_{j_{s-1}=0}^{\infty}\ldots
\sum_{j_{2}=0}^{\infty}
\int\limits_t^T\int\limits_t^{\tau}\psi_1^2(t_1)  
\left({\tilde C}_{j_{s}\ldots j_2}(\tau,t_1)\right)^2 dt_1d\tau=
\end{equation}

\vspace{5mm}

\begin{equation}
\label{molod2}
\stackrel{\hbox{(Dini Th.)}}{=}\ \ \ C_k
\sum_{j_s=p+1}^{\infty}
\sum_{j_{s-1}=0}^{\infty}\ldots
\sum_{j_{3}=0}^{\infty}
\int\limits_t^T\int\limits_t^{\tau}\psi_1^2(t_1)  
\sum_{j_{2}=0}^{\infty}
\left({\tilde C}_{j_{s}\ldots j_2}(\tau,t_1)\right)^2 dt_1d\tau=
\end{equation}

\vspace{5mm}

$$
\stackrel{\hbox{(Parseval Eq.)}}{=}\ \ \ C_k
\sum_{j_s=p+1}^{\infty}
\sum_{j_{s-1}=0}^{\infty}\ldots
\sum_{j_{3}=0}^{\infty}
\int\limits_t^T\int\limits_t^{\tau}\psi_1^2(t_1)  
\int\limits_{t_1}^{\tau}\psi_2^2(t_2)  
\left({\tilde C}_{j_{s}\ldots j_3}(\tau,t_2)\right)^2 dt_2dt_1d\tau \le
$$

\vspace{5mm}

$$
\le C_k
\sum_{j_s=p+1}^{\infty}
\sum_{j_{s-1}=0}^{\infty}\ldots
\sum_{j_{3}=0}^{\infty}
\int\limits_t^T\int\limits_t^{\tau}\psi_1^2(t_1)  
\int\limits_{t}^{\tau}\psi_2^2(t_2)  
\left({\tilde C}_{j_{s}\ldots j_3}(\tau,t_2)\right)^2 dt_2dt_1d\tau\le
$$

\vspace{5mm}

$$
\le C^{'}_k
\sum_{j_s=p+1}^{\infty}
\sum_{j_{s-1}=0}^{\infty}\ldots
\sum_{j_{3}=0}^{\infty}
\int\limits_t^T
\int\limits_{t}^{\tau}\psi_2^2(t_2)  
\left({\tilde C}_{j_{s}\ldots j_3}(\tau,t_2)\right)^2 dt_2d\tau
\le 
$$

\vspace{4mm}

$$
\le  \ldots \le
$$

\vspace{1mm}

$$
\le C^{''}_k
\sum_{j_s=p+1}^{\infty}
\int\limits_t^T\int\limits_t^{\tau}
\psi_{s-1}^2(t_{s-1})
\left({\tilde C}_{j_{s}}(\tau,t_{s-1})\right)^2 dt_{s-1} d\tau\le
$$

\vspace{4mm}

\begin{equation}
\label{la}
\le {\tilde C}_k
\sum_{j_s=p+1}^{\infty}
\int\limits_t^T\int\limits_t^{\tau}
\left(~\int\limits_{u}^{\tau}\phi_{j_s}(\theta)
\psi_s(\theta)d\theta\right)^2 du d\tau,
\end{equation}

\vspace{6mm}
\noindent
where constants $C^{'}_k,$ $C^{''}_k,$ $\tilde C_k$
depend on $k$ and $T-t.$

Let us explane more precisely how we obtain (\ref{la}).
For any function $g(s)\in L_2([t,T])$ we have the following
Parseval equality

$$
\sum\limits_{j=0}^{\infty}\left(\int\limits_{\theta}^{\tau}
\phi_j(s)g(s)ds\right)^2=
\sum\limits_{j=0}^{\infty}\left(\int\limits_t^T
{\bf 1}_{\{\theta<s<\tau\}}\phi_j(s)g(s)ds\right)^2=
$$

\begin{equation}
\label{d22}
=\int\limits_t^T
\left({\bf 1}_{\{\theta<s<\tau\}}\right)^2 g^2(s)ds=
\int\limits_{\theta}^{\tau}
g^2(s)ds.
\end{equation}

\vspace{3mm}

The equality (\ref{d22}) has been applied repeatedly when we obtaining
(\ref{la}).
Let us explane more precisely the passing from (\ref{molod1})
to (\ref{molod2}) (the same steps were used when we were
deriving (\ref{la})).

We have

$$
\int\limits_t^T\int\limits_t^{\tau}\psi_1^2(t_1)  
\sum_{j_{2}=0}^{\infty}
\left({\tilde C}_{j_{s}\ldots j_2}(\tau,t_1)\right)^2 dt_1d\tau -
\sum_{j_{2}=0}^{n}\int\limits_t^T\int\limits_t^{\tau}\psi_1^2(t_1)  
\left({\tilde C}_{j_{s}\ldots j_2}(\tau,t_1)\right)^2 dt_1d\tau =
$$

\vspace{4mm}

$$
=\int\limits_t^T\int\limits_t^{\tau}\psi_1^2(t_1)  
\sum_{j_{2}=n+1}^{\infty}
\left({\tilde C}_{j_{s}\ldots j_2}(\tau,t_1)\right)^2 dt_1d\tau =
$$

\vspace{4mm}

\begin{equation}
\label{molod3}
=\lim\limits_{N\to\infty}
\sum\limits_{j=0}^{N-1}\int\limits_t^{\tau_j}\psi_1^2(t_1)  
\sum_{j_{2}=n+1}^{\infty}
\left({\tilde C}_{j_{s}\ldots j_2}(\tau_j,t_1)\right)^2 dt_1 \Delta\tau_j,
\end{equation}

\vspace{4mm}
\noindent
where $\{\tau_j\}_{j=0}^{N}$ is the partition of the 
interval $[t, T],$ which satisfies the condition (\ref{1111}).

Since the non-decreasing functional sequence $u_n(\tau_j,t_1)$ and its
limit function $u(\tau_j,t_1)$ are continuous on the
interval $[t,\tau_j]\subseteq [t, T]$ with respect to $t_1$,
where

$$
u_n(\tau_j,t_1)=
\sum_{j_{2}=0}^{n}
\left({\tilde C}_{j_{s}\ldots j_2}(\tau_j,t_1)\right)^2,
$$

$$
u(\tau_j,t_1)=
\sum_{j_{2}=0}^{\infty}
\left({\tilde C}_{j_{s}\ldots j_2}(\tau_j,t_1)\right)^2=
\int\limits_{t_1}^{\tau_j}
\psi_2^2(t_2)
\left({\tilde C}_{j_{s}\ldots j_3}(\tau_j,t_2)\right)^2 dt_2,
$$

\vspace{4mm}

\noindent 
then by Dini Theorem we have the uniform convergence
of $u_n(\tau_j,t_1)$ to $u(\tau_j,t_1)$ at the interval $[t,\tau_j]\subseteq
[t, T]$
with respect to $t_1.$ As a result, we obtain

\begin{equation}
\label{molod4}
\sum_{j_{2}=n+1}^{\infty}
\left({\tilde C}_{j_{s}\ldots j_2}(\tau_j,t_1)\right)^2<\varepsilon,\ \ \ 
t_1\in [t,\tau_j]
\end{equation}

\vspace{4mm}
\noindent
for $n>N(\varepsilon)$ ($N(\varepsilon)$ exists
for any $\varepsilon>0$ and it does not depend on $t_1$).

From (\ref{molod3}) and (\ref{molod4}) we get

$$
\lim\limits_{N\to\infty}
\sum\limits_{j=0}^{N-1}\int\limits_t^{\tau_j}\psi_1^2(t_1)  
\sum_{j_{2}=n+1}^{\infty}
\left({\tilde C}_{j_{s}\ldots j_2}(\tau_j,t_1)\right)^2 dt_1 \Delta\tau_j
\le
\varepsilon 
\lim\limits_{N\to\infty}
\sum\limits_{j=0}^{N-1}\int\limits_t^{\tau_j}\psi_1^2(t_1)  
dt_1 \Delta\tau_j= 
$$

\vspace{3mm}
\begin{equation}
\label{molod6}
=\varepsilon \int\limits_t^T
\int\limits_t^{\tau}\psi_1^2(t_1)  
dt_1 d\tau.
\end{equation}

\vspace{4mm}

From (\ref{molod6}) we have

$$
\lim\limits_{n\to\infty}\int\limits_t^T\int\limits_t^{\tau}\psi_1^2(t_1)  
\sum_{j_{2}=n+1}^{\infty}
\left({\tilde C}_{j_{s}\ldots j_2}(\tau,t_1)\right)^2 dt_1d\tau = 0.
$$

\vspace{4mm}

This fact completes the proof of passing 
from (\ref{molod1})
to (\ref{molod2}).

Let us estimate the integral 

\vspace{-2mm}
\begin{equation}
\label{st1}
\int\limits_{u}^{\tau}\phi_{j_s}(\theta)
\psi_s(\theta)d\theta
\end{equation}

\vspace{2mm}
\noindent
from (\ref{la}) for the case when $\{\phi_j(s)\}_{j=0}^{\infty}$
is a complete orthonormal system of Legendre polynomials or
trigonometric functions in the space $L_2([t,T])$.

Note that the estimates for the integral

\vspace{-2mm}
\begin{equation}
\label{st2}
\int\limits_{t}^{\tau}\phi_{j}(\theta)\psi(\theta)d\theta,\ \ \ j\ge p+1
\end{equation}

\vspace{2mm}
\noindent
have been obtained in \cite{art-5}, \cite{arxiv-5},
where $\psi(\theta)$ is a continuously
differentiable function on the interval $[t, T]$.
The same estimates 
can also be found in early publications \cite{9}-\cite{11}, \cite{12} or
in \cite{11a}-\cite{11aaa-a} (2020-2023).

Let us estimate the integral (\ref{st1}) using the approach from
\cite{art-5}, \cite{arxiv-5}.
First consider the case of Legendre polynomials.
Then $\phi_j(s)$ looks as follows

\vspace{1mm}
\begin{equation}
\label{ogo7}
\phi_j(\theta)=\sqrt{\frac{2j+1}{T-t}}P_j\left(\left(
\theta-\frac{T+t}{2}\right)\frac{2}{T-t}\right),\ \ \ j\ge 0,
\end{equation}

\vspace{5mm}
\noindent
where $P_j(x)$ is the Legendre polynomial.

Further, we have 

$$
\int\limits_v^x\phi_{j}(\theta)\psi(\theta)d\theta=
\frac{\sqrt{T-t}\sqrt{2j+1}}{2}
\int\limits_{z(v)}^{z(x)}P_{j}(y)
\psi(u(y))dy=
$$

\vspace{2mm}
$$
=\frac{\sqrt{T-t}}{2\sqrt{2j+1}}\Biggl((P_{j+1}(z(x))-
P_{j-1}(z(x)))\psi(x)-
(P_{j+1}(z(v))-
P_{j-1}(z(v)))\psi(v)-
\Biggr.
$$

\vspace{2mm}
\begin{equation}
\label{6000}
\Biggl.-
\frac{T-t}{2}
\int\limits_{z(v)}^{z(x)}((P_{j+1}(y)-P_{j-1}(y))
{\psi}'(u(y))dy\Biggr),
\end{equation}

\vspace{5mm}
\noindent
where $x, v\in (t, T),$ $j\ge p+1,$ 
$u(y)$ and $z(x)$ are defined by the following relations

\vspace{1mm}
$$
u(y)=\frac{T-t}{2}y+\frac{T+t}{2},\ \ \
z(x)=\left(x-\frac{T+t}{2}\right)\frac{2}{T-t},
$$

\vspace{4mm}
\noindent
${\psi}'$ is a derivative of the function $\psi(\theta)$
with respect to the variable $u(y).$

Note that in (\ref{6000}) we used the following well-known property
of the Legendre polynomials

\vspace{2mm}
$$
\frac{dP_{j+1}}{dx}(x)-\frac{dP_{j-1}}{dx}(x)=(2j+1)P_j(x),\ \ \ 
j=1, 2,\ldots
$$

\vspace{5mm}

From (\ref{6000}) and the well-known estimate for the Legendre
polynomials

\vspace{1mm}
\begin{equation}
\label{200}
|P_j(y)| <\frac{K}{\sqrt{j+1}(1-y^2)^{1/4}},\ \ \ 
y\in (-1, 1),\ \ \ j\in \mathbb{N},
\end{equation}

\vspace{5mm}
\noindent
where constant $K$ does not depend on $y$ and $j$, it follows that

\vspace{2mm}
\begin{equation}
\label{101}
\left|
\int\limits_v^x\phi_{j}(\theta)\psi(\theta)d\theta
\right| <
\frac{C}{j}\Biggl(\frac{1}{(1-(z(x))^2)^{1/4}}+
\frac{1}{(1-(z(v))^2)^{1/4}}+C_1\Biggr),
\end{equation}

\vspace{5mm}
\noindent
where constants $C, C_1$ do not depend on $j$ $(j>0)$ and
$z(x), z(v)\in (-1, 1),$ $x, v\in (t, T).$

From (\ref{101}) we obtain

\vspace{1mm}
\begin{equation}
\label{102}
\left(
\int\limits_v^x\phi_{j}(\theta)\psi(\theta)d\theta
\right)^2 <
\frac{C_2}{j^2}\Biggl(\frac{1}{(1-(z(x))^2)^{1/2}}+
\frac{1}{(1-(z(v))^2)^{1/2}}+C_3\Biggr),
\end{equation}

\vspace{4mm}
\noindent
where constants $C_2, C_3$ do not depend on $j$ $(j>0)$.

Let us apply (\ref{102}) for estimation of the right-hand side
of (\ref{la}). We have

\vspace{1mm}
$$
\int\limits_t^T\int\limits_t^{\tau}
\left(~\int\limits_{u}^{\tau}\phi_{j_s}(\theta)
\psi_s(\theta)d\theta\right)^2 du d\tau\le
$$

$$
\le \frac{K_1}{j_s^2}
\left(
\int\limits_{-1}^1
\frac{dy}{\left(1-y^2\right)^{1/2}}+
\int\limits_{-1}^1\int\limits_{-1}^x
\frac{dy}{\left(1-y^2\right)^{1/2}}dx + K_2\right)\le
$$

\vspace{2mm}
\begin{equation}
\label{103}
\le \frac{K_3}{j_s^2},
\end{equation}

\vspace{5mm}
\noindent
where constants $K_1, K_2, K_3$ are independent of $j_s$ $(j_s>0).$

Now consider the trigonometric case.
The complete orthonormal system of trigonometric functions
in the space $L_2([t, T])$ has the following form

\begin{equation}
\label{trig11}
\phi_j(\theta)=\frac{1}{\sqrt{T-t}}
\left\{
\begin{matrix}
1,\ & j=0\cr\cr\cr
\sqrt{2}{\rm sin} \left(2\pi r(\theta-t)/(T-t)\right),\ & j=2r-1\cr\cr\cr
\sqrt{2}{\rm cos} \left(2\pi r(\theta-t)/(T-t)\right),\ & j=2r
\end{matrix}
,\right.
\end{equation}

\vspace{3mm}
\noindent
where $r=1, 2,\ldots $

Using the system of functions 
(\ref{trig11}), we have

\vspace{1mm}
$$
\int\limits_v^x\phi_{2r-1}(\theta)\psi(\theta)d\theta=
\sqrt{\frac{2}{T-t}}\int\limits_v^x
{\rm sin} \frac{2\pi r(\theta-t)}{T-t}\psi(\theta)d\theta=
$$

\vspace{2mm}
$$
=-\sqrt{\frac{T-t}{2}}\frac{1}{\pi r}\Biggl(
\psi(x){\rm cos}\frac{2\pi r(x-t)}{T-t}-
\psi(v){\rm cos}\frac{2\pi r(v-t)}{T-t}-\Biggr.
$$

\vspace{2mm}
\begin{equation}
\label{201}
\Biggl.-
\int\limits_v^x
{\rm cos} \frac{2\pi r(\theta-t)}{T-t}\psi'(\theta)d\theta\Biggr),
\end{equation}

\vspace{5mm}

$$
\int\limits_v^x\phi_{2r}(\theta)\psi(\theta)d\theta=
\sqrt{\frac{2}{T-t}}\int\limits_v^x
{\rm cos} \frac{2\pi r(\theta-t)}{T-t}\psi(\theta)d\theta=
$$

\vspace{2mm}
$$
=\sqrt{\frac{T-t}{2}}\frac{1}{\pi r}\Biggl(
\psi(x){\rm sin}\frac{2\pi r(x-t)}{T-t}-
\psi(v){\rm sin}\frac{2\pi r(v-t)}{T-t}-\Biggr.
$$

\vspace{2mm}
\begin{equation}
\label{202}
\Biggl.-
\int\limits_v^x
{\rm sin} \frac{2\pi r(\theta-t)}{T-t}\psi'(\theta)d\theta\Biggr),
\end{equation}

\vspace{5mm}
\noindent
where $\psi'(\theta)$ is a derivative of the function $\psi(\theta)$
with respect to the variable $\theta.$

Combining (\ref{201}) and (\ref{202}), we obtain for the
trigonometric case

\begin{equation}
\label{203}
\left(
\int\limits_v^x\phi_{j}(\theta)\psi(\theta)d\theta
\right)^2 \le 
\frac{C_4}{j^2},
\end{equation}

\vspace{4mm}
\noindent
where constant $C_4$ is independent of $j$ $(j>0).$  

From (\ref{203}) we finally have

\begin{equation}
\label{103x}
\int\limits_t^T\int\limits_t^{\tau}
\left(~\int\limits_{u}^{\tau}\phi_{j_s}(\theta)
\psi_s(\theta)d\theta\right)^2 du d\tau
\le \frac{K_4}{j_s^2},
\end{equation}

\vspace{4mm}
\noindent
where constant $K_4$ does not depend on $j_s$ $(j_s>0).$

Combibing (\ref{la}), (\ref{103}), and (\ref{103x}), we obtain

\vspace{2mm}
$$
\sum_{j_1=0}^{\infty}\ldots
\sum_{j_{s-1}=0}^{\infty}
\sum_{j_s=p+1}^{\infty}\sum_{j_{s+1}=0}^{\infty}\ldots \sum_{j_k=0}^{\infty}
C_{j_k\ldots j_1}^2\le
$$

\vspace{2mm}
\begin{equation}
\label{fff}
\le L_k
\sum_{j_s=p+1}^{\infty}\frac{1}{j_s^2} \le 
\frac{L_k}{p},
\end{equation}

\vspace{5mm}
\noindent
where constant $L_k$ depends on $k$ and $T-t.$

Obviously, the case $s=k$ can be considered absolutely analogously to the
case $s<k$. Then from (\ref{aaap}) and (\ref{fff})
we obtain

\begin{equation}
\label{ddd1}
\int\limits_{[t,T]^k}K^2(t_1,\ldots,t_k)dt_1\ldots dt_k-
\sum_{j_1=0}^{p}\ldots \sum_{j_k=0}^{p}
C_{j_k\ldots j_1}^2\le \frac{G_k}{p},
\end{equation}

\vspace{3mm}
\noindent
where constant $G_k$ depends on $k$ and $T-t.$

For the further consideration we consider the following theorem.

\vspace{2mm} 
                    
{\bf Theorem 7}\ \cite{11a} (Sect.~1.1.9, 1.11, 1.12), \cite{arxiv-1} (Sect.~6, 15, 16).\ 
{\it Under the conditions of Theorems {\rm 1, 3} the following
estimate is correct

$$
{\sf M}\left\{\biggl(J[\psi^{(k)}]_{T,t}-
J[\psi^{(k)}]_{T,t}^{p_1,\ldots,p_k}\biggr)^{2n}\right\}\le
$$

\vspace{3mm}
$$
\le
(k!)^{n}(2n-1)^{nk}\ \times
$$

\vspace{1mm}
\begin{equation}
\label{dima2ye100}
\times\ 
\left(
\int\limits_{[t,T]^k}
K^2(t_1,\ldots,t_k)
dt_1\ldots dt_k -\sum_{j_1=0}^{p_1}\ldots
\sum_{j_k=0}^{p_k}C^2_{j_k\ldots j_1}
\right)^n
\end{equation}

\vspace{5mm}
\noindent
for $n\in \mathbb{N},$ where $J[\psi^{(k)}]_{T,t}^{p_1,\ldots,p_k}$
is the expression on the right-hand side of {\rm (\ref{leto6000})} 
before passing to the limit, i.e.

$$
J[\psi^{(k)}]_{T,t}^{p_1,\ldots,p_k}=
\sum\limits_{j_1=0}^{p_1}\ldots
\sum\limits_{j_k=0}^{p_k}
C_{j_k\ldots j_1}\Biggl(
\prod_{l=1}^k\zeta_{j_l}^{(i_l)}+\sum\limits_{r=1}^{[k/2]}
(-1)^r \times
\Biggr.
$$

\vspace{4mm}
$$
\times
\sum_{\stackrel{(\{\{g_1, g_2\}, \ldots, 
\{g_{2r-1}, g_{2r}\}\}, \{q_1, \ldots, q_{k-2r}\})}
{{}_{\{g_1, g_2, \ldots, 
g_{2r-1}, g_{2r}, q_1, \ldots, q_{k-2r}\}=\{1, 2, \ldots, k\}}}}
\prod\limits_{s=1}^r
{\bf 1}_{\{i_{g_{{}_{2s-1}}}=~i_{g_{{}_{2s}}}\ne 0\}}
\Biggl.{\bf 1}_{\{j_{g_{{}_{2s-1}}}=~j_{g_{{}_{2s}}}\}}
\prod_{l=1}^{k-2r}\zeta_{j_{q_l}}^{(i_{q_l})}\Biggr),
$$

\vspace{7mm}
\noindent
where $i_1,\ldots,i_k=1,\ldots,m$ and
the remainder notations are the same as in Theorems {\rm 1, 3}}.

\vspace{2mm}

Using (\ref{ddd1}) and Theorem 7 for the case $p_1=\ldots=p_k=p$ and $n=2$ 
(see (\ref{dima2ye100})),
we obtain

\vspace{1mm}
$$
{\sf M}\left\{\biggl(J[\psi^{(k)}]_{T,t}-
J[\psi^{(k)}]_{T,t}^{p,\ldots,p}\biggr)^{4}\right\}\le
$$

\vspace{1mm}
\begin{equation}
\label{fff5}
\le C_{2,k}
\left(
\int\limits_{[t,T]^k}
K^2(t_1,\ldots,t_k)
dt_1\ldots dt_k -\sum_{j_1=0}^{p}\ldots
\sum_{j_k=0}^{p}C^2_{j_k\ldots j_1}
\right)^2
\le 
\frac{H_{2,k}}{p^2},
\end{equation}

\vspace{5mm}
\noindent
where $H_{2,k}=G_k^2{C}_{2,k}$ and

\vspace{-2mm}
$$
C_{n,k}=(k!)^{n}(2n-1)^{nk}.
$$

\vspace{5mm}

Let us consider the the well-known fact.

\vspace{2mm}

{\bf Proposition 2.}\ {\it If for the sequence of random variables
$\xi_p$ and for some real
$\alpha>0$ the number series 

\vspace{-3mm}
$$
\sum\limits_{p=1}^{\infty}{\sf M}\left\{\left|\xi_p\right|^{\alpha}\right\}
$$

\vspace{5mm}
\noindent
converges, then the sequence $\xi_p$ converges to zero w. p. {\rm 1}.}

\vspace{2mm}

Let us put
$$
\xi_p=\biggl|J[\psi^{(k)}]_{T,t}-
J[\psi^{(k)}]_{T,t}^{p,\ldots,p}\biggr|
$$

\vspace{3mm}
\noindent
and $\alpha=4.$

Then from (\ref{fff5}) we obtain

\vspace{1mm}
\begin{equation}
\label{qqq1}
\sum\limits_{p=1}^{\infty}
{\sf M}\left\{\biggl(J[\psi^{(k)}]_{T,t}-
J[\psi^{(k)}]_{T,t}^{p,\ldots,p}\biggr)^{4}\right\}
\le H_{2,k}\sum\limits_{p=1}^{\infty}\frac{1}{p^2}<\infty.
\end{equation}

\vspace{4mm}

From (\ref{qqq1}) we get

\vspace{2mm}
$$
J[\psi^{(k)}]_{T,t}^{p,\ldots,p}\ \to \ J[\psi^{(k)}]_{T,t}\ \ \ 
\hbox{if}\ \ \ p\to \infty
$$

\vspace{5mm}
\noindent
w.~p.~1, where (see Theorem 1)

\vspace{2mm}
$$
J[\psi^{(k)}]_{T,t}^{p,\ldots,p}=
\sum_{j_1=0}^{p}\ldots\sum_{j_k=0}^{p}
C_{j_k\ldots j_1}\Biggl(
\prod_{l=1}^k\zeta_{j_l}^{(i_l)}\ -
\Biggr.
$$

\vspace{3mm}
\begin{equation}
\label{kk0}
-\ \Biggl.
\hbox{\vtop{\offinterlineskip\halign{
\hfil#\hfil\cr
{\rm l.i.m.}\cr
$\stackrel{}{{}_{N\to \infty}}$\cr
}} }\sum_{(l_1,\ldots,l_k)\in {\rm G}_k}
\phi_{j_{1}}(\tau_{l_1})
\Delta{\bf w}_{\tau_{l_1}}^{(i_1)}\ldots
\phi_{j_{k}}(\tau_{l_k})
\Delta{\bf w}_{\tau_{l_k}}^{(i_k)}\Biggr)
\end{equation}

\vspace{7mm}
\noindent
or (see (\ref{leto6000}))

\vspace{2mm}
$$
J[\psi^{(k)}]_{T,t}^{p,\ldots,p}=
\sum\limits_{j_1=0}^{p}\ldots
\sum\limits_{j_k=0}^{p}
C_{j_k\ldots j_1}\Biggl(
\prod_{l=1}^k\zeta_{j_l}^{(i_l)}+\sum\limits_{r=1}^{[k/2]}
(-1)^r \times
\Biggr.
$$

\vspace{3mm}
\begin{equation}
\label{kk1}
\times
\sum_{\stackrel{(\{\{g_1, g_2\}, \ldots, 
\{g_{2r-1}, g_{2r}\}\}, \{q_1, \ldots, q_{k-2r}\})}
{{}_{\{g_1, g_2, \ldots, 
g_{2r-1}, g_{2r}, q_1, \ldots, q_{k-2r}\}=\{1, 2, \ldots, k\}}}}
\prod\limits_{s=1}^r
{\bf 1}_{\{i_{g_{{}_{2s-1}}}=~i_{g_{{}_{2s}}}\ne 0\}}
\Biggl.{\bf 1}_{\{j_{g_{{}_{2s-1}}}=~j_{g_{{}_{2s}}}\}}
\prod_{l=1}^{k-2r}\zeta_{j_{q_l}}^{(i_{q_l})}\Biggr),
\end{equation}

\vspace{9mm}
\noindent
where $i_1,\ldots,i_k=1,\ldots,m$ in (\ref{kk0}) and (\ref{kk1}).
Theorem 6 is proved.

\vspace{5mm}

\section{Mean-Square Approximation of Iterated Stratonovich 
Stochastic Integrals of Multiplicities 1 to 6}

\vspace{5mm}

This section is devoted to the mean-square approximation
of iterated Stratonovich stochastic integrals.
We consider the adaptation of Theorems 1, 3 for 
iterated Stratonovich stochastic integrals
of multiplicities 1 to 6.
Also we consider the question on the exact calculation 
of the mean-square approximation errors
for the following iterated Stratonovich stochastic integrals

\vspace{1mm}
$$
I_{(0)T,t}^{*(i_1)},\ \ \ 
I_{(1)T,t}^{*(i_1)},\ \ \ 
I_{(00)T,t}^{*(i_1i_2)},\ \ \ 
I_{(000)T,t}^{*(i_1i_2i_3)},\ \ \ 
i_1, i_2, i_3=1,\ldots,m,
$$

\vspace{4mm}
\noindent
where 
\begin{equation}
\label{k1001xxxx}
I_{(l_1\ldots l_k)T,t}^{*(i_1\ldots i_k)}
=\int\limits_t^{*T}(t-t_k)^{l_k} \ldots \int\limits_t^{*t_{2}}
(t-t_1)^{l_1} d{\bf f}_{t_1}^{(i_1)}\ldots
d{\bf f}_{t_k}^{(i_k)}
\end{equation}

\vspace{3mm}
\noindent
is the iterated Stratonovich stochastic integral;\ $i_1,\ldots,i_k = 1,\ldots,m;$\ \
$l_1,\ldots,l_k=0, 1,\ldots$

Let us first formulate some old results.

\vspace{2mm}                          

{\bf Theorem 8}\ \cite{11a}, \cite{art-5} (also see \cite{6}-\cite{11},
\cite{11aa}-\cite{art-4}, \cite{art-6}, \cite{art-9},
\cite{arxiv-2}, \cite{arxiv-4}, \cite{arxiv-5}, \cite{arxiv-9},
\cite{arxiv-11}, \cite{arxiv-14}-\cite{arxiv-19}, \cite{Kuz-Kuz}, \cite{Mikh-1}).   
{\it Suppose that 
$\{\phi_j(x)\}_{j=0}^{\infty}$ is a complete orthonormal system of 
Legendre polynomials or trigonometric functions in the space $L_2([t, T]).$
Moreover, $\psi_1(\tau),$ $\psi_2(\tau)$ are 
continuously differentiable functions on $[t, T]$. Then, 
for the iterated Stratonovich stochastic integral

\vspace{1mm}
$$
J^{*}[\psi^{(2)}]_{T,t}={\int\limits_t^{*}}^T\psi_2(t_2)
{\int\limits_t^{*}}^{t_2}\psi_1(t_1)d{\bf f}_{t_1}^{(i_1)}
d{\bf f}_{t_2}^{(i_2)}\ \ \ (i_1, i_2=1,\ldots,m)
$$

\vspace{4mm}
\noindent
the following expansion 

\begin{equation}
\label{jes}
J^{*}[\psi^{(2)}]_{T,t}=\hbox{\vtop{\offinterlineskip\halign{
\hfil#\hfil\cr
{\rm l.i.m.}\cr
$\stackrel{}{{}_{p_1,p_2\to \infty}}$\cr
}} }\sum_{j_1=0}^{p_1}\sum_{j_2=0}^{p_2}
C_{j_2j_1}\zeta_{j_1}^{(i_1)}\zeta_{j_2}^{(i_2)}
\end{equation}

\vspace{4mm}
\noindent
that converges in the mean-square
sense 
is valid, where 

\begin{equation}
\label{tupo11}
C_{j_2 j_1}=\int\limits_t^T\psi_2(t_2)\phi_{j_2}(t_2)
\int\limits_t^{t_2}\psi_1(t_1)\phi_{j_1}(t_1)dt_1dt_2
\end{equation}

\vspace{3mm}
\noindent
and
$$
\zeta_{j}^{(i)}=
\int\limits_t^T \phi_{j}(\tau) d{\bf f}_{\tau}^{(i)}
$$ 

\vspace{4mm}
\noindent
are independent
standard Gaussian random variables for various 
$i$ or $j$.}

\vspace{2mm}

{\bf Theorem 9}\ \cite{11a}, \cite{art-5} (also see \cite{6}-\cite{11},
\cite{11aa}-\cite{art-4}, \cite{art-6}, \cite{art-9},
\cite{arxiv-2}, \cite{arxiv-4}, \cite{arxiv-5}, 
\cite{arxiv-11}, \cite{arxiv-14}-\cite{arxiv-18}, \cite{Kuz-Kuz}, \cite{Mikh-1}).   
{\it Suppose that 
$\{\phi_j(x)\}_{j=0}^{\infty}$ is a complete orthonormal system of 
Legendre polynomials or trigonometric functions in the space $L_2([t, T]).$
At the same time $\psi_2(\tau)$ is a continuously dif\-ferentiable 
nonrandom function on $[t, T]$ and $\psi_1(\tau),$ $\psi_3(\tau)$ are twice
continuously differentiable nonrandom functions on $[t, T]$. 
Then, for the 
iterated Stratonovich stochastic integral of third multiplicity

\vspace{1mm}
$$
J^{*}[\psi^{(3)}]_{T,t}={\int\limits_t^{*}}^T\psi_3(t_3)
{\int\limits_t^{*}}^{t_3}\psi_2(t_2)
{\int\limits_t^{*}}^{t_2}\psi_1(t_1)
d{\bf f}_{t_1}^{(i_1)}
d{\bf f}_{t_2}^{(i_2)}d{\bf f}_{t_3}^{(i_3)}\ \ \ (i_1, i_2, i_3=1,\ldots,m)
$$

\vspace{4mm}
\noindent
the following 
expansion 

\begin{equation}
\label{feto19000a}
J^{*}[\psi^{(3)}]_{T,t}
=\hbox{\vtop{\offinterlineskip\halign{
\hfil#\hfil\cr
{\rm l.i.m.}\cr
$\stackrel{}{{}_{p\to \infty}}$\cr
}} }
\sum\limits_{j_1, j_2, j_3=0}^{p}
C_{j_3 j_2 j_1}\zeta_{j_1}^{(i_1)}\zeta_{j_2}^{(i_2)}\zeta_{j_3}^{(i_3)}
\end{equation}

\vspace{4mm}
\noindent
that converges in the mean-square sense is valid, where

$$
C_{j_3 j_2 j_1}=\int\limits_t^T\psi_3(t_3)\phi_{j_3}(t_3)
\int\limits_t^{t_3}\psi_2(t_2)\phi_{j_2}(t_2)
\int\limits_t^{t_2}\psi_1(t_1)\phi_{j_1}(t_1)dt_1dt_2dt_3
$$

\vspace{4mm}
\noindent
and
$$
\zeta_{j}^{(i)}=
\int\limits_t^T \phi_{j}(\tau) d{\bf f}_{\tau}^{(i)}
$$ 

\vspace{4mm}
\noindent
are independent standard Gaussian random variables for various 
$i$ or $j$.}

\vspace{2mm}

{\bf Theorem 10}\ \cite{11a}, \cite{art-5} (also see \cite{7}-\cite{11},
\cite{11aa}-\cite{art-4}, \cite{art-6}, \cite{art-9},
\cite{arxiv-2}, \cite{arxiv-4}, \cite{arxiv-5}, 
\cite{arxiv-11}, \cite{arxiv-14}-\cite{arxiv-18}, \cite{Kuz-Kuz}, \cite{Mikh-1}).   
{\it Suppose that
$\{\phi_j(x)\}_{j=0}^{\infty}$ is a complete orthonormal
system of Legendre polynomials or trigonometric functions
in the space $L_2([t, T])$.
Then, for the iterated 
Stratonovich stochastic integral of fourth multiplicity

\vspace{1mm}
$$
J^{*}[\psi^{(4)}]_{T,t}=
{\int\limits_t^{*}}^T
{\int\limits_t^{*}}^{t_4}
{\int\limits_t^{*}}^{t_3}
{\int\limits_t^{*}}^{t_2}
d{\bf w}_{t_1}^{(i_1)}
d{\bf w}_{t_2}^{(i_2)}d{\bf w}_{t_3}^{(i_3)}d{\bf w}_{t_4}^{(i_4)}\ \ \ 
(i_1, i_2, i_3, i_4=0, 1,\ldots,m)
$$

\vspace{4mm}
\noindent
the following 
expansion 

\begin{equation}
\label{feto1900otit}
J^{*}[\psi^{(4)}]_{T,t}=
\hbox{\vtop{\offinterlineskip\halign{
\hfil#\hfil\cr
{\rm l.i.m.}\cr
$\stackrel{}{{}_{p\to \infty}}$\cr
}} }
\sum\limits_{j_1, j_2, j_3, j_4=0}^{p}
C_{j_4 j_3 j_2 j_1}\zeta_{j_1}^{(i_1)}\zeta_{j_2}^{(i_2)}\zeta_{j_3}^{(i_3)}
\zeta_{j_4}^{(i_4)}
\end{equation}

\vspace{4mm}
\noindent
that converges in the mean-square sense is valid, where

$$
C_{j_4 j_3 j_2 j_1}=\int\limits_t^T\phi_{j_4}(t_4)\int\limits_t^{t_4}
\phi_{j_3}(t_3)
\int\limits_t^{t_3}\phi_{j_2}(t_2)\int\limits_t^{t_2}\phi_{j_1}(t_1)
dt_1dt_2dt_3dt_4
$$

\vspace{4mm}
\noindent
and
$$
\zeta_{j}^{(i)}=
\int\limits_t^T \phi_{j}(\tau) d{\bf w}_{\tau}^{(i)}
$$ 

\vspace{4mm}
\noindent
are independent standard Gaussian random variables for various 
$i$ or $j$ {\rm (}in the case when $i\ne 0${\rm ),}
${\bf w}_{\tau}^{(i)}={\bf f}_{\tau}^{(i)}$ for
$i=1,\ldots,m$ and 
${\bf w}_{\tau}^{(0)}=\tau.$}

\vspace{2mm}

Recently, a new approach to the expansion and mean-square 
approximation of iterated Stratonovich stochastic integrals has been obtained
\cite{11a} (Sect.~2.10--2.16), \cite{arxiv-5} (Sect.~13--19), 
\cite{arxiv-10} (Sect.~5--11), \cite{arxiv-11} (Sect.~7--13).
Let us formulate four theorems that were obtained using this approach.

\vspace{2mm}

{\bf Theorem 11}\ \cite{11a}, \cite{arxiv-5}, \cite{arxiv-10}, \cite{arxiv-11}.\
{\it Suppose 
that $\{\phi_j(x)\}_{j=0}^{\infty}$ is a complete orthonormal system of 
Legendre polynomials or trigonometric functions in the space $L_2([t, T]).$
Furthermore, let $\psi_1(\tau), \psi_2(\tau),$ $\psi_3(\tau)$ are continuously dif\-ferentiable 
nonrandom functions on $[t, T].$ 
Then, for the 
iterated Stra\-to\-no\-vich stochastic integral of third multiplicity

\vspace{-1mm}
$$
J^{*}[\psi^{(3)}]_{T,t}={\int\limits_t^{*}}^T\psi_3(t_3)
{\int\limits_t^{*}}^{t_3}\psi_2(t_2)
{\int\limits_t^{*}}^{t_2}\psi_1(t_1)
d{\bf w}_{t_1}^{(i_1)}
d{\bf w}_{t_2}^{(i_2)}d{\bf w}_{t_3}^{(i_3)}\ \ \ (i_1,i_2,i_3=0,1,\ldots,m)
$$

\vspace{3mm}
\noindent
the following 
relations

\begin{equation}
\label{fin1}
J^{*}[\psi^{(3)}]_{T,t}
=\hbox{\vtop{\offinterlineskip\halign{
\hfil#\hfil\cr
{\rm l.i.m.}\cr
$\stackrel{}{{}_{p\to \infty}}$\cr
}} }
\sum\limits_{j_1, j_2, j_3=0}^{p}
C_{j_3 j_2 j_1}\zeta_{j_1}^{(i_1)}\zeta_{j_2}^{(i_2)}\zeta_{j_3}^{(i_3)},
\end{equation}

\vspace{2mm}
\begin{equation}
\label{fin2}
{\sf M}\left\{\left(
J^{*}[\psi^{(3)}]_{T,t}-
\sum\limits_{j_1, j_2, j_3=0}^{p}
C_{j_3 j_2 j_1}\zeta_{j_1}^{(i_1)}\zeta_{j_2}^{(i_2)}\zeta_{j_3}^{(i_3)}\right)^2\right\}
\le \frac{C}{p}
\end{equation}

\vspace{5mm}
\noindent
are fulfilled, where $i_1, i_2, i_3=0,1,\ldots,m$ in {\rm (\ref{fin1})} and 
$i_1, i_2, i_3=1,\ldots,m$ in {\rm (\ref{fin2})},
constant $C$ is independent of $p,$

$$
C_{j_3 j_2 j_1}=\int\limits_t^T\psi_3(t_3)\phi_{j_3}(t_3)
\int\limits_t^{t_3}\psi_2(t_2)\phi_{j_2}(t_2)
\int\limits_t^{t_2}\psi_1(t_1)\phi_{j_1}(t_1)dt_1dt_2dt_3
$$

\vspace{3mm}
\noindent
and
$$
\zeta_{j}^{(i)}=
\int\limits_t^T \phi_{j}(\tau) d{\bf f}_{\tau}^{(i)}
$$ 

\vspace{2mm}
\noindent
are independent standard Gaussian random variables for various 
$i$ or $j$ {\rm (}in the case when $i\ne 0${\rm );} 
another notations are the same as in Theorem~{\rm 1}.}

\vspace{2mm}

{\bf Theorem 12}\ \cite{11a}, \cite{arxiv-5}, \cite{arxiv-10}, \cite{arxiv-11}.\ {\it Let
$\{\phi_j(x)\}_{j=0}^{\infty}$ be a complete orthonormal system of 
Legendre polynomials or trigonometric functions in the space $L_2([t, T]).$
Furthermore, let $\psi_1(\tau), \ldots,$ $\psi_4(\tau)$ be continuously dif\-ferentiable 
nonrandom functions on $[t, T].$ 
Then, for the 
iterated Stra\-to\-no\-vich stochastic integral of fourth multiplicity

\vspace{-1mm}
\begin{equation}
\label{fin0}
J^{*}[\psi^{(4)}]_{T,t}={\int\limits_t^{*}}^T\psi_4(t_4)
{\int\limits_t^{*}}^{t_4}\psi_3(t_3)
{\int\limits_t^{*}}^{t_3}\psi_2(t_2)
{\int\limits_t^{*}}^{t_2}\psi_1(t_1)
d{\bf w}_{t_1}^{(i_1)}
d{\bf w}_{t_2}^{(i_2)}d{\bf w}_{t_3}^{(i_3)}d{\bf w}_{t_4}^{(i_4)}
\end{equation}

\vspace{3mm}
\noindent
the following 
relations

\begin{equation}
\label{fin3}
J^{*}[\psi^{(4)}]_{T,t}
=\hbox{\vtop{\offinterlineskip\halign{
\hfil#\hfil\cr
{\rm l.i.m.}\cr
$\stackrel{}{{}_{p\to \infty}}$\cr
}} }
\sum\limits_{j_1, j_2, j_3,j_4=0}^{p}
C_{j_4j_3 j_2 j_1}\zeta_{j_1}^{(i_1)}\zeta_{j_2}^{(i_2)}\zeta_{j_3}^{(i_3)}\zeta_{j_4}^{(i_4)},
\end{equation}

\vspace{2mm}

\begin{equation}
\label{fin4}
{\sf M}\left\{\left(
J^{*}[\psi^{(4)}]_{T,t}-
\sum\limits_{j_1, j_2, j_3, j_4=0}^{p}
C_{j_4 j_3 j_2 j_1}\zeta_{j_1}^{(i_1)}\zeta_{j_2}^{(i_2)}\zeta_{j_3}^{(i_3)}
\zeta_{j_4}^{(i_4)}
\right)^2\right\}
\le \frac{C}{p^{1-\varepsilon}}
\end{equation}

\vspace{4mm}
\noindent
are fulfilled, where $i_1, \ldots , i_4=0,1,\ldots,m$ in {\rm (\ref{fin0}),} {\rm (\ref{fin3})} 
and $i_1, \ldots, i_4=1,\ldots,m$ in {\rm (\ref{fin4})},
constant $C$ does not depend on $p,$
$\varepsilon$ is an arbitrary
small positive real number 
for the case of complete orthonormal system of 
Legendre polynomials in the space $L_2([t, T])$
and $\varepsilon=0$ for the case of
complete orthonormal system of 
trigonometric functions in the space $L_2([t, T]),$

$$
C_{j_4 j_3 j_2 j_1}=
\int\limits_t^T\psi_4(t_4)\phi_{j_4}(t_4)
\int\limits_t^{t_4}\psi_3(t_3)\phi_{j_3}(t_3)
\int\limits_t^{t_3}\psi_2(t_2)\phi_{j_2}(t_2)
\int\limits_t^{t_2}\psi_1(t_1)\phi_{j_1}(t_1)dt_1dt_2dt_3dt_4;
$$

\vspace{4mm}
\noindent
another notations are the same as in Theorem~{\rm 11}.}

\vspace{2mm}

{\bf Theorem 13}\ \cite{11a}, \cite{arxiv-5}, \cite{arxiv-10}, \cite{arxiv-11}.\ 
{\it Assume 
that $\{\phi_j(x)\}_{j=0}^{\infty}$ is a complete orthonormal system of 
Legendre polynomials or trigonometric functions in the space $L_2([t, T])$
and $\psi_1(\tau), \ldots,$ $\psi_5(\tau)$ are continuously dif\-ferentiable 
nonrandom functions on $[t, T].$ 
Then, for the 
iterated Stra\-to\-no\-vich stochastic integral of fifth multiplicity

\vspace{-1mm}
\begin{equation}
\label{fin7}
J^{*}[\psi^{(5)}]_{T,t}={\int\limits_t^{*}}^T\psi_5(t_5)
\ldots
{\int\limits_t^{*}}^{t_2}\psi_1(t_1)
d{\bf w}_{t_1}^{(i_1)}
\ldots d{\bf w}_{t_5}^{(i_5)}
\end{equation}

\vspace{3mm}
\noindent
the following 
relations

\begin{equation}
\label{fin8}
J^{*}[\psi^{(5)}]_{T,t}
=\hbox{\vtop{\offinterlineskip\halign{
\hfil#\hfil\cr
{\rm l.i.m.}\cr
$\stackrel{}{{}_{p\to \infty}}$\cr
}} }
\sum\limits_{j_1,\ldots,j_5=0}^{p}
C_{j_5 \ldots j_1}\zeta_{j_1}^{(i_1)}\ldots \zeta_{j_5}^{(i_5)},
\end{equation}

\vspace{2mm}

\begin{equation}
\label{fin9}
{\sf M}\left\{\left(
J^{*}[\psi^{(5)}]_{T,t}-
\sum\limits_{j_1, \ldots, j_5=0}^{p}
C_{j_5 \ldots j_1}\zeta_{j_1}^{(i_1)}\ldots
\zeta_{j_5}^{(i_5)}
\right)^2\right\}
\le \frac{C}{p^{1-\varepsilon}}
\end{equation}

\vspace{5mm}
\noindent
are fulfilled, where $i_1, \ldots , i_5=0,1,\ldots,m$ in {\rm (\ref{fin7}),} {\rm (\ref{fin8})} 
and $i_1, \ldots, i_5=1,\ldots,m$ in {\rm (\ref{fin9})},
constant $C$ is independent of $p,$
$\varepsilon$ is an arbitrary
small positive real number 
for the case of complete orthonormal system of 
Legendre polynomials in the space $L_2([t, T])$
and $\varepsilon=0$ for the case of
complete orthonormal system of 
trigonometric functions in the space $L_2([t, T]),$

$$
C_{j_5 \ldots j_1}=
\int\limits_t^T\psi_5(t_5)\phi_{j_5}(t_5)\ldots
\int\limits_t^{t_2}\psi_1(t_1)\phi_{j_1}(t_1)dt_1\ldots dt_5;
$$

\vspace{3mm}
\noindent
another notations are the same as in Theorem~{\rm 11, 12}.}

\vspace{2mm}

{\bf Theorem 14}\ \cite{11a}, \cite{arxiv-5}, \cite{arxiv-10}, \cite{arxiv-11}.\ 
{\it Suppose that 
$\{\phi_j(x)\}_{j=0}^{\infty}$ is a complete orthonormal system of 
Legendre polynomials or trigonometric functions in the space $L_2([t, T]).$
Then, for the 
iterated Stratonovich stochastic integral of sixth multiplicity

\begin{equation}
\label{after10001qu1}
J_{T,t}^{*(i_1\ldots i_6)}={\int\limits_t^{*}}^T
\ldots
{\int\limits_t^{*}}^{t_2}
d{\bf w}_{t_1}^{(i_1)}
\ldots d{\bf w}_{t_6}^{(i_6)}
\end{equation}

\vspace{3mm}
\noindent
the following 
expansion 

\vspace{-1mm}
$$
J_{T,t}^{*(i_1\ldots i_6)}
=\hbox{\vtop{\offinterlineskip\halign{
\hfil#\hfil\cr
{\rm l.i.m.}\cr
$\stackrel{}{{}_{p\to \infty}}$\cr
}} }
\sum\limits_{j_1, \ldots, j_6=0}^{p}
C_{j_6 \ldots j_1}\zeta_{j_1}^{(i_1)}\ldots
\zeta_{j_6}^{(i_6)}
$$

\vspace{4mm}
\noindent
that converges in the mean-square sense is valid, where
$i_1, \ldots, i_6=0, 1,\ldots,m,$

$$
C_{j_6 \ldots j_1}=
\int\limits_t^T\phi_{j_6}(t_6)\ldots
\int\limits_t^{t_2}\phi_{j_1}(t_1)dt_1\ldots dt_6;
$$

\vspace{3mm}
\noindent
another notations are the same as in Theorems~{\rm 11--13}.}

\vspace{2mm}

The results of Theorems~8--14 were developed in 
\cite{11a} (Chapter~2).
In particular, analogues of Theorem~14 for iterated Stratonovich stochastic
integrals of multiplicities 7 and 8 were obtained in \cite{11a} (Sect.~2.36, 2.37).
In addition, the variants of Theorems 8--14
were obtained
for the case when $\{\phi_j(x)\}_{j=0}^{\infty}$ is an arbitrary complete orthonormal system
of functions in $L_2([t, T])$ \cite{11a} (Sect.~2.1.4, 2.23, 2.24, 2.31--2.34).

Consider the question on the exact calculation 
of the mean-square approximation errors
for the following iterated Stratonovich stochastic integrals

\begin{equation}
\label{dest1}
I_{(0)T,t}^{*(i_1)},\ \ \ 
I_{(1)T,t}^{*(i_1)},\ \ \ 
I_{(00)T,t}^{*(i_1i_2)},\ \ \ 
I_{(000)T,t}^{*(i_1i_2i_3)},\ \ \ 
i_1, i_2, i_3=1,\ldots,m.
\end{equation}

\vspace{4mm}

We assume that the stochastic integrals (\ref{dest1})
are approximated using Theorems 1, 8, 9 and the Legendre
polynomial system. Since
$I_{(0)T,t}^{(i_1)}=I_{(0)T,t}^{*(i_1)},$
$I_{(1)T,t}^{(i_1)}=I_{(1)T,t}^{*(i_1)}$\ w.~p.~1,
we can use the following formulas (see (\ref{a1}) for 
the case of Legendre polynomilas)

\vspace{1mm}
$$
I_{(0)T,t}^{(i_1)}=\sqrt{T-t}\zeta_0^{(i_1)},
$$

\vspace{2mm}
$$
I_{(1)T,t}^{(i_1)}=-\frac{(T-t)^{3/2}}{2}\left(\zeta_0^{(i_1)}+
\frac{1}{\sqrt{3}}\zeta_1^{(i_1)}\right)
$$

\vspace{4mm}
\noindent
to approximate the stochastic integrals 
$I_{(0)T,t}^{*(i_1)},$
$I_{(1)T,t}^{*(i_1)}.$ In this case, we will have zero
mean-square approximation errors.

To approximate the iterated Stratonovich stochastic integral 
$I_{(00)T,t}^{*(i_1i_2)}$ 
we can use the formula (see (\ref{4004}))

\vspace{-4mm}
\begin{equation}
\label{dest2}
I_{(00)T,t}^{*(i_1 i_2)p}=
\frac{T-t}{2}\left(\zeta_0^{(i_1)}\zeta_0^{(i_2)}+\sum_{i=1}^{p}
\frac{1}{\sqrt{4i^2-1}}\left(
\zeta_{i-1}^{(i_1)}\zeta_{i}^{(i_2)}-
\zeta_i^{(i_1)}\zeta_{i-1}^{(i_2)}\right)\right).
\end{equation}

\vspace{5mm}

The mean-square approximation error for (\ref{dest2})
will be determined by the formula (\ref{4007}) $(i_1\ne i_2)$.
For the case $i_1=i_2$ we can use the well-known equality \cite{Zapad-1}

\vspace{1mm}
$$
I_{(00)T,t}^{*(i_1 i_1)}
=\frac{T-t}{2}\left(\zeta_0^{(i_1)}\right)^2\ \ \ \hbox{w.~p.~1.}
$$

\vspace{4mm}

Consider now the iterated Stratonovich stochastic integral
$I_{(000)T,t}^{*(i_1i_2i_3)}$ of multiplicity 3
$(i_1, i_2, i_3=$ $1,\ldots,m)$.
For the case of pairwise different 
$i_1, i_2, i_3$ we can use the formula (\ref{dest900}).
In the case $i_1=i_2=i_3,$ to approximate the stochastic integral
$I_{(000)T,t}^{*(i_1i_1i_1)},$ we use the well-known equality \cite{Zapad-1}

\vspace{1mm}
$$
I_{(000)T,t}^{*(i_1 i_1 i_1)}
=\frac{(T-t)^{3/2}}{6}
\left(\zeta_0^{(i_1)}\right)^3\ \ \ \hbox{w.~p.~1}.
$$
                                                
\vspace{4mm}

Thus, it remains to consider the following three cases
\begin{equation}
\label{dest5}
i_1=i_2\ne i_3,
\end{equation}
\begin{equation}
\label{dest6}
i_1\ne i_2=i_3,
\end{equation}
\begin{equation}
\label{dest7}
i_1=i_3\ne i_2.
\end{equation}

\vspace{2mm}

Taking into account the standard relations between 
Ito and Stratonovich stochastic integrals \cite{Zapad-1} and
Theorem 1 (the case $k=3$) together with Theorem 9, we obtain

\vspace{1mm}
$$
{\sf M}\left\{\left(I_{(000)T,t}^{*(i_1i_2i_3)}-
I_{(000)T,t}^{*(i_1i_2i_3)p}\right)^2\right\}=
$$

\vspace{3mm}
$$
={\sf M}\left\{\left(
I_{(000)T,t}^{(i_1i_2i_3)}+
\frac{1}{2}{\bf 1}_{\{i_1=i_2\}}
\int\limits_t^T 
\int\limits_t^{\tau}dsd{\bf f}_{\tau}^{(i_3)}
+\frac{1}{2}{\bf 1}_{\{i_2=i_3\}}
\int\limits_t^T \hspace{-1mm}\int\limits_t^{\tau}d{\bf f}_{s}^{(i_1)}d\tau-
I_{(000)T,t}^{*(i_1i_2i_3)p}\right)^{2}\right\}
=
$$

\vspace{3mm}
$$
={\sf M}\left\{\left(
I_{(000)T,t}^{(i_1i_2i_3)}-I_{(000)T,t}^{(i_1i_2i_3)p}+
I_{(000)T,t}^{(i_1i_2i_3)p}+
{\bf 1}_{\{i_1=i_2\}}
\frac{1}{2}\int\limits_t^{\stackrel{~}{T}}
\int\limits_{t}^{\stackrel{~}{\tau}}
dsd{\bf f}_{\tau}^{(i_3)}
\right.\right.+
$$

\vspace{3mm}
\begin{equation}
\label{tango3}
\left.\left.
+{\bf 1}_{\{i_2=i_3\}}
\frac{1}{2}\int\limits_t^T\int\limits_t^{\tau}d{\bf f}_{s}^{(i_1)}d\tau-
I_{(000)T,t}^{*(i_1i_2i_3)p}\right)^2\right\},
\end{equation}

\vspace{5mm}
\noindent
where $I_{(000)T,t}^{(i_1i_2i_3)}$ and $I_{(000)T,t}^{(i_1i_2i_3)p}$
are defined by the relations (\ref{k1001}), (\ref{sad001}). Moreover, 
$I_{(000)T,t}^{*(i_1i_2i_3)p}$ has the form (see Theorem 9)
\begin{equation}
\label{tango2}
I_{(000)T,t}^{*(i_1i_2i_3)p}=
\sum_{j_1,j_2,j_3=0}^{p}
C_{j_3j_2j_1}
\zeta_{j_1}^{(i_1)}\zeta_{j_2}^{(i_2)}\zeta_{j_3}^{(i_3)}.
\end{equation}

\vspace{5mm}

Substituting (\ref{sad001}) and (\ref{tango2}) into (\ref{tango3}) yields

\vspace{1mm}
$$
{\sf M}\left\{\left(I_{(000)T,t}^{*(i_1i_2i_3)}-
I_{(000)T,t}^{*(i_1i_2i_3)p}\right)^2\right\}=
$$

\vspace{3mm}
$$
={\sf M}\left\{\left(I_{(000)T,t}^{(i_1i_2i_3)}-
I_{(000)T,t}^{(i_1i_2i_3)p}+
{\bf 1}_{\{i_1=i_2\}}
\left(\frac{1}{2}\int\limits_t^T
\int\limits_t^{\tau}dsd{\bf f}_{\tau}^{(i_3)}-
\sum_{j_1,j_3=0}^{q}
C_{j_3j_1j_1}
\zeta_{j_3}^{(i_3)}\right)+\right.\right.
$$

\vspace{3mm}
\begin{equation}
\label{dest10}
\left.\left.+{\bf 1}_{\{i_2=i_3\}}\left(
\frac{1}{2}\int\limits_t^T\int\limits_t^{\tau}d{\bf f}_{s}^{(i_1)}d\tau-
\sum_{j_1,j_3=0}^{p}
C_{j_3j_3j_1}
\zeta_{j_1}^{(i_1)}\right)
-{\bf 1}_{\{i_1=i_3\}}
\sum_{j_1,j_2=0}^{p}
C_{j_1j_2j_1}
\zeta_{j_2}^{(i_2)}\right)^{2}\right\}.
\end{equation}

\vspace{5mm}

Consider the case (\ref{dest5}). From (\ref{dest10}) we obtain

\vspace{1mm}
$$
{\sf M}\left\{\left(I_{(000)T,t}^{*(i_1i_2i_3)}-
I_{(000)T,t}^{*(i_1i_2i_3)p}\right)^2\right\}=
$$

\vspace{3mm}
\begin{equation}
\label{dest11}
={\sf M}\left\{\left(I_{(000)T,t}^{(i_1i_2i_3)}-
I_{(000)T,t}^{(i_1i_2i_3)p}+
\frac{1}{2}\int\limits_t^T
\int\limits_t^{\tau}dsd{\bf f}_{\tau}^{(i_3)}-
\sum_{j_1,j_3=0}^{p}
C_{j_3j_1j_1}
\zeta_{j_3}^{(i_3)}\right)^2\right\}.
\end{equation}

\vspace{5mm}

According to the formulas (\ref{yeee2}), (\ref{ttt2}), the quantity

$$
I_{(000)T,t}^{(i_1i_2i_3)}-
I_{(000)T,t}^{(i_1i_2i_3)p}
$$

\vspace{4mm}
\noindent
includes only iterated Ito stochastic integrals
of multiplicity 3. At the same time, the quantity

\vspace{1mm}
$$
\frac{1}{2}\int\limits_t^T
\int\limits_t^{\tau}dsd{\bf f}_{\tau}^{(i_3)}-
\sum_{j_1,j_3=0}^{p}
C_{j_3j_1j_1}
\zeta_{j_3}^{(i_3)}
$$

\vspace{4mm}
\noindent
contains only iterated Ito stochastic integrals
of multiplicity 1. This means that from (\ref{dest11}) we get

\vspace{1mm}
$$
{\sf M}\left\{\left(I_{(000)T,t}^{*(i_1i_2i_3)}-
I_{(000)T,t}^{*(i_1i_2i_3)p}\right)^2\right\}
={\sf M}\left\{\left(I_{(000)T,t}^{(i_1i_2i_3)}-
I_{(000)T,t}^{(i_1i_2i_3)p}\right)^2\right\}+
$$

\vspace{3mm}
\begin{equation}
\label{dest12}
+{\sf M}\left\{\left(\frac{1}{2}\int\limits_t^T
(\tau-t)d{\bf f}_{\tau}^{(i_3)}-
\sum_{j_1,j_3=0}^{p}
C_{j_3j_1j_1}
\zeta_{j_3}^{(i_3)}\right)^2\right\}.
\end{equation}

\vspace{5mm}

The relation (\ref{qq1}) implies that

\vspace{1mm}
$$
{\sf M}\left\{\left(I_{(000)T,t}^{(i_1i_2i_3)}-
I_{(000)T,t}^{(i_1i_2i_3)p}\right)^2\right\}=\frac{(T-t)^3}{6}
-
$$

\vspace{2mm}
\begin{equation}
\label{dest14}
-\sum_{j_1,j_2,j_3=0}^p C_{j_3j_2j_1}^2-
\sum_{j_1,j_2,j_3=0}^p C_{j_3j_1j_2}C_{j_3j_2j_1},
\end{equation}

\vspace{5mm}
\noindent
where $i_1=i_2\ne i_3.$

We have

\vspace{1mm}
$$
{\sf M}\left\{\left(\frac{1}{2}\int\limits_t^T
(\tau-t)d{\bf f}_{\tau}^{(i_3)}-
\sum_{j_1,j_3=0}^{p}
C_{j_3j_1j_1}
\zeta_{j_3}^{(i_3)}\right)^2\right\}=
\frac{1}{4}\int\limits_t^T
(\tau-t)^2 d\tau-
$$

\vspace{2mm}
\begin{equation}
\label{dest15}
-\sum_{j_1,j_3=0}^{p}
C_{j_3j_1j_1}\int\limits_t^T
(\tau-t)\phi_{j_3}(\tau)d\tau+
\sum_{j_3=0}^{p}\left(\sum_{j_1=0}^{p}
C_{j_3j_1j_1}\right)^2,
\end{equation}

\vspace{5mm}
\noindent
where $\{\phi_j(x)\}_{j=0}^{\infty}$ is a complete orthonormal system of 
Legendre polynomials in the space $L_2([t, T]).$

Using the orthogonality property of Legendre polynomials, we obtain

\vspace{1mm}
\begin{equation}
\label{dest16}
\int\limits_t^T
(\tau-t)\phi_{j_3}(\tau)d\tau=\frac{(T-t)^{3/2}}{2}
\left\{
\begin{matrix}
1,\ & j_3=0\cr\cr
1/\sqrt{3},\ & j_3=1\cr\cr
0,\ & j_3\ge 2
\end{matrix}
.\right.
\end{equation}

\vspace{5mm}

Combining (\ref{dest12})--(\ref{dest16}), we get

\vspace{1mm}
$$
{\sf M}\left\{\left(I_{(000)T,t}^{*(i_1i_2i_3)}-
I_{(000)T,t}^{*(i_1i_2i_3)p}\right)^2\right\}=
\frac{(T-t)^3}{4}
-\sum_{j_1,j_2,j_3=0}^p C_{j_3j_2j_1}^2-
\sum_{j_1,j_2,j_3=0}^p C_{j_3j_1j_2}C_{j_3j_2j_1}-
$$

\vspace{2mm}
\begin{equation}
\label{dest80}
-\frac{(T-t)^{3/2}}{2}
\sum_{j_1=0}^{p}
\left(C_{0j_1j_1}+\frac{1}{\sqrt{3}}C_{1j_1j_1}\right)+
\sum_{j_3=0}^{p}\left(\sum_{j_1=0}^{p}
C_{j_3j_1j_1}\right)^2,
\end{equation}

\vspace{5mm}
\noindent
where $i_1=i_2\ne i_3.$

Consider the case (\ref{dest6}). From (\ref{dest10}) we obtain

\vspace{1mm}
$$
{\sf M}\left\{\left(I_{(000)T,t}^{*(i_1i_2i_3)}-
I_{(000)T,t}^{*(i_1i_2i_3)p}\right)^2\right\}=
$$

\vspace{3mm}
$$
={\sf M}\left\{\left(I_{(000)T,t}^{(i_1i_2i_3)}-
I_{(000)T,t}^{(i_1i_2i_3)p}+
\frac{1}{2}\int\limits_t^T
\int\limits_t^{\tau}d{\bf f}_{s}^{(i_1)}d\tau-
\sum_{j_1,j_3=0}^{p}
C_{j_3j_3j_1}
\zeta_{j_1}^{(i_1)}\right)^2\right\}=
$$

\vspace{3mm}
$$
={\sf M}\left\{\left(I_{(000)T,t}^{(i_1i_2i_3)}-
I_{(000)T,t}^{(i_1i_2i_3)p}+
\frac{1}{2}\int\limits_t^T
(T-s)d{\bf f}_{s}^{(i_1)}-
\sum_{j_1,j_3=0}^{p}
C_{j_3j_3j_1}
\zeta_{j_1}^{(i_1)}\right)^2\right\}=
$$

\vspace{3mm}
$$
={\sf M}\left\{\left(I_{(000)T,t}^{(i_1i_2i_3)}-
I_{(000)T,t}^{(i_1i_2i_3)p}\right)^2\right\}+
$$

\vspace{3mm}
$$
+
{\sf M}\left\{\left(\frac{1}{2}\int\limits_t^T
(T-s)d{\bf f}_{s}^{(i_1)}-
\sum_{j_1,j_3=0}^{p}
C_{j_3j_3j_1}
\zeta_{j_1}^{(i_1)}\right)^2\right\}=
$$

\vspace{3mm}
$$
={\sf M}\left\{\left(I_{(000)T,t}^{(i_1i_2i_3)}-
I_{(000)T,t}^{(i_1i_2i_3)p}\right)^2\right\}+
$$

\vspace{3mm}
$$
+\frac{1}{4}\int\limits_t^T
(T-s)^2 ds-
\sum_{j_1,j_3=0}^{p}
C_{j_3j_3j_1}\int\limits_t^T
(T-s)\phi_{j_1}(s)ds+
$$

\vspace{3mm}
\begin{equation}
\label{dest27}
+\sum_{j_1=0}^{p}\left(\sum_{j_3=0}^{p}
C_{j_3j_3j_1}\right)^2,
\end{equation}

\vspace{5mm}
\noindent
where $\{\phi_j(x)\}_{j=0}^{\infty}$ is a complete orthonormal system of 
Legendre polynomials in the space $L_2([t, T]).$

The relation (\ref{dest1000}) implies that

\vspace{1mm}
$$
{\sf M}\left\{\left(I_{(000)T,t}^{(i_1i_2i_3)}-
I_{(000)T,t}^{(i_1i_2i_3)p}\right)^2\right\}=\frac{(T-t)^3}{6}
-
$$

\vspace{2mm}
\begin{equation}
\label{dest31}
-\sum_{j_1,j_2,j_3=0}^p C_{j_3j_2j_1}^2-
\sum_{j_1,j_2,j_3=0}^p C_{j_2j_3j_1}C_{j_3j_2j_1},
\end{equation}

\vspace{4mm}
\noindent
where $i_1\ne i_2=i_3.$

Moreover,

\vspace{-2mm}
\begin{equation}
\label{dest32}
\int\limits_t^T
(T-s)\phi_{j_1}(s)ds=\frac{(T-t)^{3/2}}{2}
\left\{
\begin{matrix}
1,\ & j_1=0\cr\cr
-1/\sqrt{3},\ & j_1=1\cr\cr
0,\ & j_1\ge 2
\end{matrix}
.\right.
\end{equation}

\vspace{5mm}

Combining (\ref{dest27})--(\ref{dest32}), we get

$$
{\sf M}\left\{\left(I_{(000)T,t}^{*(i_1i_2i_3)}-
I_{(000)T,t}^{*(i_1i_2i_3)p}\right)^2\right\}=
\frac{(T-t)^3}{4}
-\sum_{j_1,j_2,j_3=0}^p C_{j_3j_2j_1}^2-
\sum_{j_1,j_2,j_3=0}^p C_{j_2j_3j_1}C_{j_3j_2j_1}-
$$

\vspace{1mm}
\begin{equation}
\label{dest70}
-\frac{(T-t)^{3/2}}{2}
\sum_{j_3=0}^{p}
\left(C_{j_3j_3 0}-\frac{1}{\sqrt{3}}C_{j_3j_3 1}\right)+
\sum_{j_1=0}^{p}\left(\sum_{j_3=0}^{p}
C_{j_3j_3j_1}\right)^2,
\end{equation}

\vspace{4mm}
\noindent
where $i_1\ne i_2=i_3.$

Consider the case (\ref{dest7}). From (\ref{dest10}) we obtain

$$
{\sf M}\left\{\left(I_{(000)T,t}^{*(i_1i_2i_3)}-
I_{(000)T,t}^{*(i_1i_2i_3)p}\right)^2\right\}=
$$

\vspace{1mm}
$$
={\sf M}\left\{\left(I_{(000)T,t}^{(i_1i_2i_3)}-
I_{(000)T,t}^{(i_1i_2i_3)p}-
\sum_{j_1,j_2=0}^{q}
C_{j_1j_2j_1}
\zeta_{j_2}^{(i_2)}\right)^{2}\right\}=
$$

\vspace{1mm}
$$
={\sf M}\left\{\left(I_{(000)T,t}^{(i_1i_2i_3)}-
I_{(000)T,t}^{(i_1i_2i_3)p}\right)^{2}\right\}+
{\sf M}\left\{\left(\sum_{j_1,j_2=0}^{p}
C_{j_1j_2j_1}
\zeta_{j_2}^{(i_2)}\right)^{2}\right\}=
$$

\vspace{1mm}
\begin{equation}
\label{dest49}
={\sf M}\left\{\left(I_{(000)T,t}^{(i_1i_2i_3)}-
I_{(000)T,t}^{(i_1i_2i_3)p}\right)^{2}\right\}+
\sum_{j_2=0}^{p}
\left(\sum_{j_1=0}^{p}C_{j_1j_2j_1}\right)^2.
\end{equation}

\vspace{4mm}

The relation (\ref{dest1001}) implies that

\vspace{1mm}
$$
{\sf M}\left\{\left(I_{(000)T,t}^{(i_1i_2i_3)}-
I_{(000)T,t}^{(i_1i_2i_3)p}\right)^2\right\}=\frac{(T-t)^3}{6}
-
$$

\vspace{1mm}
\begin{equation}
\label{dest50}
-\sum_{j_1,j_2,j_3=0}^p C_{j_3j_2j_1}^2-
\sum_{j_1,j_2,j_3=0}^p C_{j_3j_2j_1}C_{j_1j_2j_3},
\end{equation}

\vspace{4mm}
\noindent
where $i_1=i_3\ne i_2.$

Combining (\ref{dest49}) and (\ref{dest50}), we obtain

$$
{\sf M}\left\{\left(I_{(000)T,t}^{*(i_1i_2i_3)}-
I_{(000)T,t}^{*(i_1i_2i_3)p}\right)^2\right\}=
\frac{(T-t)^3}{6}-\sum_{j_1,j_2,j_3=0}^p C_{j_3j_2j_1}^2-
\sum_{j_1,j_2,j_3=0}^p C_{j_3j_2j_1}C_{j_1j_2j_3}+
$$

\vspace{1mm}
\begin{equation}
\label{dest60}
+\sum_{j_2=0}^{p}
\left(\sum_{j_1=0}^{p}C_{j_1j_2j_1}\right)^2,
\end{equation}

\vspace{4mm}
\noindent
where $i_1=i_3\ne i_2.$

Thus, the exact calculaton of the mean-square approximation error
for the iterated Stratonovich stochastic integral 
$I_{(000)T,t}^{*(i_1i_2i_3)}$ $(i_1,i_2,i_3=1,\ldots,m)$
is given by the formulas (\ref{dest900}),
(\ref{dest80}), (\ref{dest70}), and (\ref{dest60}).

\end{document}